\documentclass[a4paper,11pt]{article}
\usepackage{a4wide}

\usepackage[T1]{fontenc}
\usepackage[utf8]{inputenc}
\usepackage{amsmath}
\usepackage{mathtools}
\usepackage{amsthm}
\usepackage{amsfonts}
\usepackage{amssymb}
\usepackage{empheq}
\usepackage{graphicx}
\usepackage{subcaption}
\usepackage{xcolor}
\usepackage{tikz}
\usepackage{hyperref}

\graphicspath{{Figures/}}

\begin{document}

\title{Isogeometric fluid-structure interaction\\ using a mixed continuous/discontinuous Galerkin scheme}
\author{R\'egis Duvigneau}
\date{Universit\'e C\^ote d'Azur, Inria, CNRS, LJAD}

\maketitle

\begin{abstract}
A mixed continuous / discontinuous Galerkin scheme is introduced for the simulation of fluid-structure interaction problems in an isogeometric analysis framework. The properties of Non-Uniform Rational B-Spline basis functions are leveraged to enable an exact transfer of the structural displacement to the fluid domain, while using different discretizations and refinements on the two sides of the coupling interface. The proposed approach is applied to the simulation of a compressible flow around an elastic wing membrane and to a classical fluid-structure benchmark involving the flow around a cylinder equipped with a hyper-elastic bar. For both cases, the results obtained are compared to those found in the literature to assess the accuracy of the proposed method.
\end{abstract}

\section{Introduction}

The accurate simulation of fluid-structure interaction problems is of upmost importance in modern engineering, in the obvious perspective of efficiency improvement, but also in terms of safety and reliability of systems. This is especially true when phenomena governed by highly nonlinear couplings occur, which cannot be understood by simulations of the fluid and the structure considered separately. Such coupled simulations can be particularly difficult to achieve, as the two disciplines may be governed by different spatial and temporal scales and are intrinsically different in their mathematical nature. The strongly coupled or monolithic approach~\cite{Blom_98,Michler_etal04} has proven its relevance in these conditions, but its implementation in an industrial context remains tedious due to its lack of flexibility. Thus, the weakly coupled or partitioned approach~\cite{Degroote_13,Felippa_eta_01} is still widely used, as it enables a highly modular framework based on different fluid and structural solvers. In this context, the numerical processing performed at the coupling interface is particularly important~\cite{Farhat_etal_98,Piperno_Farhat_01}, since different descriptions of the geometry and solution fields coexist in the fluid and structure subdomains. Therefore, many research works have focused on the effective and robust treatment of the coupling interface, giving rise to several approaches to exchange information between the subdomains, such as the nearest neighbour method or the mortar approach (see for instance~\cite{Apostolatos_19,DeBoer_07} for a detailed discussion). It should be emphasized that the geometrical aspects play a critical role, as it is necessary to control the regularity and accuracy of the structural displacement field during its transfer to the fluid domain. If this is not the case, spurious phenomena may occur in the flow, such as change in separation or reattachment location, that will in return impact the structure evolution.

Recently, the isogeometric analysis, developed by T. Hughes and co-authors, has enabled a significant qualitative leap in the treatment of geometry within a Finite-Element (FE) framework~\cite{Cottrell_etal_09}. This approach relies on Non-Uniform Rational B-Spline (NURBS) representations to define the geometry, allowing an exact description of geometries generated by Computer-Aided Design (CAD) tools. In many cases, the high regularity of the basis employed has been shown to improve solution accuracy~\cite{Bazilevs2006,Cottrell_etal_07}. Naturally, the isogeometric approach has been extended to fluid-structure interaction problems, in particular with a mortar approach~\cite{Apostolatos_19,Bazilevs_12}, using an exact matching method~\cite{Nordanger_16} or an immersed framework~\cite{Kamensky_15}.

In a previous work, a modification of the isogeometric approach has been proposed for use in a Discontinuous Galerkin (DG) framework, which is better suited to the simulation of compressible flows~\cite{Duvigneau_18}, possibly including discontinuities. This approach is based on a transformation of the NURBS computational domain to allow the use of discontinuous solutions without altering the geometry. It has been shown that the approach also facilitates the use of adaptive mesh refinement~\cite{Duvigneau_20} and can be easily extended to an Arbitrary Lagrangian Eulerian (ALE) formulation~\cite{Pezzano_Duvigneau_21}. These developments have motivated the present work, which aims at proposing an extension to fluid-structure interaction problems, by investigating the use of a mixed continuous / discontinuous Galerkin formulation.

The article is organized as follows: a first section describes the isogeometric context, in particular the NURBS basis and its important properties. In a second section, the isogeometric discretization of the structural solver is detailed, in the case of a non-linear elastic membrane and a hyper-elastic bar, based on a continuous Galerkin formulation. Then, the isogeometric discretization of the compressible Navier-Stokes equations is examined, using a discontinuous Galerkin formulation. Finally, the coupling algorithm is proposed, with a focus on the geometry treatment and the tranfer of information at the coupling interface. In section 6 and 7, two applications are presented, corresponding to the flow around a membrane wing and a classical fluid-structure interaction benchmark proposed by Turek and Hron~\cite{Turek_Hron_2006}.


\section{Isogeometric analysis}


\subsection{NURBS Basis}

\begin{figure}[!ht]
    \begin{center}
    \begin{tikzpicture}[scale = 1.5]
    \draw[->] (-1,0) -- (6,0) node[right] {$\xi$};
    \draw[dashed] (0,0) -- (0,1) ;
    \node[below] at (0,0) {$\xi_{1,2,3}$} ;
    \draw[dashed]  (1,0) -- (1,1) ;
    \node[below] at (1,0) {$\xi_{4}$} ;
    \draw[dashed]  (2,0) -- (2,1) ;
    \node[below] at (2,0) {$\xi_{5}$} ;
    \draw[dashed]  (3,0) -- (3,1) ;
    \node[below] at (3,0) {$\xi_{6}$} ;
    \draw[dashed]  (4,0) -- (4,1) ;
    \node[below] at (4,0) {$\xi_{7}$} ;
    \draw[dashed]  (5,0) -- (5,1) ;
    \node[below] at (5,0) {$\xi_{8,9,10}$} ;

    \draw[domain=0:1,smooth,variable=\x,orange,thick] plot ({\x},{(1-\x)*(1-\x)});
    \node[orange] at (0,1.2) {$N^2_1$} ;

    \draw[domain=0:1,smooth,variable=\x,black,thick] plot ({\x},{2*\x-3/2*\x*\x});
    \draw[domain=1:2,smooth,variable=\x,black,thick] plot ({\x},{0.5*(2-\x)*(2-\x)});
    \node[black] at (0.5,1) {$N^2_2$} ;

    \draw[domain=0:1,smooth,variable=\x,blue,thick] plot ({\x},{0.5*\x*\x});
    \draw[domain=1:2,smooth,variable=\x,blue,thick] plot ({\x},{-0.5+\x-(\x-1)*(\x-1)});
    \draw[domain=2:3,smooth,variable=\x,blue,thick] plot ({\x},{5/2-\x+0.5*(\x-2)*(\x-2)});
    \node[blue] at (1.5,1) {$N^2_3$} ;

    \draw[domain=1:2,smooth,variable=\x,red,thick] plot ({\x},{0.5*(\x-1)*(\x-1)});
    \draw[domain=2:3,smooth,variable=\x,red,thick] plot ({\x},{-0.5+(\x-1)-(\x-2)*(\x-2)});
    \draw[domain=3:4,smooth,variable=\x,red,thick] plot ({\x},{5/2-(\x-1)+0.5*(\x-3)*(\x-3)});
    \node[red] at (2.5,1) {$N^2_4$} ;

    \draw[domain=2:3,smooth,variable=\x,green,thick] plot ({(\x},{0.5*(\x-2)*(\x-2)});
    \draw[domain=3:4,smooth,variable=\x,green,thick] plot ({\x},{-0.5+(\x-2)-(\x-3)*(\x-3)});
    \draw[domain=4:5,smooth,variable=\x,green,thick] plot ({\x},{5/2-(\x-2)+0.5*(\x-4)*(\x-4)});
    \node[green] at (3.5,1) {$N^2_5$} ;

    \draw[domain=3:4,smooth,variable=\x,cyan,thick] plot ({\x},{0.5*(\x-3)*(\x-3)});
    \draw[domain=4:5,smooth,variable=\x,cyan,thick] plot ({\x},{-16+10*(\x-1)-3/2*(\x-1)*(\x-1)});
    \node[cyan] at (4.5,1) {$N^2_6$} ;

    \draw[domain=4:5,smooth,variable=\x,violet,thick] plot ({\x},{(\x-4)*(\x-4)});
    \node[violet] at (5,1.2) {$N^2_7$} ;

    \end{tikzpicture}
    \end{center}
\caption{Seven quadratic B-Spline basis functions}
\label{fig:bs_basis}
\end{figure}
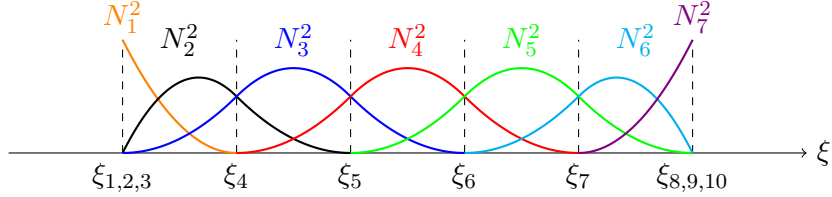

In what follows, we make the assumption that the geometry of the system is described by a set of parameterized curves. Non-Uniform Rational B-Spline (NURBS) curves are now considered as standard in CAD~\cite{DeBoor_78,Farin_89}, since they allow to represent exactly a broad class of curves commonly encountered in engineering, like conic sections. Moreover, they permit a very intuitive definition and modification of geometrical objects through the handling of \emph{control points}. Hence NURBS curves are chosen here to describe the geometry. They are defined using a so-called \emph{knot} vector  $\Xi=(\xi_1, \dots, \xi_l)\in \mathbb{R}^{l}$, which consists
of  $l$ nondecreasing real numbers. This knot vector defines a discretization of the \emph{parametric domain}
$[\xi_1 , \xi_l]$.
NURBS basis functions are derived from B-Spline functions $(N^p_{i})_{i=1, \cdots, n}$, which are polynomial and defined recursively as~\cite{DeBoor_78}:
\begin{equation}
  N^0_{i}(\xi) = \begin{cases} 1 &  \text{ if } \xi_i \leq \xi  < \xi_{i+1}\\
    0 & \text{ otherwise}
  \end{cases}
  \label{eq:nurbs0}
\end{equation}
\begin{equation}
  N^p_{i}(\xi) = \frac{\xi- \xi_i}{\xi_{i+p} - \xi_i} N^{p-1}_{i}(\xi) +
  \frac{\xi_{i+p+1}- \xi}{\xi_{i+p+1}-\xi_{i+1}} N^{p-1}_{i+1}(\xi).
  \label{eq:nurbs_eval}
\end{equation}
Note that the quotient 0/0 is assumed to be zero.
The degree of the functions $p$, the number of knots $l$ and functions $n$ are related by $l=n+p+1$~\cite{DeBoor_78}.
Open knot vectors, i.e. knot vectors with first and last knots of multiplicity $p+1$,
are usually used for representations of degree $p$ to impose interpolation and tangency conditions at both extremities~\cite{DeBoor_78}.
Therefore $\xi_1 = \ldots = \xi_{p+1}$ and $\xi_{n+1} = \ldots = \xi_{n+p+1}$.
A set of B-Spline basis functions of degree two is illustrated in Fig.~\ref{fig:bs_basis}.

Then, the NURBS basis functions $(R^p_{i})_{i=1, \cdots, n}$ are defined by associating a set of weights ($\omega_{i})_{i=1, \cdots, n}$ to the B-Spline functions according to:
\begin{equation}\label{eq:nurbsbasis}
  R^p_i(\xi) = \frac{w_i \; N^p_i(\xi)} {\sum_{j = 1}^{n} w_j \; N^p_j(\xi)} .
\end{equation}

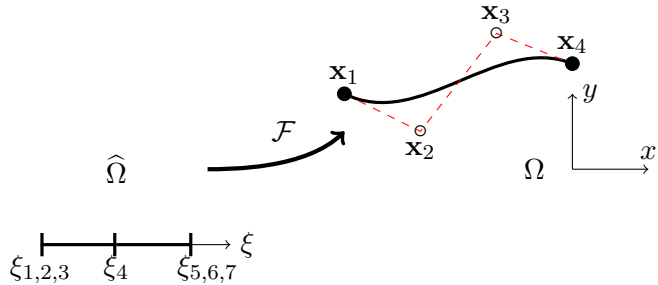
\begin{figure}[!ht]
\centering
\begin{tikzpicture}
\draw[|-|][very thick] (-2,-2) -- (-1,-2);
\draw[-|][very thick] (-1,-2) -- (0.,-2);
\draw[->] (-2,-2) -- (0.5,-2) node[right] {$\xi$};
\node[below] at (-2,-2)  {$\xi_{1,2,3}$};
\node[below] at (-1,-2)  {$\xi_{4}$};
\node[below] at (0.2,-2)  {$\xi_{5,6,7}$};
\node at (-1,-1)  {$\widehat{\Omega}$};

\draw[very thick]  (2,0) .. controls (3,-0.5) and (4,0.8) .. (5,0.4);

\draw[dashed,red] (2,0) -- (3,-0.5) -- (4,0.8) -- (5,0.4) ;
\node[circle,draw,scale=0.5,fill] at (2,0){};
\node at (3,-0.5)  {$\circ$};
\node at (4,0.8)  {$\circ$};
\node[circle,draw,scale=0.5,fill] at (5,0.4){};

\node[above] at (2,0)  {$\mathbf{x}_1$};
\node[below] at (3,-0.5)  {$\mathbf{x}_2$};
\node[above] at (4,0.8)  {$\mathbf{x}_3$};
\node[above] at (5,0.4)  {$\mathbf{x}_4$};

\node at (4.5,-1)  {$\Omega$};
\draw[ultra thick, ->]  (0.2,-1) .. controls (0.5,-1.) and (1.5,-1.) .. (2,-0.5);
\node at (1.2,-0.5)  {$\mathcal{F}$};
\draw[->] (5,-1) -- (5,0) node[right] {$y$};
\draw[->] (5,-1) -- (6,-1) node[above] {$x$};
\end{tikzpicture}
\caption{Example of quadratic NURBS curve with four control points}
\label{fig:nurbscurve}
\end{figure}

Finally, a NURBS curve of coordinates $\mathbf{x} (\xi) = (x(\xi),y(\xi))$ is obtained by positioning a set of \emph{control points} which are associated to the NURBS basis functions:
\begin{equation}\label{eq:nurbscurve1}
  \mathbf{x} (\xi) = \sum_{i=1}^{n} R^p_i(\xi) \mathbf{x}_i ,
\end{equation}
where $\mathbf{X} = (\mathbf{x}_i)_{i=1, \cdots, n}$ are the coordinates of the control points in the \emph{physical domain}. Note that the previous equation defines the mapping $\mathcal{F}$ from the parametric domain $\widehat{\Omega}$ to the physical one $\Omega$, as illustrated in Fig.~\ref{fig:nurbscurve}.

\subsection{Some properties}

\paragraph{Projection}
A first important property concerns the relationship between B-Spline and NURBS curves: it can be shown~\cite{Piegl_95} that any NURBS curve can be considered as the projection of a B-Spline curve defined in a space of higher dimension using its weights. For instance, a NURBS curve lying in the plane $(x,y)$, characterized by a set of control points  $\mathbf{X} = (\mathbf{x}_{i})_{i=1, \cdots, n} = (x_{i},y_{i})_{i=1, \cdots, n}$ and weights  $(\omega_{i})_{i=1, \cdots, n}$, can be described as the projection of a B-Spline curve lying in a 3D space and defined by the control points $(x_{i}\omega_{i},y_{i}\omega_{i},\omega_{i} )_{i=1, \cdots, n}$. This representation of NURBS curves as B-Spline ones is useful for practical manipulations because one can apply transformations directly to B-Spline curves and project the result, instead of developing procedures specific to NURBS.

\paragraph{Knot insertion}
A second property of NURBS representations is the capability to insert a new knot, and thus a new basis function, without altering the geometrical object~\cite{Piegl_95}.
This procedure can be considered as a local $h$-refinement procedure~\cite{Hughes_etal_05} and is referred as \emph{knot insertion}.

\begin{figure}[!ht]
\centering
\begin{tikzpicture}
\draw[|-|][very thick] (-2,-2) -- (-1,-2);
\draw[-|][very thick] (-1,-2) -- (-0.5,-2);
\draw[-|][very thick] (-0.5,-2) -- (0.,-2);
\draw[->] (-2,-2) -- (0.5,-2) node[right] {$\xi$};
\node[below] at (-2,-2)  {$\xi_{1,2,3}$};
\node[below] at (-1,-2)  {$\xi_{4}$};
\node[below] at (-0.5,-2)  {$\xi_{5}$};
\node[below] at (0.2,-2)  {$\xi_{6,7,8}$};
\node at (-1,-1)  {$\widehat{\Omega}$};

\draw[very thick]  (2,0) .. controls (3,-0.5) and (4,0.8) .. (5,0.4);

\draw[dashed,red] (2,0) -- (2.85,-0.4) -- (3.7,0.4) -- (4.4,0.6) -- (5,0.4) ;
\node[circle,draw,scale=0.5,fill] at (2,0){};
\node at (2.85,-0.4)  {$\circ$};
\node at (3.7,0.4)  {$\circ$};
\node at (4.4,0.6)  {$\circ$};
\node[circle,draw,scale=0.5,fill] at (5,0.4){};

\node[above] at (2,0)  {$\mathbf{x}_1$};
\node[below] at (2.85,-0.4)  {$\mathbf{x}_2$};
\node[above] at (3.7,0.4)  {$\mathbf{x}_3$};
\node[above] at (4.4,0.6)  {$\mathbf{x}_4$};
\node[above] at (5,0.4)  {$\mathbf{x}_5$};

\node at (4.5,-1)  {$\Omega$};
\draw[ultra thick, ->]  (0.2,-1) .. controls (0.5,-1.) and (1.5,-1.) .. (2,-0.5);
\node at (1.2,-0.5)  {$\mathcal{F}$};
\draw[->] (5,-1) -- (5,0) node[right] {$y$};
\draw[->] (5,-1) -- (6,-1) node[above] {$x$};
\end{tikzpicture}
\caption{Example of curve with a new inserted knot}
\label{fig:nurbsinsertion}
\end{figure}
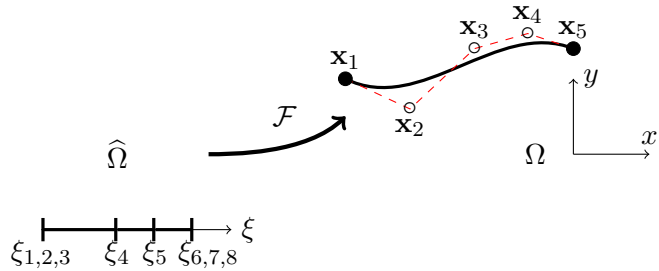

Indeed, any NURBS curve with the associated knot vector $\Xi=(\xi_1, \dots, \xi_l)\in \mathbb{R}^l$ can be identically represented using the knot vector $(\xi_1, \dots, \xi_q, \bar{\xi}, \xi_{q+1} \ldots,  \xi_l)\in \mathbb{R}^{l+1}$, that includes the additional knot $\bar{\xi}$ inserted between $\xi_q$ and $\xi_{q+1}$. The new curve is then defined by $n+1$ control points. To compute their new locations, one can use the projection property described above and simply apply the procedure to the corresponding B-Spline curve defined in the space of higher dimension. With the additional knot, the B-Spline curve can be written as~\cite{Piegl_95}:
\begin{equation}\label{eq:insertion}
  \mathbf{x} (\xi) = \sum_{i=1}^{n} N^p_{i}(\xi) \mathbf{x}_i  = \sum_{i=1}^{n+1} N^p_{i}(\xi) \mathbf{\bar{x}}_i ,
\end{equation}
with the new set of $n+1$ control points defined as:
\begin{equation}\label{eq:insertion2}
\mathbf{\bar{x}}_i = (1-\alpha_i) \mathbf{x}_{i-1} + \alpha_i \mathbf{x}_i
\quad \alpha_i = \left\{
    \begin{array}{ll}
        1 & \mbox{if } i \le q-p\\
        \frac{\bar{\xi} - \xi_i}{\xi_{i+p} - \xi_i} & \mbox{if } q-p+1 \le i \le q\\
        0 & \mbox{if } i \ge q+1
    \end{array}
\right.
\end{equation}
The figure~\ref{fig:nurbsinsertion} illustrates the knot insertion procedure applied to the curve of the previous figure.

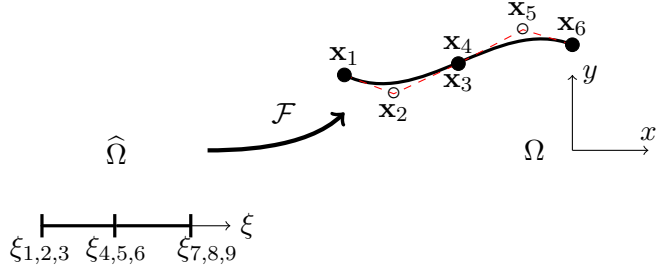
\begin{figure}[!ht]
\centering
\begin{tikzpicture}
\draw[|-|][very thick] (-2,-2) -- (-1,-2);
\draw[-|][very thick] (-1,-2) -- (0.,-2);
\draw[->] (-2,-2) -- (0.5,-2) node[right] {$\xi$};
\node[below] at (-2,-2)  {$\xi_{1,2,3}$};
\node[below] at (-1,-2)  {$\xi_{4,5,6}$};
\node[below] at (0.2,-2)  {$\xi_{7,8,9}$};
\node at (-1,-1)  {$\widehat{\Omega}$};

\draw[very thick]  (2,0) .. controls (3,-0.5) and (4,0.8) .. (5,0.4);

\draw[dashed,red] (2,0) -- (2.65,-0.25) -- (3.5,0.15) -- (4.35,0.6) -- (5,0.4) ;
\node[circle,draw,scale=0.5,fill] at (2,0){};
\node at (2.65,-0.25)  {$\circ$};
\node[circle,draw,scale=0.5,fill] at (3.5,0.15){};
\node at (4.35,0.6)  {$\circ$};
\node[circle,draw,scale=0.5,fill] at (5,0.4){};

\node[above] at (2,0)  {$\mathbf{x}_1$};
\node[below] at (2.65,-0.25)  {$\mathbf{x}_2$};
\node[below] at (3.5,0.15)  {$\mathbf{x}_3$};
\node[above] at (3.5,0.15)  {$\mathbf{x}_4$};
\node[above] at (4.35,0.6)  {$\mathbf{x}_5$};
\node[above] at (5,0.4)  {$\mathbf{x}_6$};

\node at (4.5,-1)  {$\Omega$};
\draw[ultra thick, ->]  (0.2,-1) .. controls (0.5,-1.) and (1.5,-1.) .. (2,-0.5);
\node at (1.2,-0.5)  {$\mathcal{F}$};
\draw[->] (5,-1) -- (5,0) node[right] {$y$};
\draw[->] (5,-1) -- (6,-1) node[above] {$x$};
\end{tikzpicture}
\caption{Example of curve after Bézier extraction procedure}
\label{fig:nurbsextraction}
\end{figure}

\paragraph{Bézier extraction}
It is important to underline now a particular case: if a new knot is inserted at an existing knot location, the regularity of the curve is decreased.
More generally, the curve at the knot $\xi_q$ has the regularity $C^{p-r}$, where $r$ is the multiplicity of the knot $q$~\cite{Piegl_95}.
Therefore, if one inserts $p$ knots at an existing knot location, the geometry of the curve is preserved but the curve is now divided in two independent parts. If this multiple-knot insertion is achieved for all inner knots, the curve is decomposed into a set of independent \emph{Bézier curves}, each of them composed of $p+1$ control points~\cite{Piegl_95}, as depicted in figure~\ref{fig:nurbsextraction}.

\subsection{Surfaces and volumes}

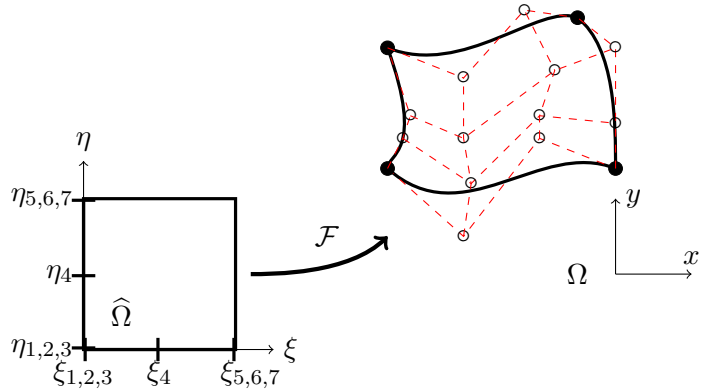
\begin{figure}[!ht]
\centering
\begin{tikzpicture}
\draw[very thick] (-2,-2) rectangle (0,0);
\draw[->] (-2,-2) -- (0.5,-2) node[right] {$\xi$};
\draw[->] (-2,-2) -- (-2,0.5) node[above] {$\eta$};
\draw[|-|][very thick] (-2,-2) -- (-1,-2);
\draw[-|][very thick] (-1,-2) -- (0.,-2);
\node[below] at (-2,-2)  {$\xi_{1,2,3}$};
\node[below] at (-1,-2)  {$\xi_{4}$};
\node[below] at (0.2,-2)  {$\xi_{5,6,7}$};
\draw[|-|][very thick] (-2,-2) -- (-2,-1);
\draw[-|][very thick] (-2,-1) -- (-2,0);
\node[left] at (-2,-2)  {$\eta_{1,2,3}$};
\node[left] at (-2,-1)  {$\eta_{4}$};
\node[left] at (-2,0)  {$\eta_{5,6,7}$};

\node at (-1.5,-1.5)  {$\widehat{\Omega}$};

\draw[very thick]  (2,0.4) .. controls (3,-0.5) and (4,0.8) .. (5,0.4);
\draw[very thick]  (2,0.4) .. controls (2,0.6) and (2.5,0.6) .. (2,2);
\draw[very thick]  (2,2) .. controls (3,1.6) and (4.,2.6) .. (4.5,2.4);
\draw[very thick]  (5,0.4) .. controls (5,1) and (5.,2.) .. (4.5,2.4);

\draw[dashed,red] (2,0.4) -- (3,-0.5) -- (4,0.8) -- (5,0.4) ;
\node[circle,draw,scale=0.5,fill] at (2,0.4){};
\node at (3,-0.5)  {$\circ$};
\node at (4,0.8)  {$\circ$};
\node[circle,draw,scale=0.5,fill] at (5,0.4){};
\draw[dashed,red] (2,0.4) -- (2.2,0.8) -- (2.3,1.1) -- (2,2);
\node at (2.2,0.8)  {$\circ$};
\node at (2.3,1.1)  {$\circ$};
\node[circle,draw,scale=0.5,fill] at (2,2){};
\draw[dashed,red]  (2,2) -- (3,1.6) -- (3.8,2.5) -- (4.5,2.4);
\node at (3,1.6)  {$\circ$};
\node at (3.8,2.5)  {$\circ$};
\node[circle,draw,scale=0.5,fill] at (4.5,2.4){};
\draw[dashed,red]  (5,0.4) -- (5,1) -- (5.,2.) -- (4.5,2.4);
\node at (5,1)  {$\circ$};
\node at (5,2)  {$\circ$};

\draw[dashed,red]  (2.2,0.8) -- (3.1,0.2) -- (4,1.1) -- (5,1);
\node at (3.1,0.2)  {$\circ$};
\node at (4,1.1)  {$\circ$};
\draw[dashed,red]  (3,-0.5) -- (3.1,0.2) -- (3.,0.8) -- (3,1.6);
\node at (3.,0.8)  {$\circ$};
\node at (4.2,1.7)  {$\circ$};
\draw[dashed,red]  (4,0.8) -- (4,1.1) -- (4.2,1.7) -- (3.8,2.5);
\draw[dashed,red]  (2.3,1.1) -- (3.,0.8) -- (4.2,1.7) -- (5,2);

\node at (4.5,-1)  {$\Omega$};
\draw[ultra thick, ->]  (0.2,-1) .. controls (0.5,-1.) and (1.5,-1.) .. (2,-0.5);
\node at (1.2,-0.5)  {$\mathcal{F}$};
\draw[->] (5,-1) -- (5,0) node[right] {$y$};
\draw[->] (5,-1) -- (6,-1) node[above] {$x$};
\end{tikzpicture}
\caption{Example of quadratic NURBS surface with $4 \times 4$ control points}
\label{fig:nurbssurf}
\end{figure}

All the previous concepts and properties can be extended to surfaces (respectively volumes) by using bivariate (respectively trivariate) tensor products. In particular, a NURBS surface of degree $p$ is defined as:
\begin{equation}\label{eq:nurbssurf}
  \mathbf{x} (\xi, \eta) = \sum_{i=1}^{n_1} \sum_{j=1}^{n_2} R^p_i(\xi) R^p_j(\eta) \mathbf{x}_{ij} ,
\end{equation}
where $(\mathbf{x}_{ij})_{i=1, \cdots, n_1 ; j=1,\cdots,n_2}$ are the coordinates of the control point indexed $ij$ in the lattice. Such a surface defined in the plane is shown in figure~\ref{fig:nurbssurf}. This tensorial construction is especially important because it will be the baseline of the CAD-consistent computational domain.

\begin{figure}[!ht]
\centering
\begin{tikzpicture}
\draw[very thick] (-2,-2) rectangle (0,0);
\draw[->] (-2,-2) -- (0.5,-2) node[right] {$\xi$};
\draw[->] (-2,-2) -- (-2,0.5) node[above] {$\eta$};
\draw[|-|][very thick] (-2,-2) -- (-1,-2);
\draw[-|][very thick] (-1,-2) -- (0.,-2);
\node[below] at (-2,-2)  {$\xi_{1,2,3}$};
\node[below] at (-1,-2)  {$\xi_{4,5,6}$};
\node[below] at (0.2,-2)  {$\xi_{7,8,9}$};
\draw[|-|][very thick] (-2,-2) -- (-2,-1);
\draw[-|][very thick] (-2,-1) -- (-2,0);
\node[left] at (-2,-2)  {$\eta_{1,2,3}$};
\node[left] at (-2,-1)  {$\eta_{4,5,6}$};
\node[left] at (-2,0)  {$\eta_{7,8,9}$};
\draw[|-|][very thick] (-1,-2) -- (-1,0);
\draw[|-|][very thick] (-2,-1) -- (0,-1);

\node at (-1.5,-1.5)  {$\widehat{\Omega}$};

\draw[very thick]  (2,0.4) .. controls (3,-0.5) and (4,0.8) .. (5,0.4);
\draw[very thick]  (2,0.4) .. controls (2.2,0.6) and (2.5,0.6) .. (2,2);
\draw[very thick]  (2,2) .. controls (3,1.6) and (4.,2.6) .. (4.5,2.4);
\draw[very thick]  (5,0.4) .. controls (5,1) and (5.,2.) .. (4.5,2.4);
\draw[very thick]  (2.2,1.2) .. controls (3,0.8) and (4,1.3) .. (4.9,1.4);
\draw[very thick]  (3.5,0.2) .. controls (3.4,1) and (3.7,1.8) .. (3.5,2.2);

 \node[circle,draw,scale=0.5,fill] at (2,0.4){};
 \node[circle,draw,scale=0.5,fill] at (5,0.4){};
 \node[circle,draw,scale=0.5,fill] at (2,2){};
 \node[circle,draw,scale=0.5,fill] at (4.5,2.4){};

 \draw[dashed,red] (2,0.4) -- (2.6,-0.15) -- (3.5,0.2) ;
 \node at (2.6,-0.15)  {$\circ$};
 \draw[dashed,red] (3.5,0.2) -- (4.4,0.55) -- (5,0.4) ;
 \node at (4.4,0.55)  {$\circ$};
 \draw[dashed,red] (2.3,0.7) -- (2.8,0.5) -- (3.45,0.6) ;
 \node at (2.3,0.7)  {$\circ$};
 \node at (2.8,0.5)  {$\circ$};
 \node at (3.45,0.6)  {$\circ$};
 \draw[dashed,red] (3.45,0.6) -- (4.3,0.8) -- (5,0.9) ;
 \node at (4.3,0.8)  {$\circ$};
 \node at (5,0.9)  {$\circ$};
 \draw[dashed,red] (2.2,1.2) -- (2.8,0.95) -- (3.48,1.08) ;
 \node at (2.8,0.95)  {$\circ$};
 \node[circle,draw,scale=0.5,fill] at (3.48,1.08){};
 \draw[dashed,red] (3.48,1.08) -- (4.3,1.3) -- (4.95,1.4) ;
 \node at (4.3,1.3)  {$\circ$};

 \draw[dashed,red] (2,0.4) -- (2.3,0.7) -- (2.2,1.2) ;
 \draw[dashed,red] (2.6,-0.15) -- (2.8,0.5) -- (2.8,0.95) ;
 \draw[dashed,red] (3.5,0.2) -- (3.45,0.6) -- (3.48,1.08) ;
 \draw[dashed,red] (4.4,0.55) -- (4.3,0.8) -- (4.3,1.3) ;
 \draw[dashed,red] (5,0.4) -- (5,0.9) -- (4.95,1.4) ;

 \draw[dashed,red] (2.15,1.6) -- (2.8,1.5) -- (3.6,1.6) ;
 \node at (2.15,1.6)  {$\circ$};
 \node at (2.8,1.5)  {$\circ$};
 \node at (3.6,1.6)  {$\circ$};
 \draw[dashed,red] (3.6,1.6) -- (4.3,1.9) -- (4.95,1.9) ;
 \node at (4.3,1.9)  {$\circ$};
 \node at (4.95,1.9)  {$\circ$};
 \draw[dashed,red] (2,2) -- (2.8,1.8) -- (3.5,2.15) ;
 \node at (2.8,1.8)  {$\circ$};
 \draw[dashed,red] (3.5,2.15) -- (4,2.4) -- (4.5,2.4) ;
 \node at (4,2.4)  {$\circ$};

 \draw[dashed,red] (2.2,1.2) -- (2.15,1.6) -- (2,2) ;
 \draw[dashed,red] (2.8,0.95) -- (2.8,1.5) -- (2.8,1.8) ;
 \draw[dashed,red] (3.48,1.08) -- (3.6,1.6) -- (3.5,2.15) ;
 \draw[dashed,red] (4.3,1.3) -- (4.3,1.9) -- (4,2.4) ;
 \draw[dashed,red] (4.95,1.4) -- (4.95,1.9) -- (4.5,2.4) ;

\node[circle,draw,scale=0.5,fill] at (2.2,1.2){};
\node[circle,draw,scale=0.5,fill] at (4.95,1.4){};
\node[circle,draw,scale=0.5,fill] at (3.5,0.2){};
\node[circle,draw,scale=0.5,fill] at (3.5,2.15){};

\node at (4.5,-1)  {$\Omega$};
\draw[ultra thick, ->]  (0.2,-1) .. controls (0.5,-1.) and (1.5,-1.) .. (2,-0.5);
\node at (1.2,-0.5)  {$\mathcal{F}$};
\draw[->] (5,-1) -- (5,0) node[right] {$y$};
\draw[->] (5,-1) -- (6,-1) node[above] {$x$};
\end{tikzpicture}
\caption{Example of quadratic NURBS surface split to rational Bézier patches.}
\label{fig:beziersurf}
\end{figure}
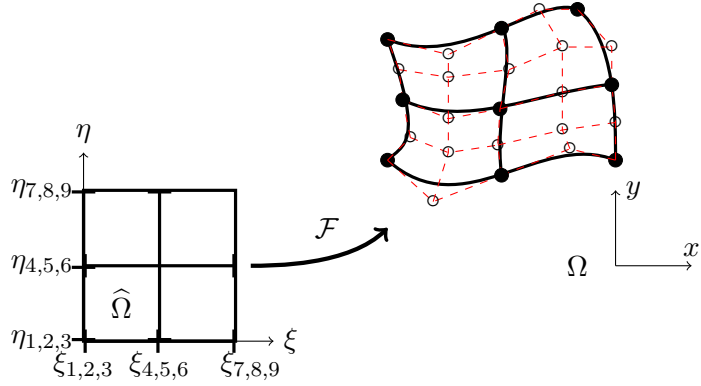

The knot insertion and the Bézier extraction procedures can be achieved for NURBS surfaces by applying the one-dimensional algorithm to each direction sequentially, thanks to the tensorial construction. Fig.~(\ref{fig:beziersurf}) depicts a NURBS surface split into four rational Bézier \emph{patches}. In the next sections, we see how these representations can be used to define isogeometric discretizations for the fluid and structure systems.

\section{Structural solver}


\subsection{Membrane wing case}

\subsubsection{Model}

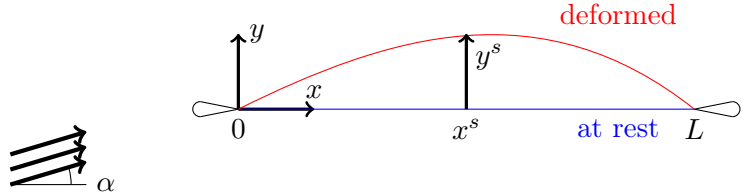
\begin{figure}[!ht]
\centering
\begin{tikzpicture}

\draw[->,very thick] (0,0) -- (0,1) node[right] {$y$};
\draw[->, very thick] (0,0) -- (1,0) node[above] {$x$};

\draw (0,0) node[below] {$0$};
\draw[-,blue] (0,0) -- (6,0) ;
\draw (6,0) node[below] {$L$};

\draw[-,red]  (0,0) .. controls (2,1) and (4,1.6) .. (6,0);

\draw (3,0) node[below] {$x^s$};
\draw (3.3,1) node[below] {$y^s$};
\draw[->, very thick] (3,0) -- (3,1);

\draw[blue] (5,0) node[below] {at rest};
\draw[red] (5,1) node[above] {deformed};

\draw[->, ultra thick] (-3,-0.6) -- (-2,-0.3);
\draw[->, ultra thick] (-3,-0.8) -- (-2,-0.5);
\draw[->, ultra thick] (-3,-1.) -- (-2,-0.7);
\draw[-] (-3,-1) -- (-2,-1) node[right] {$\alpha$};

\draw[-] (-2.2,-1) arc (0:15:0.8);

\draw[-] (0,0) -- (-0.5,0.1);
\draw[-] (0,0) -- (-0.5,-0.1);
\draw[-] (-0.5,0.1) arc (90:270:0.1);

\draw[-] (6,0) -- (6.5,0.1);
\draw[-] (6,0) -- (6.5,-0.1);
\draw[-] (6.5,-0.1) arc (-90:90:0.1);

\end{tikzpicture}
\caption{Configuration of the membrane wing.}
\label{fig:membrane}
\end{figure}

One considers as first application a one-dimensional membrane wing, held by two supports at a distance $L$ from one edge to the other, of thickness $h \ll L$ and oriented along $x$ at rest, as depicted in Fig.~(\ref{fig:membrane}). The membrane is immersed in steady free-stream flow of incidence $\alpha$. The details of the configuration, such as the geometry of the supports, can be found in~\cite{Gordnier_09}.

The structural model adopted for the membrane originates from the non-linear one proposed by Smith and Shyy in~\cite{Smith_Shyy_95}, which has been used in several numerical studies~\cite{Gordnier_09,Molki_10,Serrano_16}. Note that some other non-linear models may also be employed, such as the Neo-Hookean formulation~\cite{Frohle_persson_14}, interested readers can consult the review on membrane wing models proposed in~\cite{Tiomkin_21}.

In this context, the normal displacement $y^s$ of a membrane point initially located at $(x^s,0)$, subject to the normal pressure force $\Delta p$ exerted on both sides, is governed by~\cite{Smith_Shyy_95}:
\begin{equation}
\rho^s h \frac{\partial^2 y^s}{\partial t^2} - T \frac{\partial^2 y^s}{\partial {x^s}^2} \left [ 1 + \left (\frac{\partial y^s}{\partial x^s} \right )^2 \right ]^{-3/2} = \Delta p,
\label{eq:membrane}
\end{equation}
where $\rho^s$ is the membrane density. The tension $T$ depends on the pretension $T_0$ and the length increase $\delta$:
\begin{equation}
T = T_0 + E h \delta \qquad \delta = \int_0^L \sqrt{1+ \left ( \frac{\partial y^s}{\partial x^s} \right )^2 } \mbox{d}x^s - L ,
\label{eq:tension}
\end{equation}
with $E$ the modulus of elasticity. As explained in~\cite{Molki_10}, small deviations of the membrane make usually the non-linear term between square brackets in Eq.~(\ref{eq:membrane}) negligible. Then, the model adopted is simplified to:
\begin{equation}
\rho^s h \frac{\partial^2 y^s}{\partial t^2} - T \frac{\partial^2 y^s}{\partial {x^s}^2} = \Delta p .
\label{eq:membrane2}
\end{equation}
This model remains however non-linear due to the dependency of the tension $T(y)$ on the displacement.
The membrane is supposed to be pinned at both extremities yielding zero displacement for $x^s = 0$ and $x^s = L$. Note also that no dissipation term is included in the model.


\subsubsection{Isogeometric discretization}

The problem is solved according to the isogeometric analysis paradigm, by using NURBS basis functions. Thus, the membrane is described by a single NURBS curve of coordinates $\mathbf{x}^s (\xi^s) = (x^s(\xi^s),y^s(\xi^s))$, obtained by defining a set of control points which are associated to the NURBS basis functions:
\begin{equation}\label{eq:nurbscurve2}
    \mathbf{x}^s (\xi^s) = \sum_{i=1}^{n} R^p_i(\xi^s) \mathbf{x}^s_i
\end{equation}
where $\mathbf{X}^s = (X^s,Y^s) = (x^s_i, y^s_i)_{i=1, \cdots, n} = (\mathbf{x}_i^s)_{i=1, \cdots, n}$ are the coordinates of the control points in the physical domain, as illustrated in Fig.~(\ref{fig:nurbscurve}). According to the model chosen, the vector $Y^s$ represents the degrees of freedom for the normal displacement of the membrane along time, whereas $X^s$ is fixed. The boundary conditions (fixed extremities) correspond simply to $y^s_1 = y^s_n = 0$ thanks to the use of an open knot vector.


\subsubsection{Galerkin method}

Eq.~(\ref{eq:membrane2}) is solved using a variational formulation, which writes in the parametric domain, after integration by part:
\begin{equation}
\int_{\hat{\Gamma}^s} \varphi \; \rho^s h \, \ddot{y}^s \; |\mbox{J}| \mbox{d}\xi^s +
\int_{\hat{\Gamma}^s} T(y) \, \frac{\partial \varphi}{\partial \xi^s} \, \frac{\partial y^s}{\partial \xi^s} \; |\mbox{J}|^{-1} \mbox{d}\xi^s =
\int_{\hat{\Gamma}^s} \varphi \; \Delta p  \; |\mbox{J}| \mbox{d}\xi^s \qquad \forall \varphi \in H^1_0(\hat{\Gamma}^s),
\label{eq:variational}
\end{equation}
where $(\dot{ })$ stands for the time derivative, $|\mbox{J}| = | \frac{\partial x^s}{\partial \xi^s} |$ is the Jacobian of the isogeometric mapping $\mathcal{F}$ and $\varphi$ a test function. By adopting a Galerkin approach and choosing $\varphi = R^p_j$, one obtains:
\begin{equation}
\sum_{i=1}^n \left [ \int_{\hat{\Gamma}^s} \rho^s h R^p_j R^p_i \, |\mbox{J}| \mbox{d}\xi^s \right ] \ddot{y}^s_i  +
\sum_{i=1}^n \left [ \int_{\hat{\Gamma}^s} T(y)  \frac{\partial R^p_j}{\partial \xi^s} \; \frac{\partial R^p_i}{\partial \xi^s}  |\mbox{J}|^{-1} \mbox{d}\xi^s \right ] y^s_i =
\int_{\hat{\Gamma}^s} \Delta p  \; R^p_j \; |\mbox{J}| \mbox{d}\xi^s ,
\label{eq:variational_galerkin}
\end{equation}
where $R^p_i$ are the NURBS trial functions. By introducing the mass matrix $\mathcal{M}^s$ and the stifness matrix $\mathcal{K}^s$, the problem reads:
\begin{equation}
\mathcal{M}^s \ddot{Y}^s  + \mathcal{K}^s(Y^s) Y^s = F^s ,
\label{eq:variational_matrix}
\end{equation}
with $F^s$ the right-hand side term in Eq.~(\ref{eq:variational_galerkin}).
The computation of the matrices and the right-hand side vector is achieved using Gauss-Legendre quadratures, for each knot interval due to the piecewise definition of the NURBS basis functions. In this context, the knot intervals can therefore be considered as elements for the discretization.

\subsection{Turek's case}

\subsubsection{Model}

\begin{figure}[!ht]
    \centering
    \includegraphics[width=0.6\textwidth]{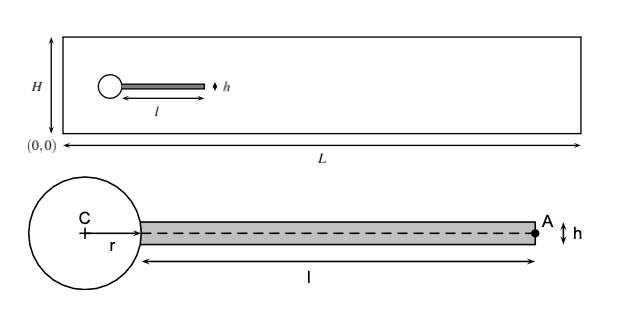}
  \caption{Configuration of Turek's Benchmark}
  \label{fig:turek_config}
\end{figure}

We consider as second application a well-known benchmark problem proposed by Turek and Hron~\cite{Turek_Hron_2006}, dealing with the flow around a cylinder with an attached flexible bar. The geometrical configuration of the problem is depicted in Fig.(\ref{fig:turek_config}).

The flexible bar is assumed to be elastic and compressible. The mechanical balance equation written in Lagrangian description writes:
\begin{equation}
\rho^s \frac{\partial^2 \mathbf{u}^s}{\partial t^2} = \nabla \cdot (J \boldsymbol\sigma \mathbf{F}^{-T}) + \rho^s \mathbf{g} ,
\label{eq:balance}
\end{equation}
where $\mathbf{u}^s$ denotes the displacement field, $\rho^s$ the material density and $\mathbf{g}$ the gravity vector. $\mathbf{F} = \mathbf{I} + \nabla \mathbf{u}^s$ is the deformation gradient tensor, $J = \det \mathbf{F}$ and $\boldsymbol\sigma$ the Cauchy stress tensor. The Saint Venant-Kirchhoff constitutive law is employed to specify the material, via the second Piola-Kirchhoff stress tensor $\mathbf{S} = J \mathbf{F}^{-1} \boldsymbol\sigma \mathbf{F}^{-T}$:
\begin{equation}
\mathbf{S} = \lambda^s (tr \mathbf{E}) \mathbf{I} + 2 \mu^s \mathbf{E} ,
\end{equation}
where $\mathbf{E} = \frac{1}{2}(\mathbf{F}^T\mathbf{F}-\mathbf{I})$ is the Green-Lagrange strain tensor. To obtain the weak form within a total Lagrangian formulation, Eq.~(\ref{eq:balance}) is integrated over the undeformed reference configuration $\Omega_0^f$, yielding after integration by parts:
\begin{equation}
\label{eq:hyperelasticity_weak}
\int_{\Omega_0^s} \rho^s \frac{\partial^2 \mathbf{u}^s}{\partial t^2} \cdot \boldsymbol\varphi \,\mbox{d} \Omega_0 +
  \int_{\Omega_0^s} (\mathbf{F}\mathbf{S}) : \nabla\boldsymbol\varphi  \, \mbox{d} \Omega_0 = \int_{\Gamma_0^s} \boldsymbol\tau^s \cdot \boldsymbol\varphi \, \mbox{d} \Gamma_0 + \int_{\Omega_0^s} \rho^s \mathbf{g} \cdot \boldsymbol\varphi  \, \mbox{d} \Omega_0 \qquad \forall \boldsymbol\varphi \in H^1(\Omega^s_0),
\end{equation}
where $\boldsymbol\varphi$ is the test function and $\boldsymbol\tau^s$ the traction exerted on the Neumann boundary $\Gamma_0^s$.


\subsubsection{Isogeometric discretization}

The undeformed geometry $\mathbf{x}_0^s$ as well as the displacement field $\mathbf{u}^s$ are both represented by NURBS surfaces in the $(x,y)$ plane, defined according to the control points $({\mathbf{x}_0^s}_{ij})_{i=1, \cdots, n_1 ; j=1,\cdots,n_2}$ and $({\mathbf{u}^s}_{ij})_{i=1, \cdots, n_1 ; j=1,\cdots,n_2}$. As a result, the deformed geometry $\mathbf{x}^s$ is also described by the same manner. Due to the simplicity of the geometry, a single patch can be used, parameterized by the coordinates $(\xi^s,\eta^s)$:
\begin{equation}
    \begin{split}\label{eq:nurbsstruct}
      \mathbf{x}^s (\xi^s, \eta^s) &= \sum_{i=1}^{n_1} \sum_{j=1}^{n_2} R^p_i(\xi^s) R^p_j(\eta^s) \mathbf{x}^s_{ij} =
      \sum_{i=1}^{n_1} \sum_{j=1}^{n_2} R^p_i(\xi^s) R^p_j(\eta^s) ({\mathbf{x}_0^s}_{ij} + \mathbf{u}^s_{ij}) \\
       &= \sum_{i=1}^{n_1} \sum_{j=1}^{n_2} R^p_i(\xi^s) R^p_j(\eta^s) {\mathbf{x}_0^s}_{ij} + \sum_{i=1}^{n_1} \sum_{j=1}^{n_2} R^p_i(\xi^s) R^p_j(\eta^s) \mathbf{u}^s_{ij} \\
       &= \mathbf{x}_0^s (\xi^s, \eta^s) + \mathbf{u}^s (\xi^s, \eta^s) .
    \end{split}
\end{equation}
For the sake of readability, the notations are simplified in what follows, as $R^p_i R^p_j = R_k$, the global index $k$ denoting the tensorial product of the functions of indices $(i,j)$.

\subsubsection{Galerkin method}

Semi-discrete equations are obtained by injecting the isogeometric representation in Eq.~(\ref{eq:hyperelasticity_weak}) and choosing a NURBS basis function as test function $\boldsymbol\varphi = R_j$, the spatial integration being carried out in the parametric domain.  One obtains a set of differential ordinary equations:
\begin{align}
\label{eq:hyperelasticity_weak2}
\sum_{i=1}^n \left [ \int_{\hat{\Omega}^s} \rho^s  R_j \, R_i \, |\mbox{J}_{\Omega^s}| \,\mbox{d} \hat{\Omega} \right ] \ddot{\mathbf{u}_i^s} &+
  \int_{\hat{\Omega}^s} (\mathbf{F}\mathbf{S})(\mathbf{u}^s) : \nabla R_j  \, |\mbox{J}_{\Omega^s}| \,\mbox{d} \hat{\Omega} = \\ &\int_{\hat{\Gamma}^s} \boldsymbol\tau^s \cdot R_j \, |\mbox{J}_{\Gamma^s}| \, \mbox{d}  \hat{\Gamma} + \int_{\hat{\Omega}^s} \rho^s \mathbf{g} \cdot \nabla R_j  \, |\mbox{J}_{\Omega^s}| \,\mbox{d} \hat{\Omega},
\end{align}
where $\mbox{J}_{\Omega^s}$ and $\mbox{J}_{\Gamma^s}$ are the Jacobians of the geometric transformation from the parametric domain to the physical domain. After evaluation of the spatial integrals by Gauss-Legendre quadratures, the previous equations can be summarized as:
\begin{equation}
\mathcal{M}^s \ddot{U}^s  + \mathcal{K}^s (U^s) = F^s ,
\label{eq:variational_matrix2}
\end{equation}
which is formally similar to Eq.~(\ref{eq:variational_matrix}) for the membrane.

\subsection{Time integration}

The set of non-linear ordinary differential equations in Eq.~(\ref{eq:variational_matrix}) and Eq.~(\ref{eq:variational_matrix2}) are integrated in time using the Newmark scheme~\cite{RAO2011241}
which defines a second-order accurate unconditionnaly stable integration (for linear problems).

The mild non-linearity in Eq.~(\ref{eq:variational_matrix}), due to the depencency of the tension coefficient with respect to the membrane length, can be easily managed by evaluating the stiffness matrix at the previous time step, which transforms Eq.~(\ref{eq:variational_matrix}) to a linear system. For the Turek's case, the non-linear system in Eq.~(\ref{eq:variational_matrix2}) is solved directly using a Newton-Raphson method at each time step.

\clearpage


\section{Flow solver}\label{sec:scheme}


\subsection{Model}

In order to simulate fluid flow problems with moving geometries, we consider the compressible Navier-Stokes equations in Arbitrary Lagrangian-Eulerian (ALE) form, with a computational domain which moves with a generic velocity $\mathbf{v}^f$. The divergence form of the equations is the following:
\begin{equation}
\frac{\partial \mathbf{w}}{\partial t}  + \nabla \cdot \big( \mathbf{f}_c - \mathbf{f}_v \big)  - \mathbf{v}^f \cdot \nabla \mathbf{w} = 0,
\label{eq:cons}
\end{equation}
with $\mathbf{w}$ the vector of conservative variables, $\mathbf{f}_c$ the convective flux, and $\mathbf{f}_v$ the viscous flux:
\begin{equation}
\mathbf{w} = \begin{pmatrix}
\rho^f \\
\rho^f u_1 \\
\rho^f u_2 \\
\rho^f e
\end{pmatrix},
\quad
\mathbf{f}_{c,i} = \begin{pmatrix}
\rho^f u_i \\
\rho^f u_1 u_i + p \delta_{1i} \\
\rho^f u_2 u_i  + p \delta_{2i} \\
\rho^f u_i \big( e + \frac{p}{\rho^f} \big)
\end{pmatrix},
\quad
\mathbf{f}_{v,i} = \begin{pmatrix}
0 \\
\tau_{1i} \\
\tau_{2i} \\
u_k \tau_{ki}  - q_i
\end{pmatrix},
\label{eq:navierstokes1}
\end{equation}
with $\rho^f$ the fluid density, $\rho^f u_1$ and $\rho^f u_2$ the momentum components and $\rho^f e$ the total energy.
$\tau_{ij}$ is the viscous stress tensor and $q_i$ is the thermal conduction flux, defined as:
\begin{equation}
\tau_{ij} = \mu \bigg(	\frac{\partial u_i}{\partial x_j} + \frac{\partial u_j}{\partial x_i} \bigg) - \frac{2}{3} \mu \frac{\partial u_k}{\partial x_k}\delta_{ij},
\label{eq:tau}
\end{equation}
\begin{equation}
q_i = - \gamma \frac{\mu}{Pr} \frac{\partial e}{\partial x_i},
\label{eq:heat}
\end{equation}
where $\gamma = 1.4$, $Pr = 0.72$ and $\mu$ is the dynamic viscosity.
To account for the second-order derivatives, the state equations in ALE form are written as a system of first order equations:
\begin{equation}
	\left\{
	\begin{aligned}
		&\frac{\partial \mathbf{w}}{\partial t}  + \nabla \cdot \mathbf{f}_c (\mathbf{w}) - \nabla \cdot \mathbf{f}_v (\mathbf{w},\mathbf{g}) - \mathbf{v}^f \cdot \nabla \mathbf{w} = 0, \\
		&\mathbf{g} - \nabla \mathbf{w} = 0.
	\end{aligned}
	\right.
\label{eq:cons_final}
\end{equation}
where $\mathbf{g}$ is the gradient of the conservative variables.


\subsection{Isogeometric discretization}\label{sec:iga-dg}

To define a fluid domain $\Omega^f$ whose boundary exactly matches the geometry of the structure $\Omega^s$, the same NURBS representation has to be adopted. More precisely, we make the assumption that the fluid domain is composed of a set of bivariate NURBS surfaces, as illustrated in Fig.~(\ref{fig:nurbssurf}). Each NURBS surface of the fluid domain parameterized by the coordinate system $(\xi^f,\eta^f)$ is therefore defined using a tensorial representation:
\begin{equation}\label{eq:nurbscurve3}
  \mathbf{x}^f (\xi^f,\eta^f) = \sum_{i=1}^{n} \sum_{j=1}^{n} R^p_i(\xi^f) R^p_j(\eta^f) \mathbf{x}^f_{ij} \, ,
\end{equation}
with $\mathbf{x}^f = (x^f,y^f)$.
However, contrary to the structure case, this representation is not employed directly to solve the flow problem in Eq.~(\ref{eq:cons}). Indeed, we intend to use a DG method, well suited to convection-dominant problems, which requires to define for each element a set of basis functions, discontinuous at the interfaces. Therefore, each NURBS surface is first split into a set of rational Bézier patches using the Bézier extraction procedure. Each resulting rational Bézier patch is then defined by its own basis of degree $p$ and, therefore, can be employed as a discontinuous element in a DG framework. It counts $(p+1) \times (p+1)$ control points. The set of rational Bézier patches is geometrically identical to the original NURBS surface, although described by a different parameterization. The whole process to generate a DG-compliant computational domain from NURBS surfaces is detailed in~\cite{Duvigneau_18}.
As result of the Bézier extraction, each rational Bézier patch $\Omega^f_k$, considered as element for the fluid analysis, is represented by the rational Bernstein basis~\cite{Piegl_95}:
\begin{equation}
\tilde{R}_{i_1 i_2}^p(\xi^f,\eta^f) = \frac{B_{i_1}^p(\xi^f) \, B_{i_2}^p(\eta^f) \, \omega_{i_1 i_2}}{\sum_{j_1 = 1}^{p+1} \sum_{j_2 = 1}^{p+1} B_{j_1}^p(\xi^f) \, B_{j_2}^p(\eta^f) \, \omega_{j_1 j_2}},
\label{eq:bezier}
\end{equation}
defined in the parametric domain $\widehat{\Omega}^f = [0,1]^2$. The coefficients $\omega_{i_1 i_2}$ are the weights and $(B_{i}^p)_{i=1,\ldots,p+1}$ are the Bernstein polynomials of degree $p$.


\subsection{Discontinuous Galerkin method}

The computational mesh is thus a tessellation of rational B\'ezier patches $(\Omega^f_k)_{k=1,\ldots,N}$. The first equation in Eq.~(\ref{eq:cons_final}) yields the following weak formulation in each patch $\Omega_k^f$ after integration by parts:
\begin{align}
	\int_{\Omega^f_k}  \varphi \frac{\partial \mathbf{w}}{\partial t}  \,\mbox{d} \Omega -
	\int_{\Omega^f_k} \nabla\varphi \cdot & (\mathbf{f}_c - \mathbf{f}_v)  \,\mbox{d} \Omega +
	\int_{\Omega^f_k} \nabla \cdot \big( \varphi \mathbf{v}^f \big)  \mathbf{w} \,\mbox{d} \Omega \notag \\
	&+ \oint_{\partial \Omega^f_k} \varphi \big( \mathbf{f}_{ale}^* - \mathbf{f}_v^* \big) \, \mbox{d} \Gamma = 0 \qquad \forall \varphi \in \hat{H}^1(\Omega^f_k),
\label{eq:ale_base}
\end{align}
where $\hat{H}^1(\Omega^f)$ is the broken Sobolev space defined over the fluid domain. The fluxes $\mathbf{f}_{ale}^* = (\mathbf{f}_c - \mathbf{v}^f \mathbf{w})^* \cdot \mathbf{n}$ and $\mathbf{f}_v^*$ are respectively the interfacial ALE and viscous normal fluxes. By applying the Reynolds transport theorem for moving domains, one obtains:
\begin{align}
	\frac{\mbox{d}}{\mbox{d}t} \int_{\Omega^f_k} \varphi \, \mathbf{w} \, \mbox{d} \Omega -
	\int_{\Omega^f_k} \nabla \varphi  \cdot \big( \mathbf{f}_{ale} - \mathbf{f}_v \big) \, \mbox{d} \Omega +
	\oint_{\partial \Omega^f_k} \varphi \big( \mathbf{f}_{ale}^* - \mathbf{f}_v^* \big) \, \mbox{d} \Gamma = 0 ,
\label{eq:Nguyen}
\end{align}
with $\mathbf{f}_{ale} = \mathbf{f}_c - \mathbf{v}^f \mathbf{w}$ the ALE flux.
The weak formulation for the second equation in Eq.~(\ref{eq:cons_final}) writes:
\begin{align}
    \int_{\Omega^f_k} \varphi \, \mathbf{g} \, \mbox{d} \Omega +
	\int_{\Omega^f_k} \nabla \varphi \cdot \mathbf{f}_w \, \mbox{d} \Omega -
	\oint_{\partial \Omega^f_k} \varphi \, \mathbf{f}_w^*  \, \mbox{d} \Gamma = 0 ,
\label{eq:Nguyen2}
\end{align}
where $\mathbf{f}_w$ and $\mathbf{f}_w^*$ are respectively the flux and the interfacial normal flux for the gradient fields.

According to the isogeometric paradigm, a unified mathematical representation is adopted for the geometry of the fluid domain $\mathbf{x}^f$, its velocity $\mathbf{v}^f$ and the solution fields $\mathbf{w} $ and $\mathbf{g}$~\cite{Pezzano_Duvigneau_21}. Thus, in each patch $\Omega^f_k$, these fields are expressed using the rational Bernstein basis:
\begin{equation}
\begin{pmatrix}
\mathbf{x}^f \\
\mathbf{v}^f \\
\mathbf{w} \\
\mathbf{g}
\end{pmatrix} =
\sum_{i=1}^{(p+1)^2} \tilde{R}_i(\xi,\eta) \,
\begin{pmatrix}
\mathbf{x}^f_i \\
\mathbf{v}^f_i \\
\mathbf{w}_i \\
\mathbf{g}_i
\end{pmatrix},
\end{equation}
with a simplification of notation and a change of indices to facilitate reading. The time-dependent geometry is described by the control points $\mathbf{x}^f_i$ and their velocity $\mathbf{v}^f_i$, whereas $\mathbf{w}_i$ and $\mathbf{g}_i$ are respectively the degrees of freedom (DOF) of the conservative variables and their gradients.
Finally, the system composed of Eq.~(\ref{eq:Nguyen}) and Eq.~(\ref{eq:Nguyen2}) is discretized using this representation and choosing $\varphi = \tilde{R}_j$, as detailed in~\cite{Pezzano_Duvigneau_21}:
\begin{subequations}
	\begin{empheq}[left=\empheqlbrace]{alignat=2}
	&\frac{\mbox{d}}{\mbox{d}t} \bigg( \sum_{i=1}^{(p+1)^2} \left [ \int_{\widehat{\Omega}^f} \tilde{R}_j  \tilde{R}_i \, | \mbox{J}_{\Omega^f_k} |  \,\mbox{d} \widehat{\Omega} \right ] \mathbf{w}_i  \bigg) &&=
	\int_{\widehat{\Omega}^f} \nabla \tilde{R}_j  \cdot \big( \mathbf{f}_{ale} - \mathbf{f}_v \big) \, | \mbox{J}_{\Omega^f_k} |  \,\mbox{d} \widehat{\Omega} \notag \\
	& &&- \oint_{\partial \widehat{\Omega}^f} \tilde{R}_j \big( \mathbf{f}_{ale}^* - \mathbf{f}_v^* \big) \, | \mbox{J}_{\Gamma^f_k} |  \,\mbox{d} \widehat{\Gamma}, \label{eq:weak_w} \\
	& \sum_{i=1}^{(p+1)^2} \left [ \int_{\widehat{\Omega}^f} \tilde{R}_j  \tilde{R}_i \, | \mbox{J}_{\Omega^f_k} |  \,\mbox{d} \widehat{\Omega} \right ] \mathbf{g}_i =
	- \int_{\widehat{\Omega}^f} &&\nabla \tilde{R}_j \, \mathbf{f}_w \, | \mbox{J}_{\Omega^f_k} |  \,\mbox{d} \widehat{\Omega} +
	\oint_{\partial \widehat{\Omega}^f} \tilde{R}_j \, \mathbf{f}_w^* \, | \mbox{J}_{\Gamma^f_k} |  \,\mbox{d} \widehat{\Gamma}.
	\label{eq:weak_g}
	\end{empheq}
	\label{eq:weak_ale}
\end{subequations}
The convective numerical flux $\mathbf{f}_{ale}^*$ is computed through a modified HLL Riemann solver~\cite{Pezzano_Duvigneau_21}, whereas the LDG flux~\cite{Cockburn_Shu_98} is used for $\mathbf{f}_v^*$ and $\mathbf{f}_w^*$. Integration is carried out in the parametric domain $\widehat{\Omega}^f$ using Gauss-Legendre formulas. The isogeometric mapping from the parametric domain to the physical domain yields the non-linear Jacobians $\mbox{J}_{\Omega^f_k}$ and $\mbox{J}_{\Gamma^f_k}$, which are specific for each patch. Eq.~\eqref{eq:weak_g} is solved within each time iteration in a decoupled manner and the system of equations in Eq.~\eqref{eq:weak_ale} can then be rewritten as:
\begin{equation}
\dot{\big( \mathcal{M}^f \mathbf{W} \big)} = \mathcal{R}(\mathbf{W},\mathbf{G},\mathbf{V}^f),
\label{eq:compact}
\end{equation}
where $\mathcal{M}^f$ is the mass matrix, $\mathbf{W} = (\mathbf{w}_i)_{i=1,\ldots,I}$ the vector of DOF, $\mathbf{G} = (\mathbf{g}_i)_{i=1,\ldots,I}$ the vector of gradients, $\mathbf{V}^f = (\mathbf{v}_i)_{i=1,\ldots,I}$ the vector of grid velocity, with $I = (p+1)^2 N$ the total number of DOFs, and $\mathcal{R}$ the non-linear right-hand side vector of Eq.~\eqref{eq:weak_w}. The resulting set of ordinary differential equations in Eq.~\eqref{eq:compact} is integrated in time using a 4th-order four-stage explicit Runge-Kutta (RK) scheme. Note that the mass matrix is included in the time derivative, therefore it should be computed and inverted at each sub-iteration of the RK scheme. Hopefully, it is block-diagonal thanks to the compactness of the DG formulation.


\subsection{Local refinement}

To accurately capture flow features, such as boundary layers or vortex shedding phenomena, the fluid grid is locally refined. From geometrical viewpoint, this can be easily achieved by splitting some selected rational Bézier patches using the knot insertion procedure presented in section~(\ref{sec:iga-dg}). Thanks to the NURBS properties, this refinement strategy does not alter the geometry, which is a critical property for the refinement in the vicinity of curved boundaries.

From analysis viewpoint, the numerical scheme has to cope with hanging nodes at the faces neighbouring refined patches (see Fig.~(\ref{fig:turek_mesh}) for an illustration). Hopefully, discontinuous Galerkin schemes can naturally deal with such situations, by computing the integral of fluxes along each of the sub-faces generated with the hanging node. Details about this procedure can be found in~\cite{Duvigneau_20}.


\section{Coupling algorithm}\label{sec:coupling_algorithm}


The computation of the coupled fluid-structure system requires to solve the equations governing the fluid and the structure, as well as the following coupling conditions at the interface:
\begin{equation}
	\left\{
	\begin{aligned}
		& \mathbf{x}^s = \mathbf{x}^f \\
		&\boldsymbol\tau^s = \boldsymbol\tau^f = (-p \mathbf{I} + \boldsymbol{\tau}) \cdot \mathbf{n} ,
	\end{aligned}
	\right.
\label{eq:coupling}
\end{equation}
The first condition expresses the compatibility between the displacement of the fluid and the structure (kinematic condition), whereas the second one represents the mechanical equilibrium of the interface (dynamic condition). For the specific case of the membrane, the later degenerates to the equality of the fluid pressure and the normal forces applied to the membrane.

To ensure these coupling conditions, it is necessary to tranfer information between the two disciplines, specifically the displacement of the interface and the forces exerted on each side, and define a strategy to organize this transfer at discrete time steps. These aspects are detailed in the following sections.

\subsection{Transfer of displacement}

\subsubsection{Interface displacement}

We adopt here a classical approach in fluid-structure interaction, that consists in computing the kinematics of the structural system, by solving the elasticity problem subject to the fluid load, and then transfering the position and velocity of the interface to the fluid system. The position is used to update the fluid mesh, whereas the velocity is necessary to compute the ALE flux in~Eq.(\ref{eq:weak_w}), as detailed in the previous section.

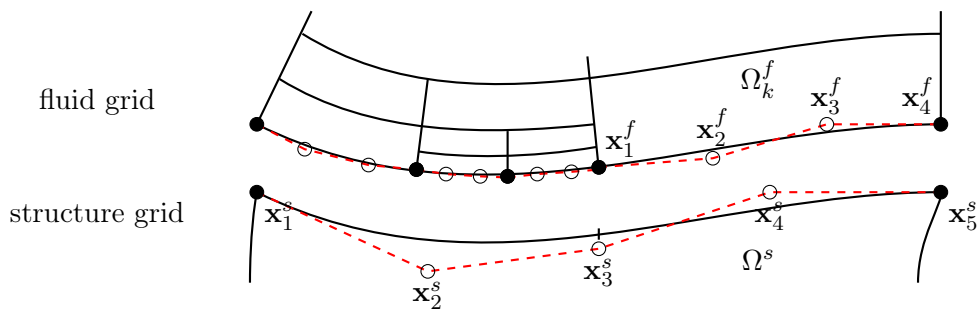
\begin{figure}[!ht]
\centering
\begin{center}
\begin{tikzpicture}[scale=3]
\draw[thick]  (0,0) .. controls (1,-0.5) and (2,0) .. (3,0);
\draw[dashed,red,thick]  (0,0) -- (0.75,-0.35) -- (1.5,-0.25) -- (2.25,0) -- (3,0);
\node[circle,draw,scale=0.5,fill] at (0,0){};
\node[below] at (0.1,0) {$\mathbf{x}^s_{1}$} ;
\node[circle,draw,scale=0.5] at (0.75,-0.35){};
\node[below] at (0.75,-0.35) {$\mathbf{x}^s_{2}$} ;
\node[circle,draw,scale=0.5] at (1.5,-0.25){};
\node[below] at (1.5,-0.25) {$\mathbf{x}^s_{3}$} ;
\node[circle,draw,scale=0.5] at (2.25,0){};
\node[below] at (2.25,0) {$\mathbf{x}^s_{4}$} ;
\node[circle,draw,scale=0.5,fill] at (3,0){};
\node[below] at (3.1,0) {$\mathbf{x}^s_{5}$} ;
\draw[thick]  (1.5,-0.21) -- (1.5,-0.16);

\node at (2.2,-0.3) {$\Omega^s$} ;

\draw[thick]  (0,0) .. controls (-0.01,-0.05) and (-0.03,-0.2) .. (-0.03,-0.4);
\draw[thick]  (3,0) .. controls (3.01,-0.05) and (2.9,-0.2) .. (2.9,-0.4);

\draw[thick]  (0,0+0.3) .. controls (1,-0.5+0.3) and (2,0+0.3) .. (3,0+0.3);
\draw[thick]  (0.2,0+0.7) .. controls (1,-0.5+0.7) and (2,0+0.7) .. (3,0+0.7);
\draw[thick]  (0,0+0.3) -- (0.24,0+0.8);
\draw[thick]  (3,0+0.3) -- (3,0+0.8);
\draw[thick]  (1.5,0.1) -- (1.45,0.6);

\draw[thick]  (0.1,0+0.5) .. controls (0.5,0.25) and (1,0.25) .. (1.48,0+0.3);
\draw[thick]  (0.7,0.1) -- (0.75,0.5);
\draw[thick]  (0.71,0.18) .. controls (0.9,0.15) and (1.2,0.15) .. (1.49,0.2);
\draw[thick]  (1.1,0.07) -- (1.1,0.27);

\draw[dashed,red,thick] (1.5,0.11) -- (2,0.15) -- (2.5,0+0.3) -- (3,0+0.3);
\node[circle,draw,scale=0.5,fill] at (1.5,0.11){};
\node[circle,draw,scale=0.5] at (2,0.15){};
\node[circle,draw,scale=0.5] at (2.5,0.3){};
\node[circle,draw,scale=0.5,fill] at (3,0.3){};

\draw[dashed,red,thick]  (0,0.3) -- (0.21,0.19) -- (0.49,0.12) -- (0.7,0+0.1);
\node[circle,draw,scale=0.5,fill] at (0,0.3){};
\node[circle,draw,scale=0.5] at (0.21,0.19){};
\node[circle,draw,scale=0.5] at (0.49,0.12){};
\node[circle,draw,scale=0.5,fill] at (0.7,0.1){};

\draw[dashed,red,thick]  (0.7,0.1) -- (0.83,0.08) -- (0.98,0.07) -- (1.1,0.065);
\node[circle,draw,scale=0.5,fill] at (0.7,0.1){};
\node[circle,draw,scale=0.5] at (0.83,0.08){};
\node[circle,draw,scale=0.5] at (0.98,0.07){};
\node[circle,draw,scale=0.5,fill] at (1.1,0.07){};

\draw[dashed,red,thick]  (1.1,0.065) -- (1.23,0.08) -- (1.38,0.09) -- (1.5,0.1);
\node[circle,draw,scale=0.5,fill] at (1.1,0.07){};
\node[circle,draw,scale=0.5] at (1.23,0.08){};
\node[circle,draw,scale=0.5] at (1.38,0.09){};
\node[circle,draw,scale=0.5,fill] at (1.5,0.11){};

\node at (2.2,0.5) {$\Omega^f_k$} ;
\node[above] at (1.6,0.11) {$\mathbf{x}^f_{1}$} ;
\node[above] at (2,0.16) {$\mathbf{x}^f_{2}$} ;
\node[above] at (2.5,0.3) {$\mathbf{x}^f_{3}$} ;
\node[above] at (2.9,0.3) {$\mathbf{x}^f_{4}$} ;

\node[] at (-0.7,0.4) {fluid grid} ;
\node[] at (-0.7,-0.1) {structure grid} ;

\end{tikzpicture}
\end{center}
\caption{Configuration of the fluid-structure interface}
\label{fig:transfer_displ}
\end{figure}

In the proposed isogeometric framework, the geometry of the structure is represented by NURBS curves or surfaces, whereas the fluid grid is composed by a set of disconnected rational Bézier patches, with possible local refinement. We make the assumption that the baseline NURBS surfaces used to generate the fluid patches (by Bézier extraction) coincide initially with the structure geometry at the interface. This condition is easy to fulfil at the construction of the computational domain, and will be illustrated in application sections. After Bézier extraction and subsequent refinement steps, one obtains the situation depicted in~Fig.(\ref{fig:transfer_displ}).

\begin{figure}[!ht]
\centering
\begin{tikzpicture}[scale=3]
\draw[thick]  (0,0) .. controls (1,-0.5) and (2,0) .. (3,0);
\draw[dashed,red,thick]  (0,0) -- (0.75,-0.35) -- (1.5,-0.25) -- (2.25,0) -- (3,0);
\node[circle,draw,scale=0.5,fill] at (0,0){};
\node[circle,draw,scale=0.5] at (0.75,-0.35){};
\node[circle,draw,scale=0.5] at (1.5,-0.25){};
\node[circle,draw,scale=0.5] at (2.25,0){};
\node[circle,draw,scale=0.5,fill] at (3,0){};
\draw[thick]  (1.5,-0.21) -- (1.5,-0.16);

\draw[thick]  (0,0+0.3) .. controls (1,-0.5+0.3) and (2,0+0.3) .. (3,0+0.3);

\draw[dashed,red,thick] (1.5,0.11) -- (2,0.15) -- (2.5,0+0.3) -- (3,0+0.3);
\node[circle,draw,scale=0.5,fill] at (1.5,0.11){};
\node[circle,draw,scale=0.5] at (2,0.15){};
\node[circle,draw,scale=0.5] at (2.5,0.3){};
\node[circle,draw,scale=0.5,fill] at (3,0.3){};

\draw[dashed,red,thick]  (0,0.3) -- (0.5,0.05) -- (1,0.05) -- (1.5,0.11);
\node[circle,draw,scale=0.5,fill] at (0,0.3){};
\node[circle,draw,scale=0.5] at (0.5,0.05){};
\node[circle,draw,scale=0.5] at (1,0.05){};

\begin{scope}[shift={(0,0.3)},rotate=0]
\draw[thick]  (0,0+0.3) .. controls (1,-0.5+0.3) and (2,0+0.3) .. (3,0+0.3);

\draw[dashed,red,thick] (1.5,0.11) -- (2,0.15) -- (2.5,0+0.3) -- (3,0+0.3);
\node[circle,draw,scale=0.5,fill] at (1.5,0.11){};
\node[circle,draw,scale=0.5] at (2,0.15){};
\node[circle,draw,scale=0.5] at (2.5,0.3){};
\node[circle,draw,scale=0.5,fill] at (3,0.3){};

\draw[dashed,red,thick]  (0,0.3) -- (0.21,0.19) -- (0.49,0.12) -- (0.7,0+0.1);
\node[circle,draw,scale=0.5,fill] at (0,0.3){};
\node[circle,draw,scale=0.5] at (0.21,0.19){};
\node[circle,draw,scale=0.5] at (0.49,0.12){};
\node[circle,draw,scale=0.5,fill] at (0.7,0.1){};

\draw[dashed,red,thick]  (0.7,0.1) -- (0.95,0.07) -- (1.25,0.08) -- (1.5,0.11);
\node[circle,draw,scale=0.5,fill] at (0.7,0.1){};
\node[circle,draw,scale=0.5] at (0.95,0.07){};
\node[circle,draw,scale=0.5] at (1.25,0.08){};

\end{scope}

\begin{scope}[shift={(0,0.6)},rotate=0]
\draw[thick]  (0,0+0.3) .. controls (1,-0.5+0.3) and (2,0+0.3) .. (3,0+0.3);

\draw[dashed,red,thick] (1.5,0.11) -- (2,0.15) -- (2.5,0+0.3) -- (3,0+0.3);
\node[circle,draw,scale=0.5,fill] at (1.5,0.11){};
\node[circle,draw,scale=0.5] at (2,0.15){};
\node[circle,draw,scale=0.5] at (2.5,0.3){};
\node[circle,draw,scale=0.5,fill] at (3,0.3){};

\draw[dashed,red,thick]  (0,0.3) -- (0.21,0.19) -- (0.49,0.12) -- (0.7,0+0.1);
\node[circle,draw,scale=0.5,fill] at (0,0.3){};
\node[circle,draw,scale=0.5] at (0.21,0.19){};
\node[circle,draw,scale=0.5] at (0.49,0.12){};
\node[circle,draw,scale=0.5,fill] at (0.7,0.1){};

\draw[dashed,red,thick]  (0.7,0.1) -- (0.83,0.08) -- (0.98,0.07) -- (1.1,0.065);
\node[circle,draw,scale=0.5,fill] at (0.7,0.1){};
\node[circle,draw,scale=0.5] at (0.83,0.08){};
\node[circle,draw,scale=0.5] at (0.98,0.07){};
\node[circle,draw,scale=0.5,fill] at (1.1,0.07){};

\draw[dashed,red,thick]  (1.1,0.065) -- (1.23,0.08) -- (1.38,0.09) -- (1.5,0.1);
\node[circle,draw,scale=0.5,fill] at (1.1,0.07){};
\node[circle,draw,scale=0.5] at (1.23,0.08){};
\node[circle,draw,scale=0.5] at (1.38,0.09){};
\node[circle,draw,scale=0.5,fill] at (1.5,0.11){};
\end{scope}

\node[] at (-0.7,0.9) {second splitting} ;
\node[] at (-0.7,1.05) {fluid boundary} ;
\node[] at (-0.7,0.6) {first splitting} ;
\node[] at (-0.7,0.3) {Bézier curves} ;
\node[] at (-0.7,0) {NURBS curve} ;
\node[] at (-0.7,-0.15) {structure boundary} ;

\end{tikzpicture}
\caption{Bézier extraction and refinement step by step}
\label{fig:transfer_displ2}
\end{figure}
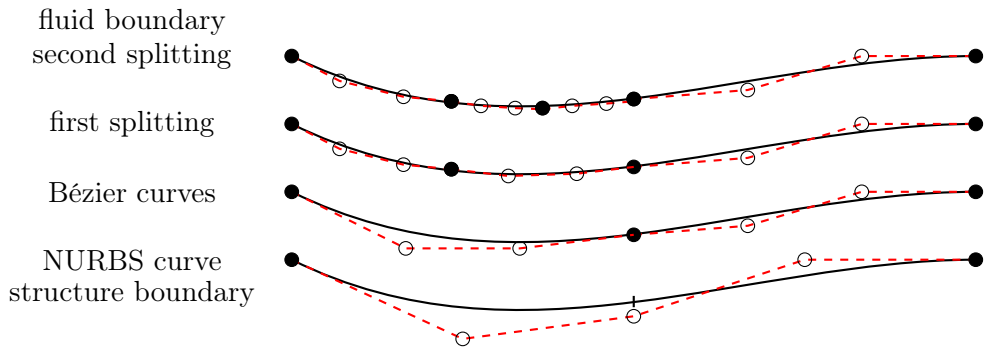

To ensure an exact matching of the fluid and structure domains, for any displacement of the structure, one just has to apply the Bézier extraction procedure to the structure boundary, followed by some splittings, until the current refinement level of the fluid boundary is reached. This is illustrated in~Fig.(\ref{fig:transfer_displ2}), which corresponds to the interface in Fig.(\ref{fig:transfer_displ}), for a NURBS curve of degree $p=3$ describing the structure boundary with $n=5$ control points. The curve counts a total number of $k=9$ knots, including four ($p+1$) equal knots at each extremity plus one interior knot. The NURBS curve is first transformed in two rational Bézier curves by inserting $p$ knots at the existing interior knot. These Bézier curves count $p+1$ control points and join at the point corresponding previously to the interior knot. Then, two successive splitting steps are achieved by inserting $p+1$ new knots at the mid-parameter of the curves to be split. Finally, the new fluid boundary is obtained and matches exactly the structure boundary, from geometrical point of view. In particular, the geometrical regularity of the membrane is maintained, even if each rational Bézier patch of the fluid grid is independent from one another.

\subsubsection{Interior grid displacement}

\begin{figure}[!ht]
\centering
\begin{tikzpicture}[scale=1.5]

\draw[ultra thick]  (0,0) .. controls (1,0.2) and (2,0) .. (3,0);
\draw[thick]  (-0.2,1.5) .. controls (1.1,1.5) and (2,1.8) .. (3.5,1.8);
\draw[thick]  (0,0) .. controls (0,0.5) and (-0.2,1) .. (-0.2,1.5);
\draw[thick]  (3,0) .. controls (3,0.6) and (3.5,1.2) .. (3.5,1.8);

\draw[thick]  (1.5,0.07) .. controls (1.5,0.7) and (1.6,1.2) .. (1.6,1.65);

\draw[thick]  (1.55,0.95) .. controls (2,0.92) and (2.5,0.9) .. (3.25,0.9);

\draw[thick]  (2.25,0.05) .. controls (2.3,0.3) and (2.35,0.7) .. (2.4,0.9);

\draw[dashed,red,thick]  (0,0) -- (1,0.2) -- (2,0) -- (3,0);

\node[circle,draw,scale=0.5,fill] at (0,0){};
\node[circle,draw,scale=0.5] at (1,0.2){};
\node[circle,draw,scale=0.5] at (2,0){};
\node[circle,draw,scale=0.5,fill] at (3,0){};

\draw[->,ultra thick,blue]  (1,0.2) -- (1,0.5);

\end{tikzpicture}
\caption{Initial Bézier patches and structure control points displacement}
\label{fig:grid_move0}
\end{figure}
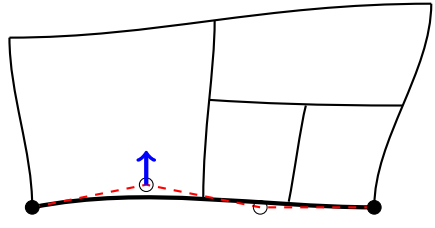

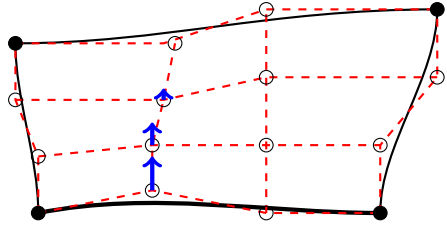
\begin{figure}[!ht]
\centering
\begin{tikzpicture}[scale=1.5]

\draw[ultra thick]  (0,0) .. controls (1,0.2) and (2,0) .. (3,0);
\draw[thick]  (-0.2,1.5) .. controls (1.1,1.5) and (2,1.8) .. (3.5,1.8);
\draw[thick]  (0,0) .. controls (0,0.5) and (-0.2,1) .. (-0.2,1.5);
\draw[thick]  (3,0) .. controls (3,0.6) and (3.5,1.2) .. (3.5,1.8);

\draw[dashed,red,thick]  (0,0) -- (1,0.2) -- (2,0) -- (3,0);
\draw[dashed,red,thick]  (0,0.5) -- (1,0.6) -- (2,0.6) -- (3,0.6);
\draw[dashed,red,thick]  (-0.2,1) -- (1.1,1) -- (2,1.2) -- (3.5,1.2);
\draw[dashed,red,thick]  (-0.2,1.5) -- (1.1,1.5) -- (2,1.8) -- (3.5,1.8);

\draw[dashed,red,thick]  (0,0) -- (0,0.5) -- (-0.2,1) -- (-0.2,1.5);
\draw[dashed,red,thick]  (1,0.2) -- (1,0.6) -- (1.1,1) -- (1.2,1.5);
\draw[dashed,red,thick]  (2,0.) -- (2,0.6) -- (2,1.2) -- (2,1.8);
\draw[dashed,red,thick]  (3,0.) -- (3,0.6) -- (3.5,1.2) -- (3.5,1.8);

\node[circle,draw,scale=0.5,fill] at (0,0){};
\node[circle,draw,scale=0.5] at (1,0.2){};
\node[circle,draw,scale=0.5] at (2,0){};
\node[circle,draw,scale=0.5,fill] at (3,0){};

\node[circle,draw,scale=0.5,fill] at (-0.2,1.5){};
\node[circle,draw,scale=0.5] at (1.2,1.5){};
\node[circle,draw,scale=0.5] at (2,1.8){};
\node[circle,draw,scale=0.5,fill] at (3.5,1.8){};

\node[circle,draw,scale=0.5] at (0,0.5){};
\node[circle,draw,scale=0.5] at (1,0.6){};
\node[circle,draw,scale=0.5] at (2,0.6){};
\node[circle,draw,scale=0.5] at (3,0.6){};

\node[circle,draw,scale=0.5] at (-0.2,1){};
\node[circle,draw,scale=0.5] at (1.1,1){};
\node[circle,draw,scale=0.5] at (2,1.2){};
\node[circle,draw,scale=0.5] at (3.5,1.2){};

\draw[->,ultra thick,blue]  (1,0.2) -- (1,0.5);
\draw[->,ultra thick,blue]  (1,0.6) -- (1,0.8);
\draw[->,ultra thick,blue]  (1.1,1) -- (1.1,1.1);

\end{tikzpicture}
\caption{Initial underlying NURBS surface and control points displacement with damping}
\label{fig:grid_move1}
\end{figure}

\begin{figure}[!ht]
\centering
\begin{tikzpicture}[scale=1.5]

\draw[ultra thick]  (0,0) .. controls (1,0.7) and (2,0) .. (3,0);
\draw[thick]  (-0.2,1.5) .. controls (1.1,1.5) and (2,1.8) .. (3.5,1.8);
\draw[thick]  (0,0) .. controls (0,0.5) and (-0.2,1) .. (-0.2,1.5);
\draw[thick]  (3,0) .. controls (3,0.6) and (3.5,1.2) .. (3.5,1.8);

\draw[thick]  (1.5,0.24) .. controls (1.5,0.7) and (1.6,1.2) .. (1.6,1.65);

\draw[thick]  (1.55,1) .. controls (2,0.9) and (2.5,0.9) .. (3.25,0.9);

\draw[thick]  (2.25,0.075) .. controls (2.3,0.3) and (2.35,0.7) .. (2.4,0.9);

\end{tikzpicture}
\caption{Final Bézier patches after refinement}
\label{fig:grid_move3}
\end{figure}
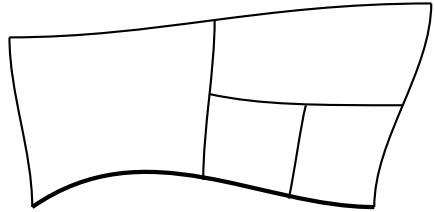

To propagate the displacement of the grid boundary to the interior of the fluid computational domain, an explicit geometrical procedure is considered. The properties of the NURBS representation make the approach robust although very simple. The control points that define the Bézier patches are moved according to the closest boundary control point with a linear damping of the displacement, to deform the grid in the vicinity of the structure only. Any linear transformation applied to the control points is equivalent to its application to the whole NURBS object (curve, surface, etc.)~\cite{Piegl_95}. Therefore, one can ensure that the patches remain admissible and the deformation does not lead to tangling: thanks to the convex hull property~\cite{Piegl_95}, a patch remains admissible if its control points are not overlaid. Moreover, the deformation is performed in a hierarchical manner, accounting for the local refinement. The move of the control points is first achieved for the baseline NURBS surface, and then Bézier extraction and splitting steps are applied to the whole grid, as illustrated in~Fig.(\ref{fig:grid_move0}-\ref{fig:grid_move3}). This approach has been detailed and tested on various problems in~\cite{Pezzano_Duvigneau_21}.

\begin{figure}[!ht]
\centering
\begin{tikzpicture}[scale=1.5]

\draw[ultra thick]  (0,0) .. controls (1,0.7) and (2,0) .. (3,0);
\draw[thick]  (-0.2,1.5) .. controls (1.1,1.5) and (2,1.8) .. (3.5,1.8);
\draw[thick]  (0,0) .. controls (0,0.5) and (-0.2,1) .. (-0.2,1.5);
\draw[thick]  (3,0) .. controls (3,0.6) and (3.5,1.2) .. (3.5,1.8);

\draw[dashed,red,thick]  (0,0) -- (1,0.5) -- (2,0) -- (3,0);
\draw[dashed,red,thick]  (0,0.5) -- (1,0.8) -- (2,0.6) -- (3,0.6);
\draw[dashed,red,thick]  (-0.2,1) -- (1.1,1.1) -- (2,1.2) -- (3.5,1.2);
\draw[dashed,red,thick]  (-0.2,1.5) -- (1.1,1.5) -- (2,1.8) -- (3.5,1.8);

\draw[dashed,red,thick]  (0,0) -- (0,0.5) -- (-0.2,1) -- (-0.2,1.5);
\draw[dashed,red,thick]  (1,0.5) -- (1,0.8) -- (1.1,1.1) -- (1.2,1.5);
\draw[dashed,red,thick]  (2,0.) -- (2,0.6) -- (2,1.2) -- (2,1.8);
\draw[dashed,red,thick]  (3,0.) -- (3,0.6) -- (3.5,1.2) -- (3.5,1.8);

\node[circle,draw,scale=0.5,fill] at (0,0){};
\node[circle,draw,scale=0.5] at (1,0.5){};
\node[circle,draw,scale=0.5] at (2,0){};
\node[circle,draw,scale=0.5,fill] at (3,0){};

\node[circle,draw,scale=0.5,fill] at (-0.2,1.5){};
\node[circle,draw,scale=0.5] at (1.2,1.5){};
\node[circle,draw,scale=0.5] at (2,1.8){};
\node[circle,draw,scale=0.5,fill] at (3.5,1.8){};

\node[circle,draw,scale=0.5] at (0,0.5){};
\node[circle,draw,scale=0.5] at (1,0.8){};
\node[circle,draw,scale=0.5] at (2,0.6){};
\node[circle,draw,scale=0.5] at (3,0.6){};

\node[circle,draw,scale=0.5] at (-0.2,1){};
\node[circle,draw,scale=0.5] at (1.1,1.1){};
\node[circle,draw,scale=0.5] at (2,1.2){};
\node[circle,draw,scale=0.5] at (3.5,1.2){};

\end{tikzpicture}
\caption{Final NURBS surface}
\label{fig:grid_move2}
\end{figure}
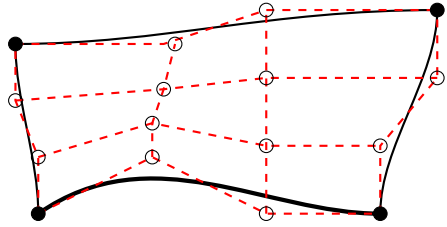

\subsection{Transfer of forces}

\begin{figure}[!ht]
\centering
\begin{tikzpicture}[scale=1.5]

\draw[ultra thick]  (0,0) .. controls (1,0.5) and (2,0) .. (3,0);
\draw[thick]  (-0.2,1.5) .. controls (1.1,1.5) and (2,1.8) .. (3.5,1.8);
\draw[thick]  (0,0) .. controls (0,0.5) and (-0.2,1) .. (-0.2,1.5);
\draw[thick]  (3,0) .. controls (3,0.6) and (3.5,1.2) .. (3.5,1.8);

\draw[thick]  (1.5,0.2) .. controls (1.5,0.7) and (1.6,1.2) .. (1.6,1.65);

\draw[thick]  (1.55,1) .. controls (2,0.9) and (2.5,0.9) .. (3.25,0.9);

\draw[thick]  (2.25,0.075) .. controls (2.3,0.3) and (2.35,0.7) .. (2.4,0.9);

\draw[thick,red,->]  (0,-0.3) .. controls (1,0.2) and (2,-0.3) .. (3,-0.3);
\node[red,right] at (3,-0.3){$\xi^s$};
\draw[thick,red,-]  (0.3,-0.25) -- (0.3,-0.1);
\node[red] at (0.45,-0.4){$\xi_q^s$};

\draw[thick,blue,->]  (0,0.3) .. controls (0.5,0.6) and (1,0.5) .. (1.5,0.5);
\node[blue,right] at (1.5,0.5){$\xi^f$};
\draw[thick,blue,-]  (0.3,0.35) -- (0.3,0.5);
\node[blue,above] at (0.3,0.5){$\xi_q^f$};

\node[draw,scale=0.5,fill] at (0.3,0.13){};
\node[] at (0.5,0.){$\mathbf{x}_q$};
\node[draw,scale=0.5,fill] at (0.7,0.22){};
\node[draw,scale=0.5,fill] at (1.1,0.22){};

\end{tikzpicture}
\caption{Parametric pairing of quadrature points for pressure integration}
\label{fig:pressure_pair}
\end{figure}
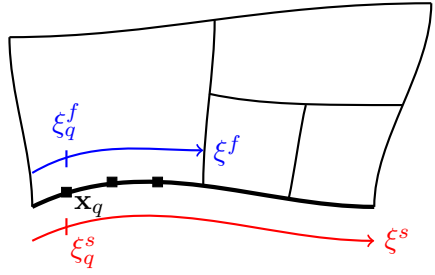

The second information to transfer is the fluid forces exerted on the structure surface. Specifically, the integral along the structure boundary
 $\Gamma^s$ of the external forces $\boldsymbol\tau^s$ multiplied by a structural test function is required to compute the right-hand side in~Eq.(\ref{eq:variational_galerkin}) and Eq.(\ref{eq:hyperelasticity_weak}). The pressure and velocity gradient fields are defined in the fluid domain, for each Bézier patch, and are discontinuous between two patches. Therefore, the integral over the structural parametric domain $\hat{\Gamma}^s$ has to be decomposed into a sum of integrals over subdomains $\hat{\Gamma}^s_k$, that correspond to the boundary of the different fluid patches:
\begin{equation}
\int_{\hat{\Gamma}^s} \boldsymbol\tau^s  \cdot R_j \; |\mbox{J}| \mbox{d}\hat{\Gamma} = \sum_{k \in \mathcal{B}^f} \int_{ \widehat{\Gamma}_k^s} \boldsymbol\tau^f \cdot R_j \, | \mbox{J} |  \,\mbox{d}\hat{\Gamma} ,
\label{eq:pressure_integral}
\end{equation}
where $\mathcal{B}^f$ corresponds to the set of rational Bézier curves at the fluid boundary. The integrals are evaluated by numerical quadratures, which necessitate to evaluate the quantities at $N_q$ quadrature points $\mathbf{x}_q$, whose parametric coordinates are  $\xi^s_q$ in the structural domain and $\xi^f_q$ in the fluid domain, yielding:
\begin{equation}
    \begin{aligned}
        \int_{\hat{\Gamma}^s} \boldsymbol\tau^s  \cdot R_j \; |\mbox{J}| \mbox{d}\hat{\Gamma} \approx
         & \sum_{k \in \mathcal{B}^f} \sum_{q=1}^{N_q} \omega_q \, \boldsymbol\tau^f(\xi^f_q) \cdot R_j(\xi^s_q)\, | \mbox{J}(\xi^s_q) | .
    \end{aligned}
\label{eq:pressure_integral_gauss}
\end{equation}
The evaluations of $\xi_q^s$ and $\xi^f_q$ are achieved by exploiting the underlying NURBS surface in the fluid domain that coincide with the NURBS curve representing the structure boundary, and applying again Bézier extraction and refinement operators. For each quadrature point, this allows to pair the parametric coordinate $\xi_q^f$ on the rational Bézier curve, where the flow fields are defined, and the parametric coordinate $\xi_q^s$ on the NURBS curve, where the structural test functions and Jacobians are defined, as depicted in~Fig.(\ref{fig:pressure_pair}).


\subsection{Coupling strategy}

\begin{figure}[!ht]
\centering
\begin{tikzpicture}[scale=1.7]

\draw[thick,->]  (0,0) -- (5,0);
\node[right] at (5,0.){$t$};

\draw[thick,->,red]  (-0.05,3) -- (-0.05,1.1);

\draw[very thick,->]  (0,1) -- (1.9,1);
\node[below] at (0,1){$U^s_k$};
\node[below] at (2,1){$U^s_{k+1}$};

\draw[thick,->,blue]  (2,1) -- (0,1.9);
\node[left] at (0,2){$y^f_k$};
\draw[thick,->,blue]  (2,1) -- (1,1.9);
\node[right] at (1,2){$y^f_{k+1/2}$};
\draw[thick,->,blue]  (2,1) -- (2,1.9);
\node[right] at (2.1,2){$y^f_{k+1}$};

\draw[thick,->,cyan]  (0,2) -- (0,2.9);
\node[above] at (0,3){$\mathbf{w}_k$};
\draw[very thick,->]  (0,3) -- (0.9,3);
\draw[thick,->,cyan]  (1,2) -- (1,2.9);
\node[above] at (1,3){$\mathbf{w}_{k+1/2}$};
\draw[very thick,->]  (1,3) -- (1.9,3);
\draw[thick,->,cyan]  (2,2) -- (2,2.9);
\node[above] at (2,3){$\mathbf{w}_{k+1}$};

\draw[thick,->,red]  (2.05,3) -- (2.05,1.1);

\draw[very thick,->]  (2.1,1) -- (3.9,1);
\node[below] at (4,1){$y^s_{k+2}$};

\draw[thick,->,blue]  (4,1) -- (2.1,1.9);
\draw[thick,->,blue]  (4,1) -- (3,1.9);
\node[right] at (3,2){$y^f_{k+3/2}$};
\draw[thick,->,blue]  (4,1) -- (4,1.9);
\node[right] at (4,2){$y^f_{k+2}$};

\draw[thick,->,cyan]  (2,2) -- (2,2.9);
\draw[very thick,->]  (2.1,3) -- (2.9,3);
\draw[thick,->,cyan]  (3,2) -- (3,2.9);
\node[above] at (3,3){$\mathbf{w}_{k+3/2}$};
\draw[very thick,->]  (3,3) -- (3.9,3);
\draw[thick,->,cyan]  (4,2) -- (4,2.9);
\node[above] at (4,3){$\mathbf{w}_{k+2}$};

\node[left] at (-1,1){structure};
\node[left] at (-1,2){fluid grid};
\node[left] at (-1,3){fluid variables};

\node[left,blue] at (0.9,1.4){displ.};
\node[left,blue] at (0.9,1.2){transfer};

\node[left,cyan] at (1.6,2.5){ALE};
\node[left,cyan] at (1.7,2.3){terms};
\node[left,red] at (0.,2.5){forces};
\node[left,red] at (0.,2.3){transfer};

\node[] at (2,0.5){\bf{Newmark integration}};
\node[] at (2,3.5){\bf{Runge-Kutta integration}};

\end{tikzpicture}
\caption{Fluid-structure coupling along time}
\label{fig:coupling}
\end{figure}
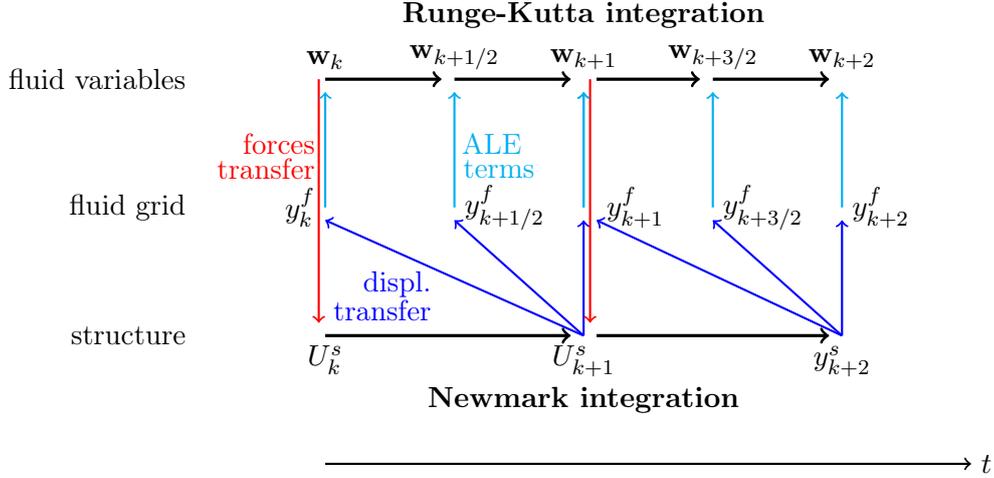

The transfer techniques being defined, one can now describe the coupling strategy, i.e. the organisation of the computations and transfers along time. A \emph{partitioned approach} is adopted, which is illustrated in Fig.(\ref{fig:coupling}) and summarized as follows at each time $t$:
\begin{enumerate}
    \item Compute the time-step $\delta t$ according to the stability condition for the fluid
    \item Given the flow fields at $t$, compute the structure load using Eq.(\ref{eq:pressure_integral_gauss})
    \item Compute the structure displacement $U^s_{k+1}$ and velocity $\dot{U}^s_{k+1}$ at time $t+\delta t$ by Newmark algorithm using Eq.(\ref{eq:variational_matrix}) or Eq.(\ref{eq:variational_matrix2})
    \item Compute the displacement $\mathbf{x}^f$ and velocity $\mathbf{v}^f$ of the fluid grid, by Bézier extraction and refinement from the structure displacements $U^s_{k}$ and $U^s_{k+1}$, at the sub-times $t+\delta t_k$ required by Runge-Kutta integrator
    \item Compute ALE terms in Eq.(\ref{eq:weak_w}), at the sub-times $t+\delta t_k$
    \item Compute the fluid solution $\mathbf{w}_{k+1}$ at time $t+\delta t$ using Runge-Kutta integrator by solving Eq.((\ref{eq:weak_w}-\ref{eq:weak_g}) at each sub-iteration
\end{enumerate}

 This algorithm corresponds to a \emph{loosely coupling} between the two disciplines. Indeed, there is no guarantee that the kinematic and dynamic conditions are fulfilled at each time step. Specifically, the displacements of the structure and the fluid domains are equal at the interface by construction, but the load on the structure is different from the forces exerted by the fluid at the end of each time step.
Consequently, the energy is not conserved during the coupling. This could be controled by introducing some sub-iterations at each time step, for instance based on a fixed-point approach, to converge towards the fulfillness of the two coupling conditions. This would correspond to a \emph{strong coupling} algorithm. This strategy has not been adopted in this work, due to its computational cost and its uselessness for the problems studied, as will be shown in the next sections, in the context of small time steps imposed by the explicit fluid solver.


\section{Membrane wing case}


\subsection{Computational setup}

The first testcase considered deals with the flow around an elastic membrane wing, whose physical configuration is depicted in Fig.~(\ref{fig:membrane}) and corresponds to the one studied experientally by  Rojratsirikul \textit{et al.}~\cite{Rojratsirikul2010} and numerically by Gordnier~\cite{Gordnier_09}.

The geometry is normalized by setting the membrane length $L=1$ and thickness $h=0.001$. The structure and fluid characteristics are chosen to respect the non-dimensional aeroelastic number $Ae = 100$, the mass ratio $m^\star = 0.589$, the freestream Reynolds number $R_\infty = 2500$ and Mach number $M_\infty = 0.1$:
\begin{equation}
Ae = \frac{Eh}{\frac{1}{2} \rho_\infty^f U_\infty^2 L} \quad m^\star = \frac{\rho^s h}{\rho_\infty^f L} \quad R_\infty = \frac{\rho_\infty^f U_\infty L}{\mu^f} \quad M_\infty = \frac{U_\infty}{\sqrt{\gamma p_\infty / \rho_\infty^f}} .
\end{equation}
Specifically, the following values are employed for the computations: $E=700$, $\rho^s = 825$, $\gamma=1.4$, $\rho^f_\infty = 1.4$, $U_\infty = 0.1$, $p_\infty = 1$, $\mu^f = 5.6 10^{-5}$.
No pretension is applied to the membrane, i.e. $T_0 = 0$. The angle of attack $\alpha$ varies between $4^\circ$ and $20^\circ$.

The structural computational domain is just a single NURBS curve, with unitary weights (i.e. a B-Spline curve), whose knots are uniformly distributed in the parametric domain and control points chosen to refine the membrane extremities.
The fluid computational domain is composed by a set of NURBS patches, with unitary weights (i.e. a set of B-Spline patches), constructed to match the membrane curve at rest and possibly refined afterwards, as shown in Fig.~(\ref{fig:grid}). A quite uniform grid area is defined around the membrane to accurately capture the flow detachements and vortices, whereas the grid becomes progressively coarser far from the membrane. The distance between the membrane and the exterior boundary is about 50 times the membrane length in all directions.

Regarding the boundary conditions, a zero displacement condition is prescribed for both membrane extremities. A no-slip condition is imposed for the fluid at the membrane surface, whereas a far-field condition is chosen for the exterior boundary, based on Riemann invariants.


\subsection{Grid convergence study}

\begin{figure}[!ht]
\centering
    \begin{subfigure}[b]{0.7\linewidth}
      \centering
      \includegraphics[width=\textwidth]{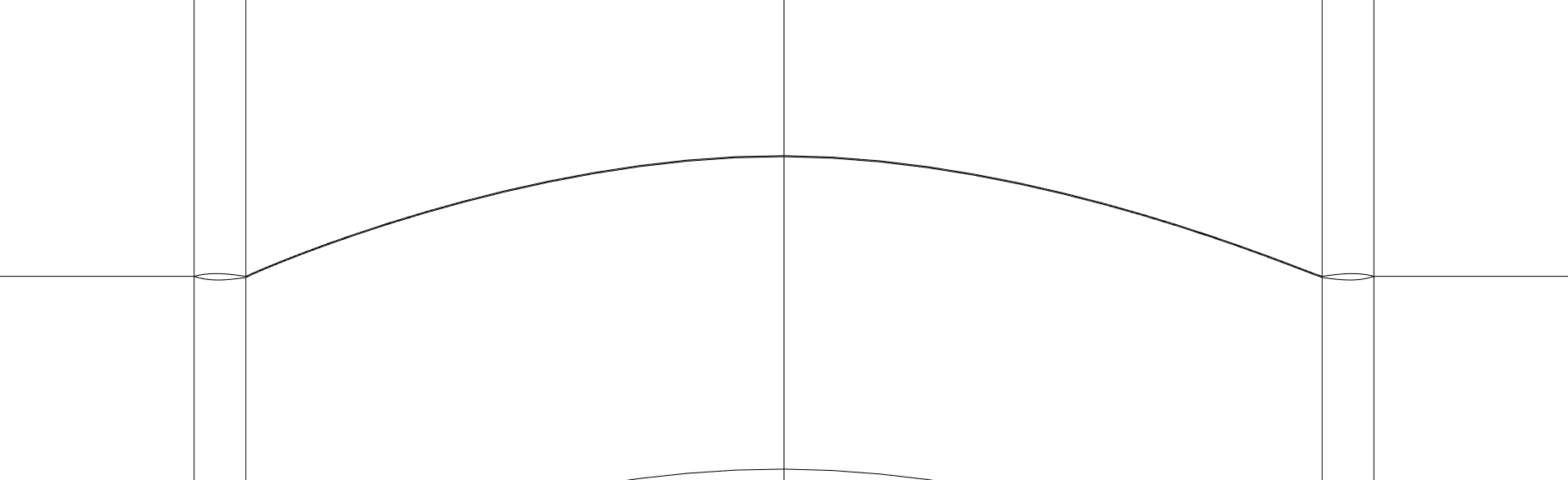}
      \caption{Very coarse grid}
      \label{fig:grid_verycoarse}
    \end{subfigure}\\
  \begin{subfigure}[b]{0.7\linewidth}
    \centering
    \includegraphics[width=\textwidth]{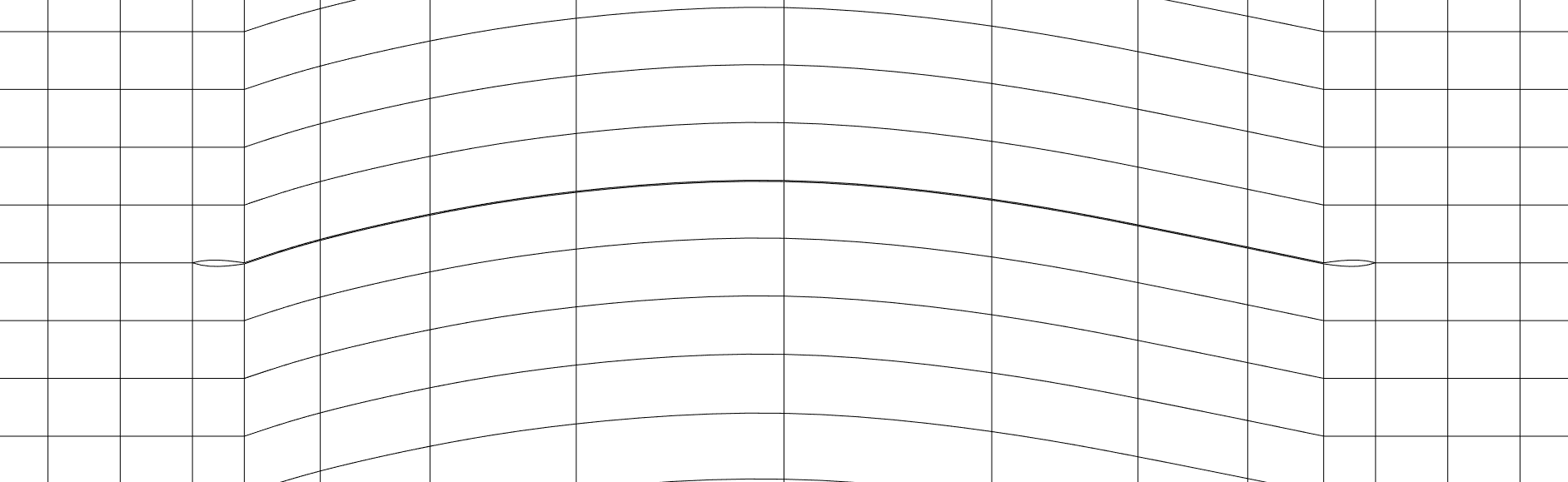}
    \caption{Coarse grid}
    \label{fig:grid_coarse}
  \end{subfigure}\\
  \begin{subfigure}[b]{0.7\linewidth}
    \centering
    \includegraphics[width=\textwidth]{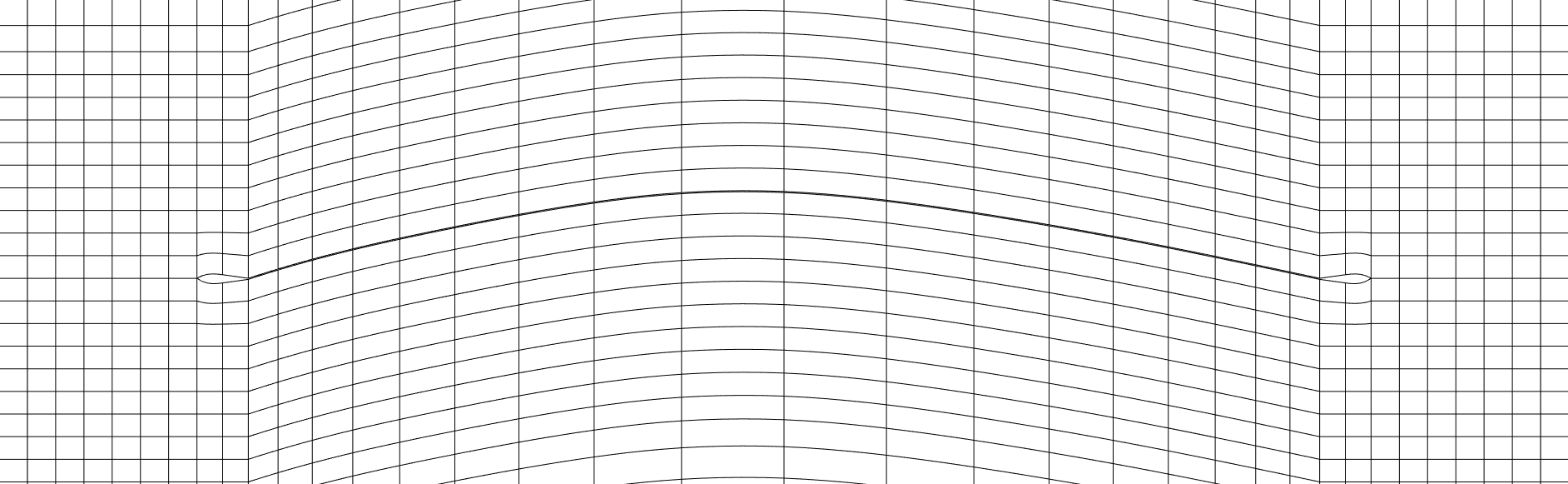}
    \caption{Medium grid}
    \label{fig:grid_medium}
  \end{subfigure}\\
  \begin{subfigure}[b]{0.7\linewidth}
    \centering
    \includegraphics[width=\textwidth]{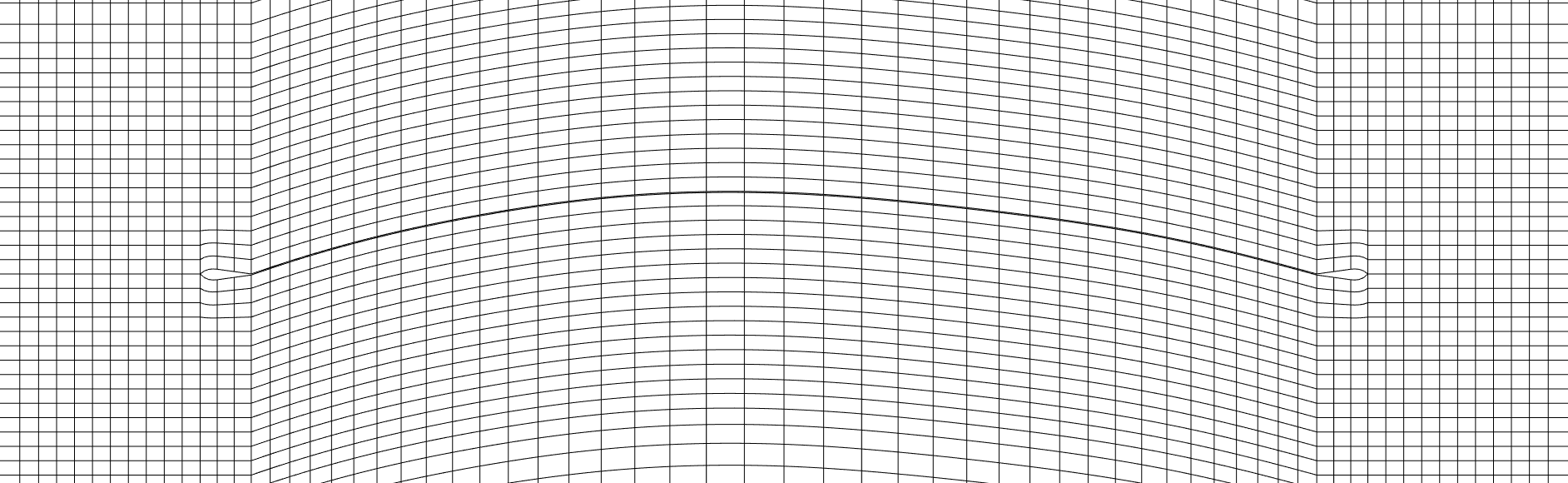}
    \caption{Fine grid}
    \label{fig:grid_fine}
  \end{subfigure}\\
  \begin{subfigure}[b]{0.7\linewidth}
    \centering
    \includegraphics[width=\textwidth]{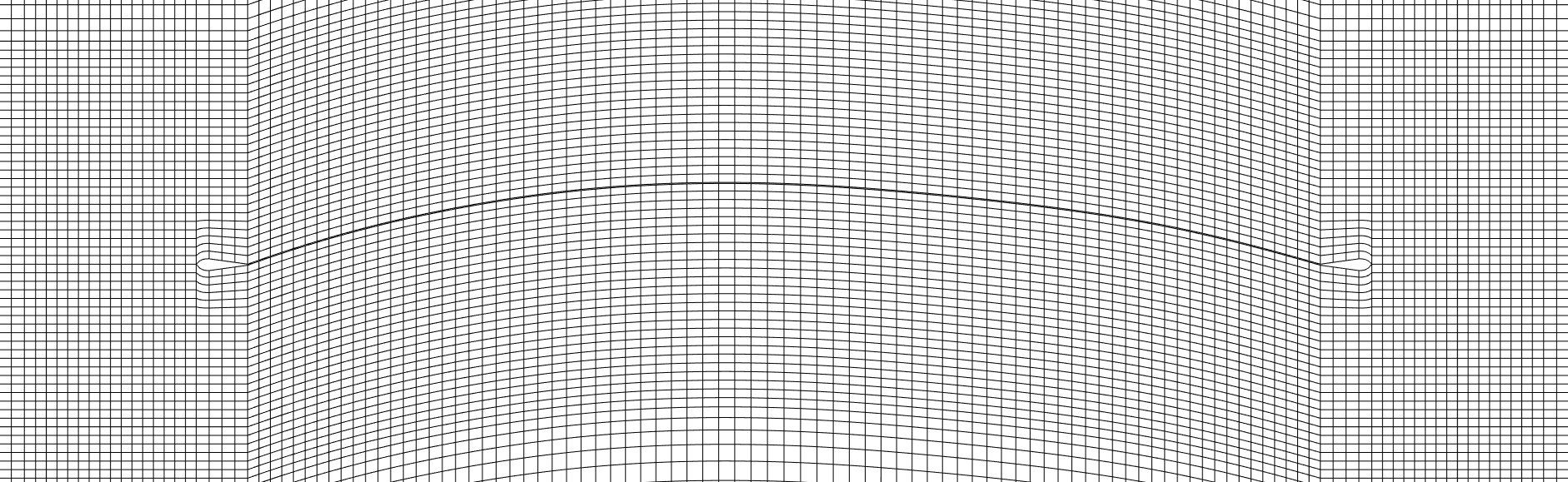}
    \caption{Very fine grid}
    \label{fig:grid_veryfine}
  \end{subfigure}
  \caption{Grid sequence}
  \label{fig:grid}
\end{figure}

A set of computations are carried out using five different grids, shown in Fig.~(\ref{fig:grid}), whose characteristics are provided in Tab.~(\ref{tab:grid}). For all cases, a cubic basis is employed. Note that the coarsest grid counts only two elements and five DOFs to represent the membrane deformation, but the membrane remains smooth thanks to the $C^2$ regularity of the basis.

\begin{table}[!ht]
    \centering
    \begin{tabular}{|c|c|c|}
        \hline
        & fluid DOF number & structure DOF number\\
        \hline
        very coarse & 33,792 & 5 \\
        coarse & 97,280 & 10 \\
        medium & 179,200 & 20 \\
        fine & 270,976 & 40 \\
        very fine & 479,232 & 80 \\
        \hline
    \end{tabular}
    \caption{Size of the grids}
    \label{tab:grid}
\end{table}

\begin{figure}[!ht]
  \begin{subfigure}[b]{0.49\linewidth}
    \centering
    \includegraphics[width=\textwidth]{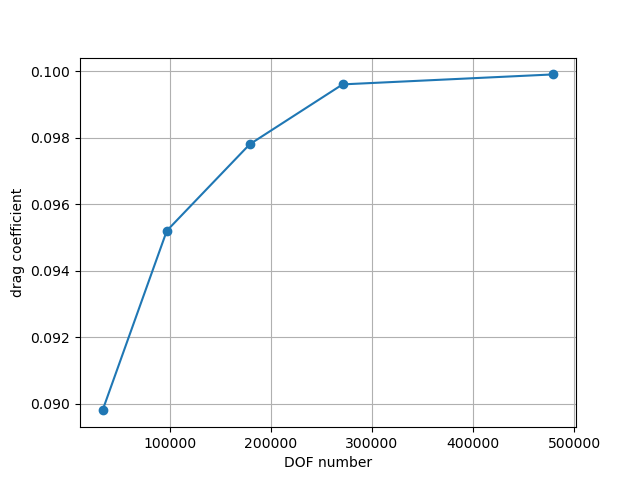}
    \caption{Drag coefficient}
    \label{fig:cd_conv}
  \end{subfigure}
  \begin{subfigure}[b]{0.49\linewidth}
    \centering
    \includegraphics[width=\textwidth]{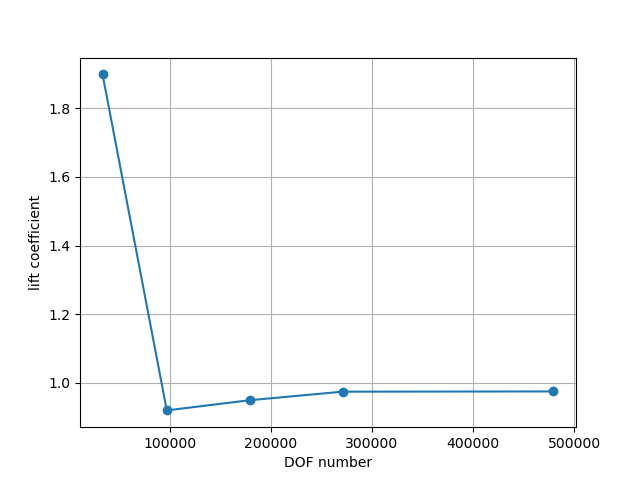}
    \caption{lift coefficient}
    \label{fig:cl_conv}
  \end{subfigure}
  \caption{Grid convergence study.}
  \label{fig:grid_conv}
\end{figure}

The convergence assessment is performed for an incidence of value $\alpha=8^\circ$, for which a periodic solution is found, regarding the time-averaged lift and drag coefficients. The time averaging process is carried out after the transcient solution has vanished. The results are shown in Fig.~(\ref{fig:grid_conv}). As can be seen, the fine grid provides almost converged values. The difference between the coefficient values found using the fine grid and the very fine grid is $0.3\%$ for the drag coefficient and $0.07\%$ for the lift coefficient. Consequently, all forthcoming results are provided for the fine grid.

\subsection{Validation of the coupling}

As explained in the method description, only a loosely coupling approach is used in this study. In particular, the load on the structure is different from the forces exerted by the fluid at the end of a time step. As a result, some energy is lost in the coupling during each time step. Therefore, it is mandatory to quantify this loss to assess the results. In this perspective, the energy loss at the fluid-structure interface during each time step is computed by:
\begin{equation}
\mathcal{L}_E = \left | \int_{\Delta t} \int_{\hat{\Gamma}^s} \Delta p \dot{y^s} |\mbox{J}_{\Gamma^s}| \mbox{d}\xi^s -  \int_{\Delta t} \int_{\hat{\Gamma}^f} p \mathbf{v}^f \cdot \mathbf{n}_y |\mbox{J}_{\Gamma^f}| \mbox{d}\xi^f \right | .
\end{equation}

\begin{figure}[!ht]
    \centering
    \includegraphics[width=0.6\textwidth]{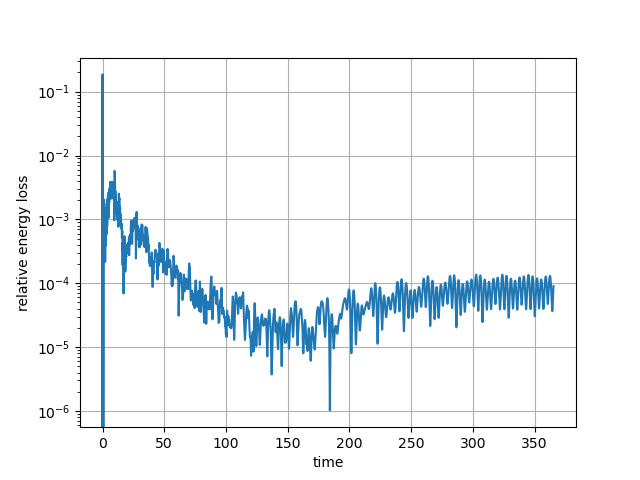}
  \caption{Evolution of relative energy loss along time}
  \label{fig:energy_loss}
\end{figure}

The time evolution of the relative energy loss $\mathcal{L}/\bar{\mathcal{E}}$, where $\bar{\mathcal{E}}$ is the time-averaged value of the energy transfer, is shown in Fig.~(\ref{fig:energy_loss}). As can be seen, the energy loss is large at the beginning of the computation, due to the fast displacements observed at these times, but becomes then very small in comparison with the magnitude of the energy transfered during the coupling. This is an \textit{a posteriori} validation of the coupling strategy.

\subsection{Results}

A study of the membrane wing behavior is then achieved for a set of different incidence values, ranging from $\alpha = 4^\circ$ to $\alpha = 20^\circ$, including some comparisons with results found in the literature. More precisely, the time-averaged lift and drag coefficients are compared to the values obtained in the work of Gordnier~\cite{Gordnier_09}, Li \textit{et al.}~\cite{Li_etal_20} and Sun \textit{et al.}~\cite{Sun_etal_17}. Then, the time-averaged distribution of the membrane deviation and pressure coefficient are compared to those obtained by Serrano-Galiano~\cite{Serrano_16}.

\begin{figure}[!ht]
  \begin{subfigure}[b]{0.49\linewidth}
    \centering
    \includegraphics[width=\textwidth]{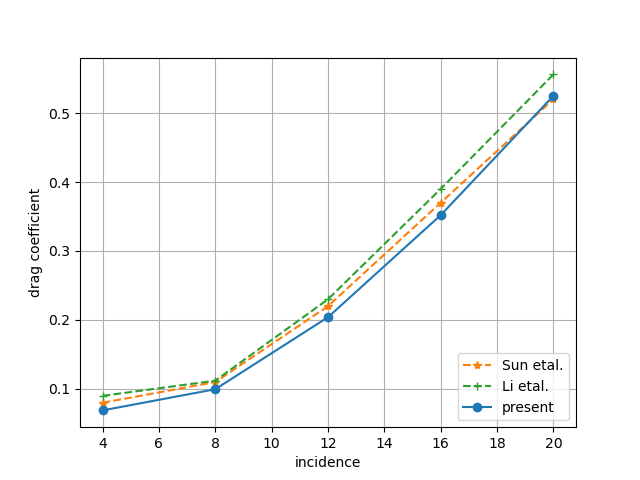}
    \caption{Drag coefficient}
    \label{fig:cd_conv2}
  \end{subfigure}
  \begin{subfigure}[b]{0.49\linewidth}
    \centering
    \includegraphics[width=\textwidth]{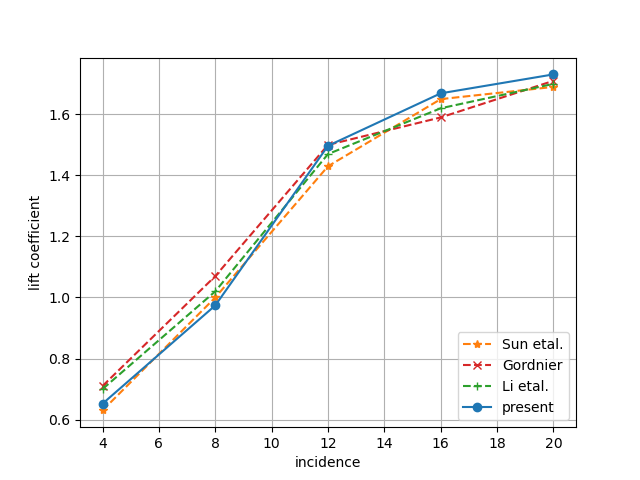}
    \caption{lift coefficient}
    \label{fig:cl_conv2}
  \end{subfigure}
  \caption{Aerodynamic coefficients w.r.t incidence}
  \label{fig:coef_inc}
\end{figure}

The time-averaged lift and drag coefficients with respect to the incidence angle are shown in Fig.~(\ref{fig:coef_inc}). The results obtained look similar to those found in the literature. In particular, one can notice the strong increase of the lift coefficient between $4^\circ$ and $12^\circ$ whereas the increase is lower at higher incidences. As explained below, this corresponds to a change of solution regime. The discrepancies in the aerodynamic coefficients may be explained by the use of different models and numerical configurations. In particular, in~\cite{Li_etal_20} the flow is assumed to be fully incompressible, whereas the computations in \cite{Gordnier_09} and \cite{Serrano_16} are based on compressible flows with Mach number 0.05 and 0.2 respectively.

\begin{figure}[!ht]
  \begin{subfigure}[b]{0.49\linewidth}
    \centering
    \includegraphics[width=\textwidth]{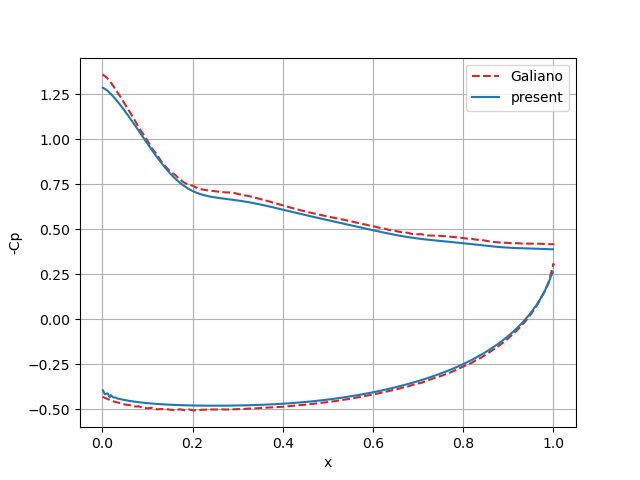}
    \caption{Time-averaged pressure coefficient}
    \label{fig:mean_cp}
  \end{subfigure}
  \begin{subfigure}[b]{0.49\linewidth}
    \centering
    \includegraphics[width=\textwidth]{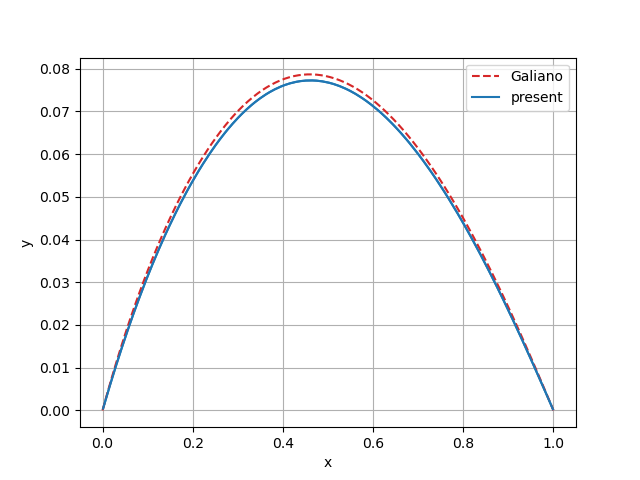}
    \caption{Time-averaged deformation}
    \label{fig:mean_y}
  \end{subfigure}
  \caption{Comparison of membrane pressure and displacement at 8 degrees}
  \label{fig:comp_displ}
\end{figure}

Then, a comparison of the time-averaged pressure coefficient and displacement of the membrane is achieved in Fig.~(\ref{fig:comp_displ}). It corresponds to an incidence of value $\alpha=8^\circ$, for which a periodic solution is found and facilitates the computation of the time-averaged quantities. As can be seen, the overall shape of the pressure and displacement are satisfactory. One can notice however a slight shift in the pressure, that explains a small discrepamcy in the displacement at the maximum deviation location.

\begin{figure}[!ht]
    \centering
    \includegraphics[width=0.25\textwidth]{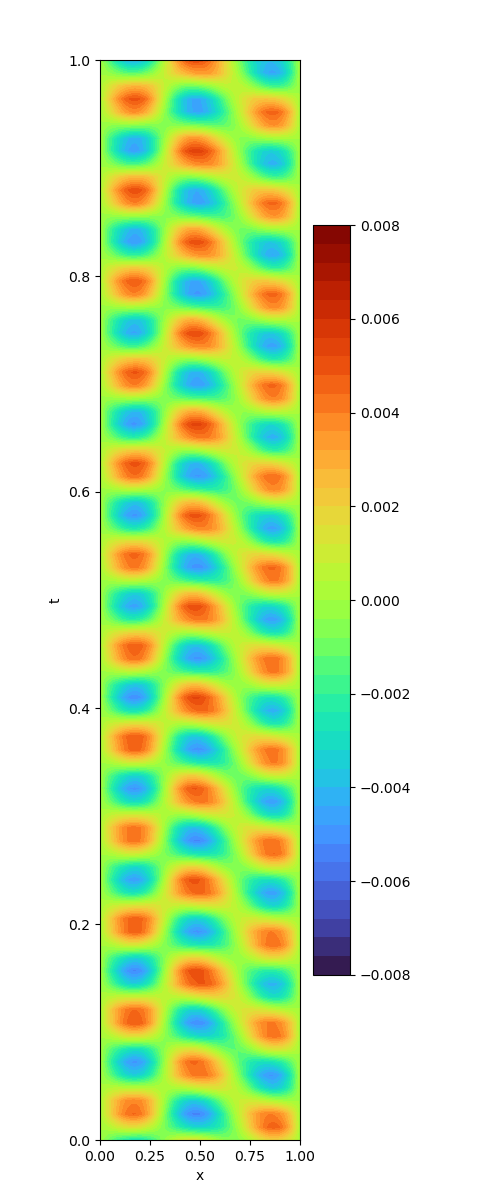}
  \caption{$y(x,t)-\bar{y}(x)$ in $(x,t)$ plane}
  \label{fig:map_displ}
\end{figure}

It is also interesting to compare the membrane displacement dynamics. In this perspective, the difference between the displacement at time $t$ at the location $x$ and the time-averaged displacement $\bar{y}(x)$ is plotted in Fig.~(\ref{fig:map_displ}). Again, we consider here the case for an incidence $8^\circ$, for which the solution is periodic. A good agreement is observed with the work of Serrano-Galiano~\cite{Serrano_16}, which indicates that the membrane dynamics is correctly captured.

\begin{figure}[!ht]
    \centering
    \includegraphics[width=0.8\textwidth]{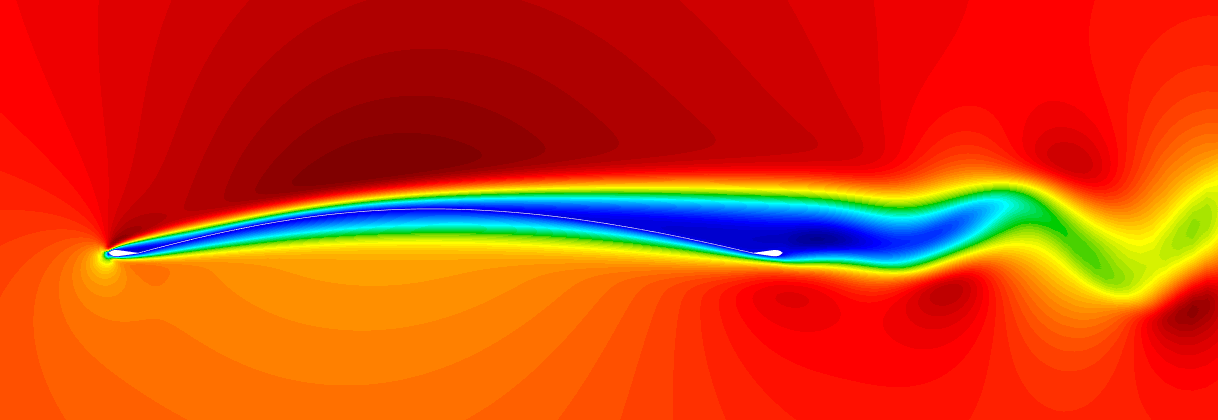}
    \includegraphics[width=0.1\textwidth]{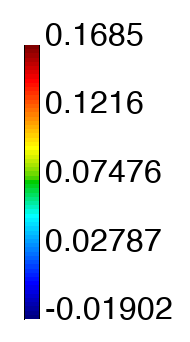}\\
    \includegraphics[width=0.8\textwidth]{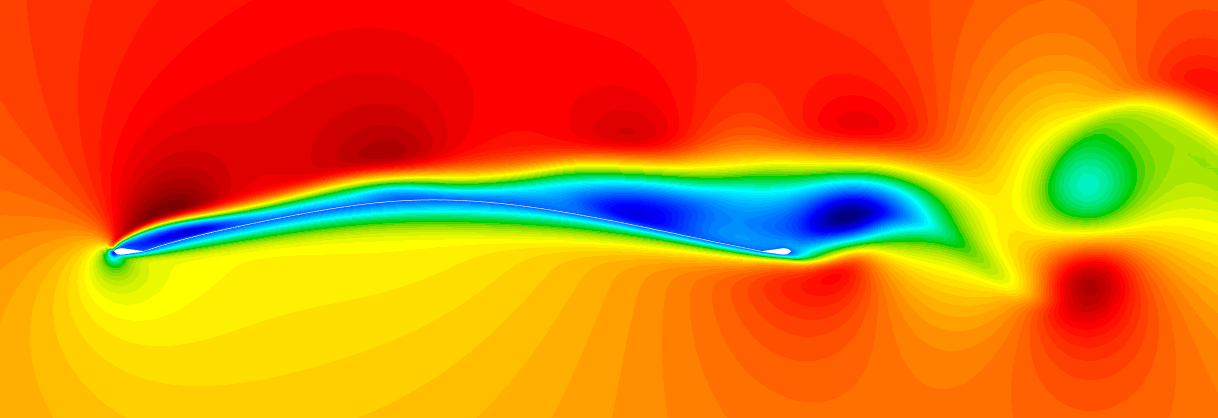}
    \includegraphics[width=0.1\textwidth]{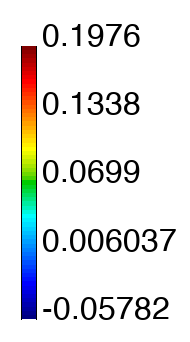}\\
    \includegraphics[width=0.8\textwidth]{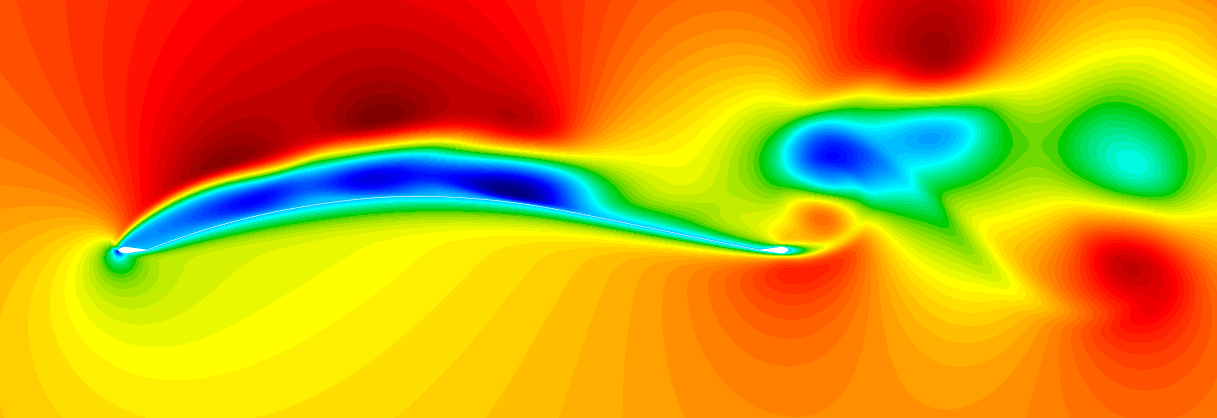}
    \includegraphics[width=0.1\textwidth]{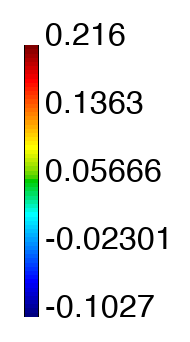}\\
    \includegraphics[width=0.8\textwidth]{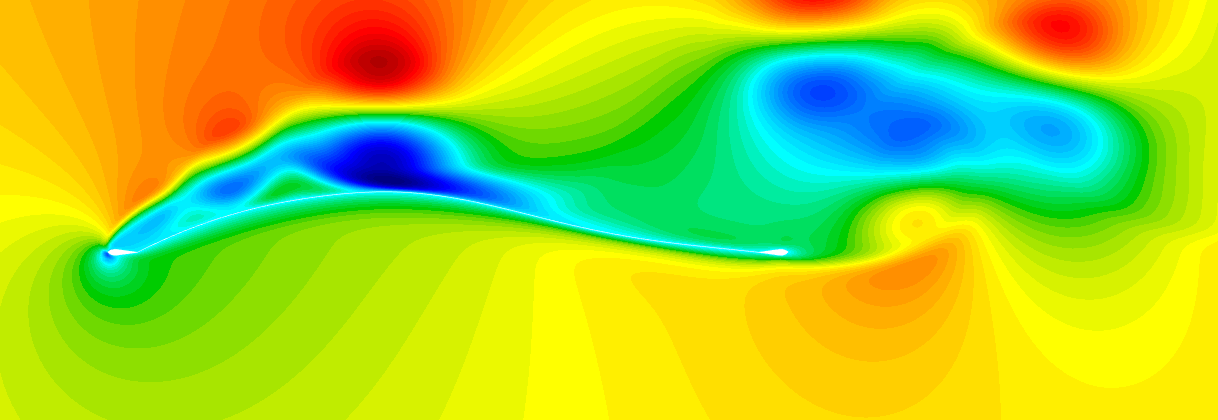}
    \includegraphics[width=0.1\textwidth]{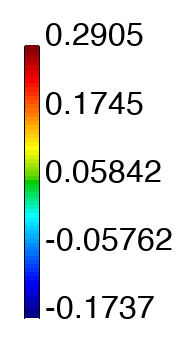}\\
    \includegraphics[width=0.8\textwidth]{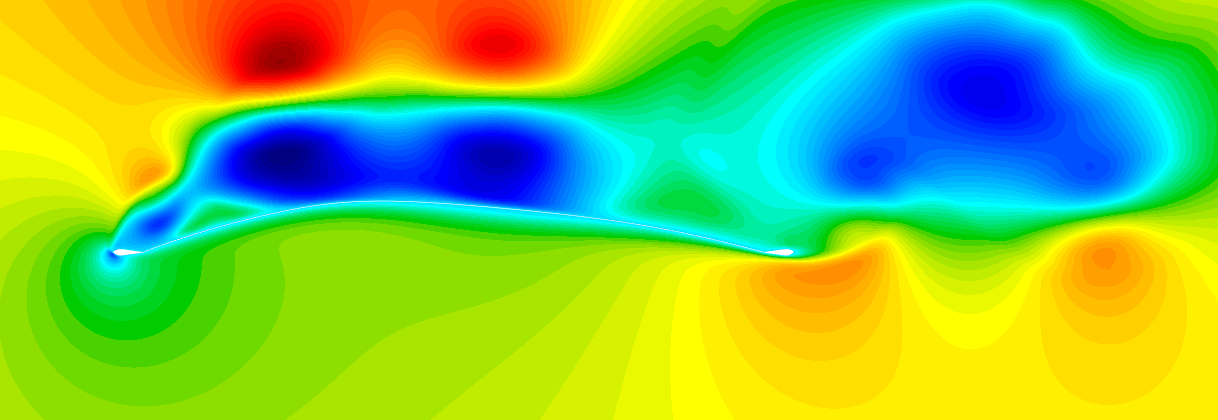}
    \includegraphics[width=0.1\textwidth]{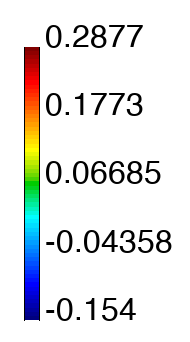}
    \caption{Snapshots of the momentum-x fields at 4, 8, 12, 16 and 20 incidence degrees}
    \label{fig:flow_inc}
\end{figure}

Finally, some snapshots of the velocity field are shown in Fig.~(\ref{fig:flow_inc}) for different incidence values. This allows to understand the dynamics of the fluid-structure interactions and how this evolves with the flow incidence.
\begin{itemize}
    \item $\alpha=4^\circ$: a quasi steady solution is observed for the membrane displacement ; the flow is attached over most of the membrane upper side and a vortex shedding occurs at the trailing edge.
    \item $\alpha=8^\circ$: some vortices are generated at the leading edge and are travelling along the membrane, yielding unsteady membrane displacements ; the flow and the membrane exhibit a periodic behavior. The lift-to-drag ratio is maximal.
    \item $\alpha=12^\circ$: the vortices grow and interact, yielding less regular membrane displacements~; this corresponds to the end of the strong lift increase.
    \item $\alpha=16^\circ$ and higher: the flow and the membrane displacements become chaotic ; the flow exhibits large separation areas and the membrane abrupt moves including possible change of the curvature sign.
\end{itemize}

\clearpage


\section{Turek's case}


\subsection{Computational setup}

We consider then a well-known benchmark for fluid-structure interaction methods, presented in~\cite{Turek_Hron_2006} by Turek and Hron. The geometry is defined by a fixed circular cylinder connected to a flexible bar, located in a channel, as depicted in Fig.~(\ref{fig:turek_config}). The size of the channel is defined by $H=0.41$ and $L=2.5$. The cylinder radius is $r=0.05$ and its center is located at the coordinates $(0.2,0.2)$. The bar has a length $l=0.35$ and a thickness $h=0.02$.

No-slip boundary conditions are prescribed on the walls, including the channel borders. At the left boundary, one imposes the inflow via a parabolic velocity profile:
\begin{equation}
    u_{in}(y) = 1.5 \bar{U} \frac{y(H-y)}{(H/2)^2} .
\end{equation}
As suggested in~\cite{Turek_Hron_2006}, a smooth increase of the flow is imposed to avoid any uncontrolled transcient effects at the beginning of the simulation:
\begin{equation}
    u(0,y,t) = \begin{cases}
      u_{in}(y) \frac{1-\cos(\pi/2 t)}{2} & \text{if $t<2$}\\
      u_{in}(y) & \text{otherwise} .\\
    \end{cases}
\end{equation}
The inlet density has a value $\rho^f_{in} = 1000$ and $\bar{U} = 1$. The kinematic viscosity is set to $\nu^f=10^{-3}$, yielding a Reynolds number of value 100.
At the right boundary, a pressure value $p_{out}=10^5$ is chosen. As the previous case, inlet and outlet boundary conditions are weakly imposed using Riemann invariants.

The density of the elastic beam is $\rho^s=10^3$, the Poisson coefficient and Young modulus of the material are respectively $\nu^s = 0.4$ and $E=1.4 \, 10^6$.

\begin{figure}[!ht]
  \begin{subfigure}[b]{\linewidth}
    \centering
    \includegraphics[width=\textwidth]{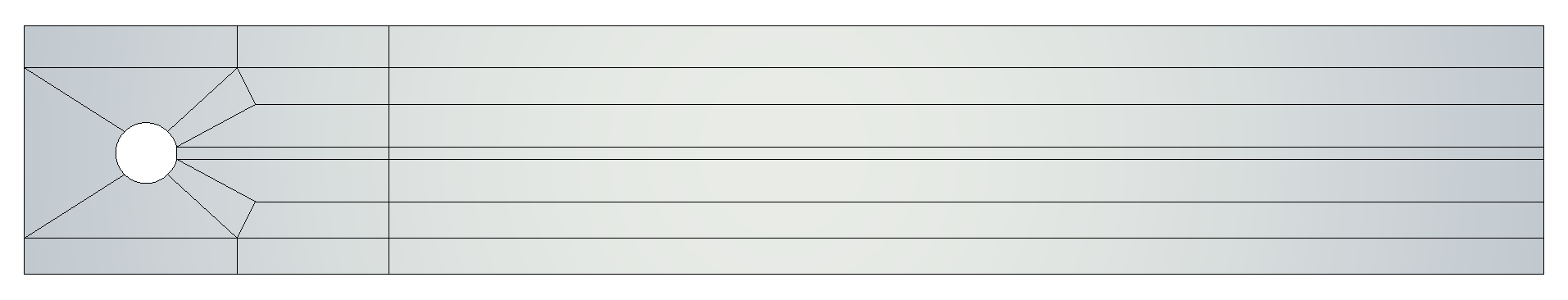}
    \caption{NURBS surfaces that define the fluid and structure domains}
    \label{fig:turek_mesh_ini}
  \end{subfigure}
  \begin{subfigure}[b]{\linewidth}
    \centering
    \includegraphics[width=\textwidth]{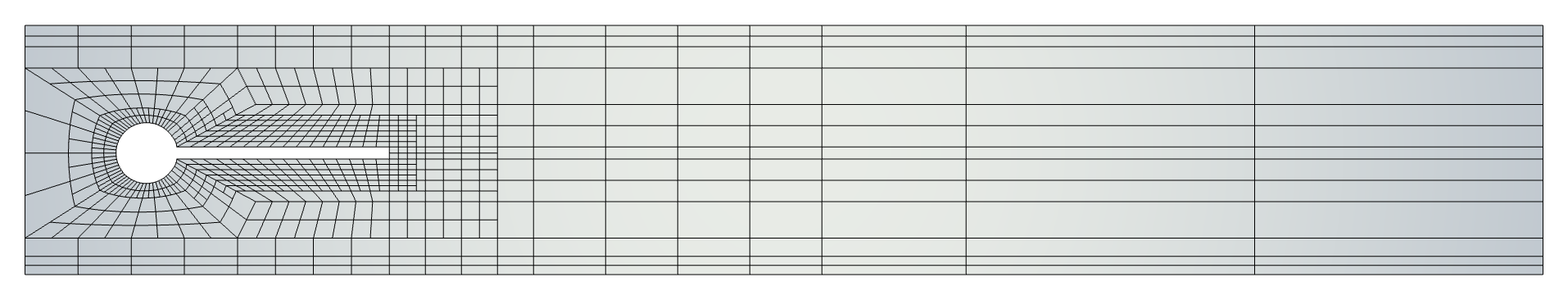}
    \caption{Fluid grid obtained after Bézier extraction and refinement}
    \label{fig:turek_mesh}
  \end{subfigure}
  \begin{subfigure}[b]{\linewidth}
    \centering
    \includegraphics[width=\textwidth]{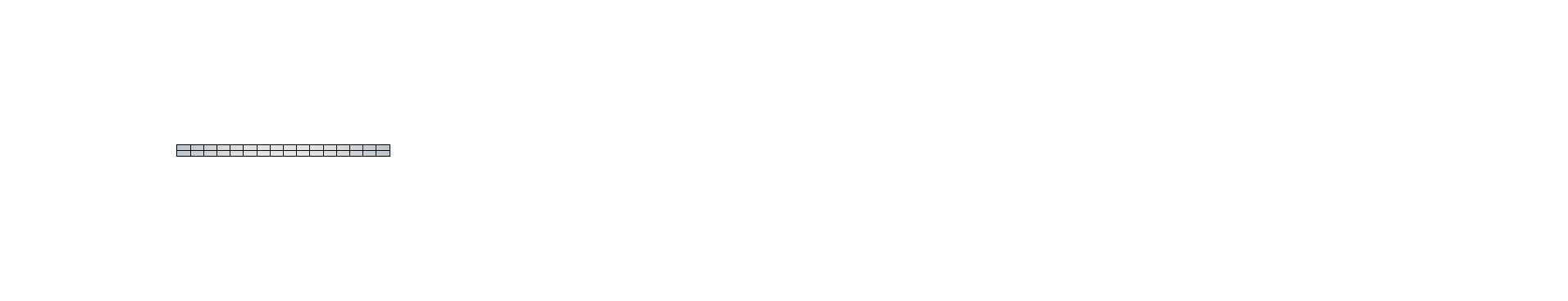}
    \caption{Structure grid obtained after refinement}
    \label{fig:turek_patch}
  \end{subfigure}
  \caption{Fluid and structure grids}
  \label{fig:grid_turek}
\end{figure}

The geometry being more complex than the previous case, the computational domain has be to defined according to a set of connected quadratic NURBS surfaces, shown in Fig.~(\ref{fig:turek_mesh_ini}). Note that non-unitary weights are employed for control points located on the cylinder, in order to represent exactly the geometry. For the fluid part, the NURBS surfaces are then transformed to rational Bézier patches and locally refined to provide the fluid grid, as exposed in Fig~(\ref{fig:turek_mesh}). For the structure part, the grid is simply obtained by refining the NURBS surface corresponding to the bar, see Fig.~(\ref{fig:turek_patch}).

\subsection{Validation of the fluid simulation}

One first validates the simulation of the flow around the cylinder for a fixed bar. The configuration corresponds to the case CFD2 in Turek's benchmark and yields a steady solution. The results in terms of drag and lift coefficients are compared to those obtained by Turek and Hron~\cite{Turek_Hron_2006} and Chabannes~\cite{Chabannes_13}, for grids of similar size for which the values are converged or almost converged. They are summarized in Tab.~(\ref{tab:grid_cfd}).

\begin{table}[!ht]
    \centering
    \begin{tabular}{|c|c|c|c|c|}
        \hline
        & element number & DOF number & $C_D$ & $C_L$\\
        \hline
        \hline
        coarse & 232 & 2,088 & 144.6 & 12.16 \\
        medium & 704 & 6,336 & 138.7 & 10.99 \\
        fine & 2,816 & 25,344 & 138.4 & 10.75 \\
        very fine & 11,264 & 101,376 & 138.4 & 10.66 \\
        \hline
        \hline
        Chabannes & 7,026 & 78,758 & 136.7 & 10.54 \\
        \hline
        Turek-Hron & 9,216 & 177,472 & 136.7 & 10.53 \\
        \hline
    \end{tabular}
    \caption{CFD results}
    \label{tab:grid_cfd}
\end{table}

Four grids of degree two are tested to provide converged results. However, one can notice a discrepancy between the present coefficients and those found in the literature. These results being confirmed by using other types of grid, the discrepancy is attributed to the use of different fluid models. Indeed, the reference values from the literature rely on an incompressible assumption, whereas a compressible flow is considered in this study, using a low Mach number $M=0.085$.


\subsection{Validation of the structure simulation}

A similar validation is carried out for the structural model alone. The deformation of the hyper-elastic beam is computed without the fluid, for a vertical gravity field of value $\mathbf{g}=(0,-2)$. This corresponds to the CSM3 case in Turek's benchmark, for which a periodic solution is obtained. The results are analysed in terms of displacement of the bar extremity (point A in~Fig.(\ref{fig:turek_config})). One provides the mean and amplitude values of the displacement in~Tab.(\ref{tab:grid_csm}) for four different grids with a quadratic basis. Again, a comparison is achieved with the results obtained by Turek and Hron~\cite{Turek_Hron_2006} and Chabannes~\cite{Chabannes_13}. As can be seen, a good agreement is obtained, although a lower number of degrees of freedom is used. This may be a gain resulting from the use of bases of higher regularity, thanks to the isogeometric formulation.

\begin{table}[!ht]
    \centering
    \begin{tabular}{|c|c|c|c|c|}
        \hline
        & element number & DOF number & $u_x$ ($\times 10^{-3}$) & $u_y$ ($\times 10^{-3}$)\\
        \hline
        \hline
        coarse & 32 & 72 & $-14.042 \pm 14.043$ & $-62.520 \pm 64.302$ \\
        medium & 128 & 204 & $-14.273 \pm 14.277$ & $-63.078 \pm 64.876$ \\
        fine & 512 & 660 & $-14.374 \pm 14.379$ & $-63.307 \pm 65.109$ \\
        very fine & 2,048 & 2,340 & $-14.405 \pm 14.411$ & $-63.394 \pm 65.194$ \\
        \hline
        \hline
        Chabannes & 4,199 & 17,536 & $-14.59 \pm 14.59$ & $-63.98 \pm 65.52$ \\
        \hline
        Turek-Hron & 320 & 6,468 & $-14.39 \pm 14.40$ & $-64.27 \pm 64.60$ \\

        \hline

    \end{tabular}
    \caption{CSM results}
    \label{tab:grid_csm}
\end{table}


\subsection{Fluid-structure interaction results}

Finally, the coupled fluid-structure system is simulated, by immersing the cylinder and the elastic bar in the flow. The case corresponds to the FSI2 problem in Turek's benchmark, for which the largest deformations are observed. Computations are performed for three grids with quadratic bases, whose characterstics are detailed in~Tab.(\ref{tab:grid_nb}).

\begin{table}[!ht]
    \centering
    \begin{tabular}{|c|c|c|c|c|}
        \hline
        & fluid elements & fluid DOFs & structure elements & structure DOFs \\
        \hline
        \hline
         coarse & 402 & 3,618 & 8 & 30 \\
        \hline
         medium & 728 & 6,552 & 32 & 72 \\
        \hline
        fine & 2,578 & 23,202 & 128 & 204\\
        \hline
        \hline
        Turek-Hron & 992 & 19,488 & 320 & 6,468 \\
        \hline
        Nordanger \textit{etal.} & 8,020 & 14,862 & 192 & 560 \\
        \hline
    \end{tabular}
    \caption{Grids for FSI case}
    \label{tab:grid_nb}
\end{table}

\begin{figure}[!ht]
  \begin{subfigure}[b]{0.49\linewidth}
    \centering
    \includegraphics[width=\textwidth]{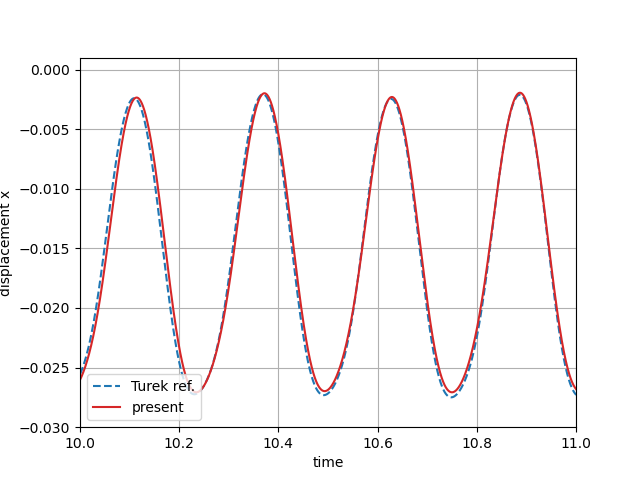}
    \label{fig:displ_x_quad}
  \end{subfigure}
  \begin{subfigure}[b]{0.49\linewidth}
    \centering
    \includegraphics[width=\textwidth]{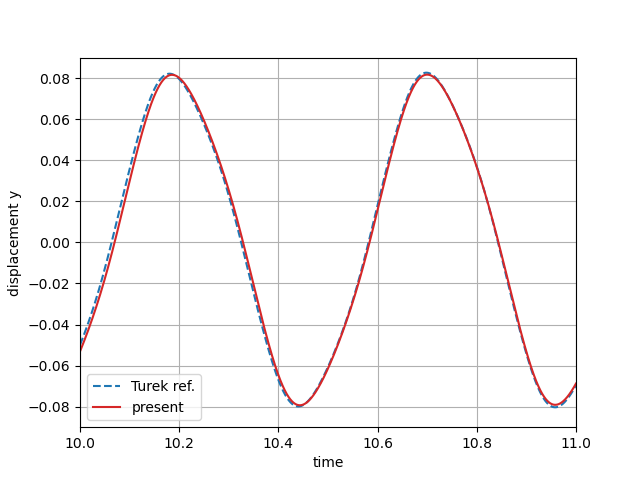}
    \label{fig:displ_y_quad}
  \end{subfigure}
  \caption{Displacement of the bar extremity}
  \label{fig:displ_quad}
\end{figure}

\begin{figure}[!ht]
  \begin{subfigure}[b]{0.49\linewidth}
    \centering
    \includegraphics[width=\textwidth]{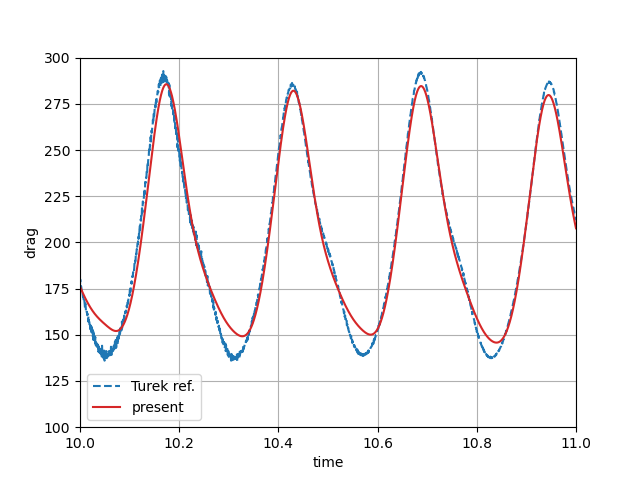}
    \label{fig:turek_fx_quad}
  \end{subfigure}
  \begin{subfigure}[b]{0.49\linewidth}
    \centering
    \includegraphics[width=\textwidth]{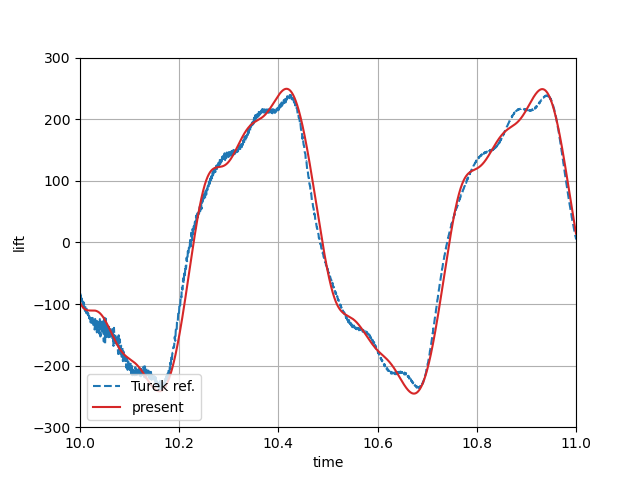}
    \label{fig:turek_fy_quad}
  \end{subfigure}
  \caption{Forces exerted on the cylinder and bar}
  \label{fig:forces_quad}
\end{figure}

The results in terms of drag and lift coefficients, as well as extremity displacements, are reported in~Tab.(\ref{tab:grid_fsi}). The results obtained by Turek and Hron~\cite{Turek_Hron_2006} are considered as reference. The present results are also compared with those obtained by Nordanger \textit{etal.}~\cite{Nordanger_16}, who employed a classical isogeometric finite-element method for the fluid and the structure, with an exact grid matching at the interface. To have a deeper analysis of the dynamics, the time evolution of the displacements and forces are also provided for the fine grid in Fig.~(\ref{fig:displ_quad}) and (\ref{fig:forces_quad}), with comparison with the values obtained by Turek and Hron using their finest discretization.

\begin{table}[!ht]
    \centering
    \begin{tabular}{|c|c|c|c|c|}
        \hline
        &  $u_x$ ($\times 10^{-3}$) & $u_y$ ($\times 10^{-3}$) & $C_D$ & $C_L$\\
        \hline
        \hline
        coarse & $-12.37 \pm 10.88$ & $-0.469 \pm 73.25$ & $214.7 \pm 63.14$ & $2.956 \pm 229.2$\\
        \hline
        medium & $-13.95 \pm 12.15$ & $-1.246 \pm 78.75$ & $217.88 \pm 68.98$ & $1.53 \pm 243.7$\\
        \hline
        fine & $-14.55 \pm 12.55$ & $-1.267 \pm 80.63$ & $215.9 \pm 68.78$ & $3.00 \pm 247.7$\\
        \hline
        \hline
        Turek-Hron & $-14.02 \pm 12.03$ & $-1.250 \pm 79.30$ & $210.1 \pm 72.62$ & $0.25 \pm 227.9$\\
        \hline
        Nordanger \textit{etal.} & $-14.75 \pm 12.80$ & $-1.300 \pm 81.70$ & $214.5 \pm 78.7$ & $1.19 \pm 229.7$\\
        \hline

    \end{tabular}
    \caption{FSI results}
    \label{tab:grid_fsi}
\end{table}

Globally, the results presented are in line with the reference values found in the literature. In particular, a very good agreement is obtained for the time-evolution of the displacement of the bar extremity. Larger discrepancies are observed for the drag and lift coefficients, in particular in the extrema values for the drag. This could be related to the use of different flow models, as already mentioned for the fluid simulation alone. Note that some instabilities in the coupling are reported by Turek and Hron, yielding noisy forces that can be seen in~Fig.(\ref{fig:forces_quad}). However, no such instabilities are observed in the present computations.

\begin{figure}[!ht]
  \begin{subfigure}[b]{0.52\linewidth}
    \centering
    \includegraphics[width=0.9\textwidth]{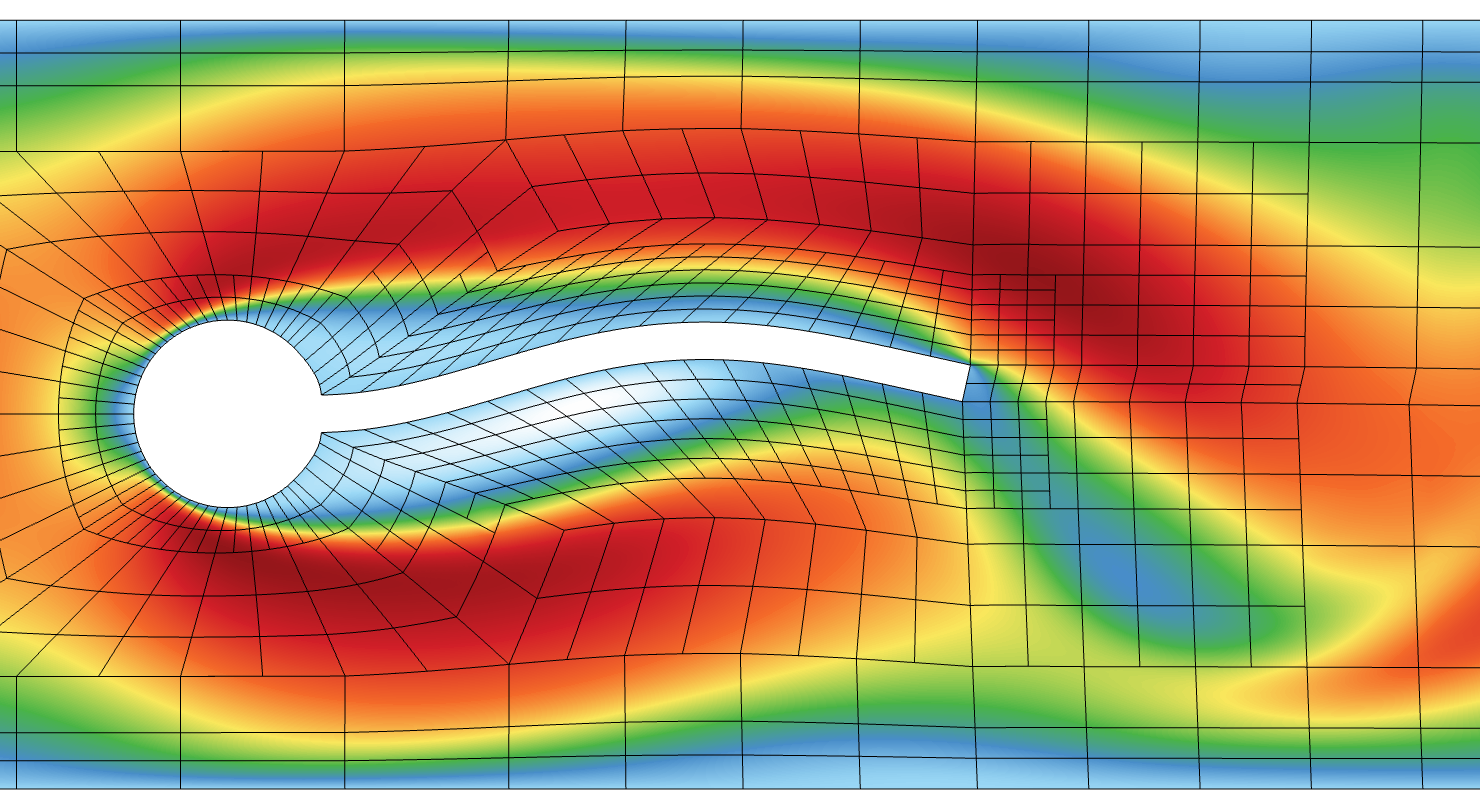}\\
    \includegraphics[width=0.9\textwidth]{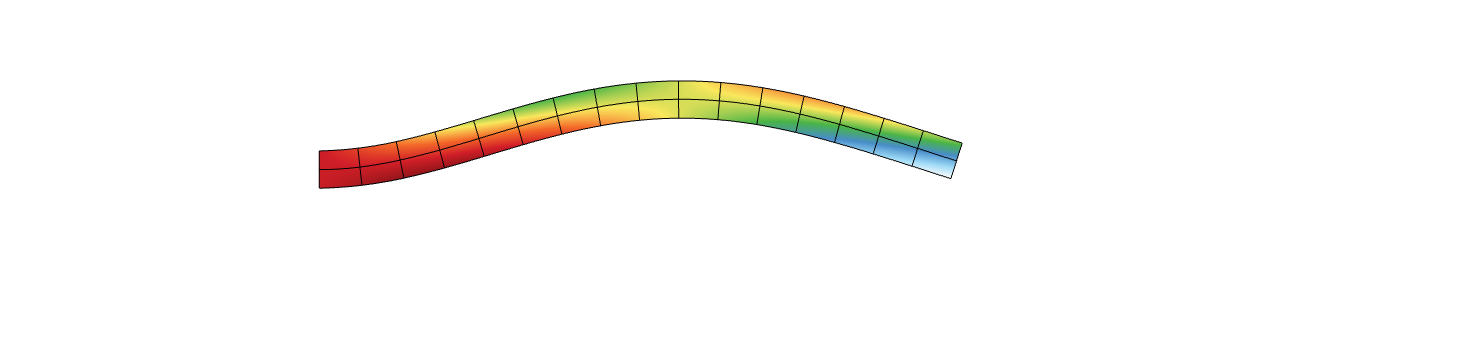}
    \caption{$t = t_0$}
    \label{fig:img500}
\end{subfigure}\hspace{-4mm}
  \begin{subfigure}[b]{0.52\linewidth}
    \centering
    \includegraphics[width=0.9\textwidth]{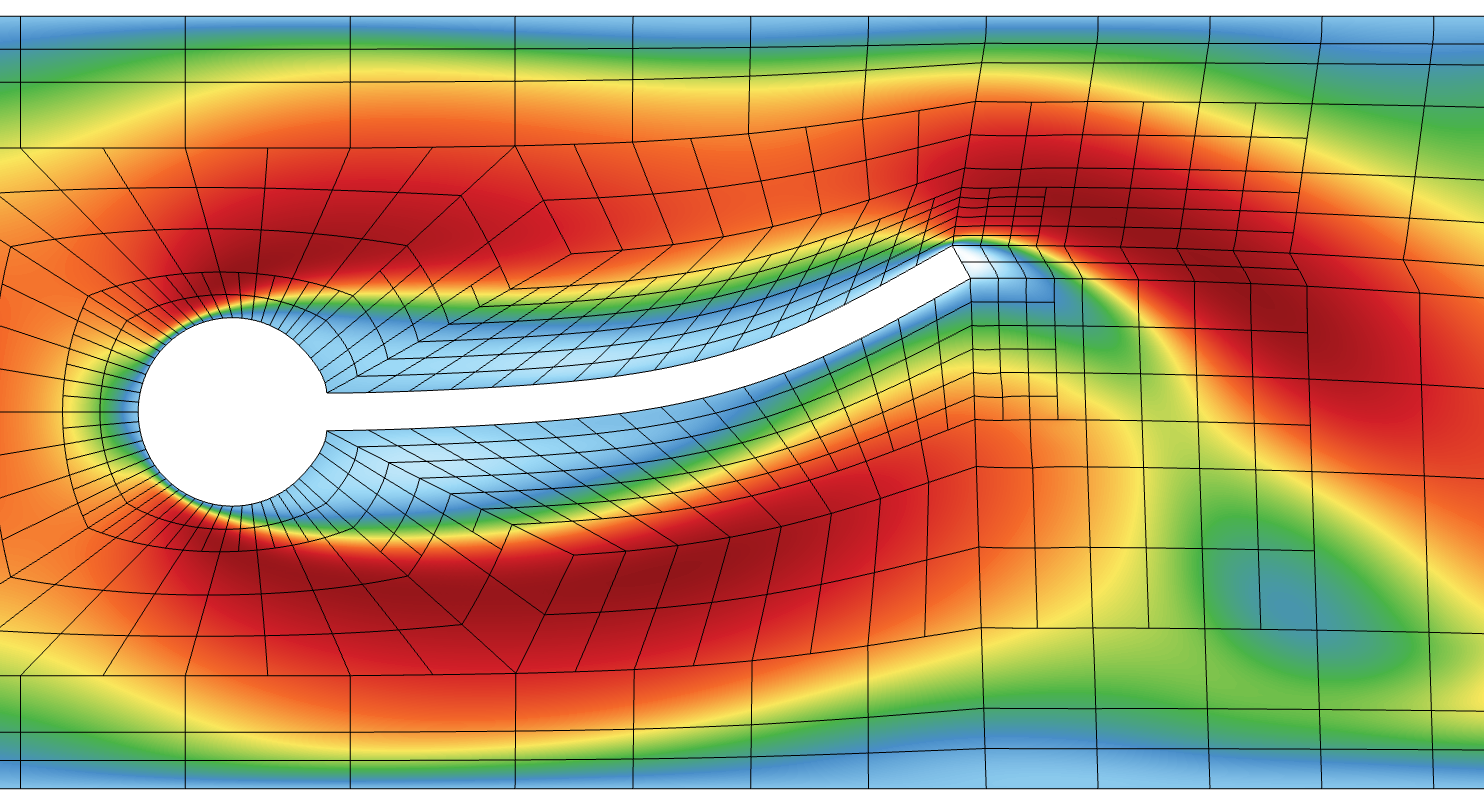}\hspace{-1mm}
    \includegraphics[width=0.09\textwidth]{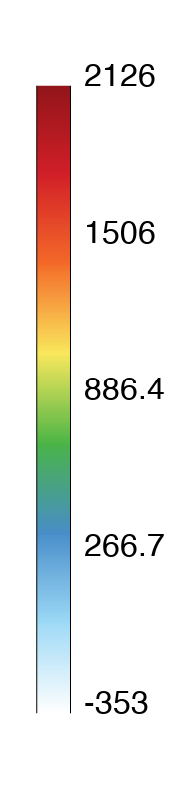}\\
    \includegraphics[width=0.9\textwidth]{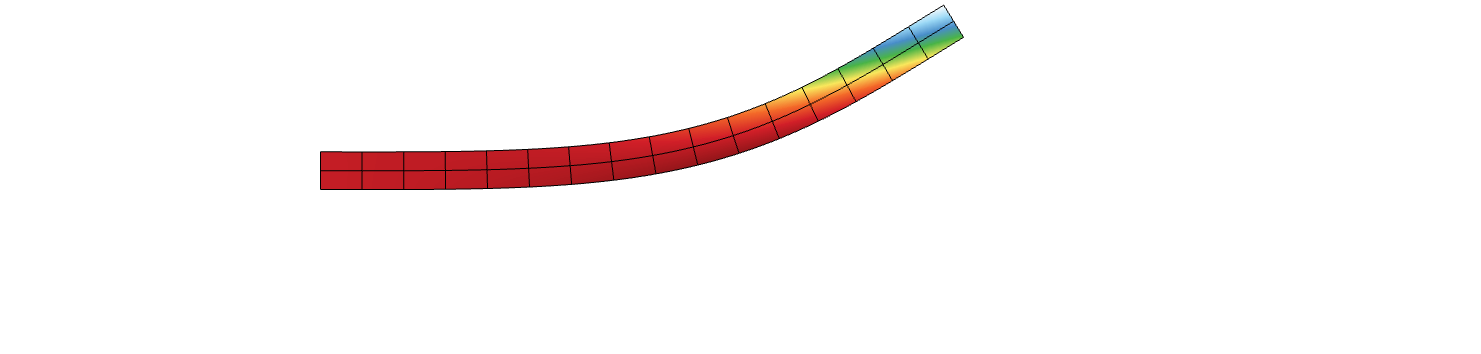}\hspace{-1mm}
    \includegraphics[width=0.09\textwidth]{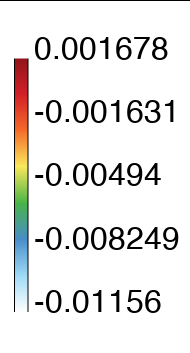}
    \caption{$t = t_0 + T/6$}
    \label{fig:img504}
  \end{subfigure}\\

  \begin{subfigure}[b]{0.52\linewidth}
    \centering
    \includegraphics[width=0.9\textwidth]{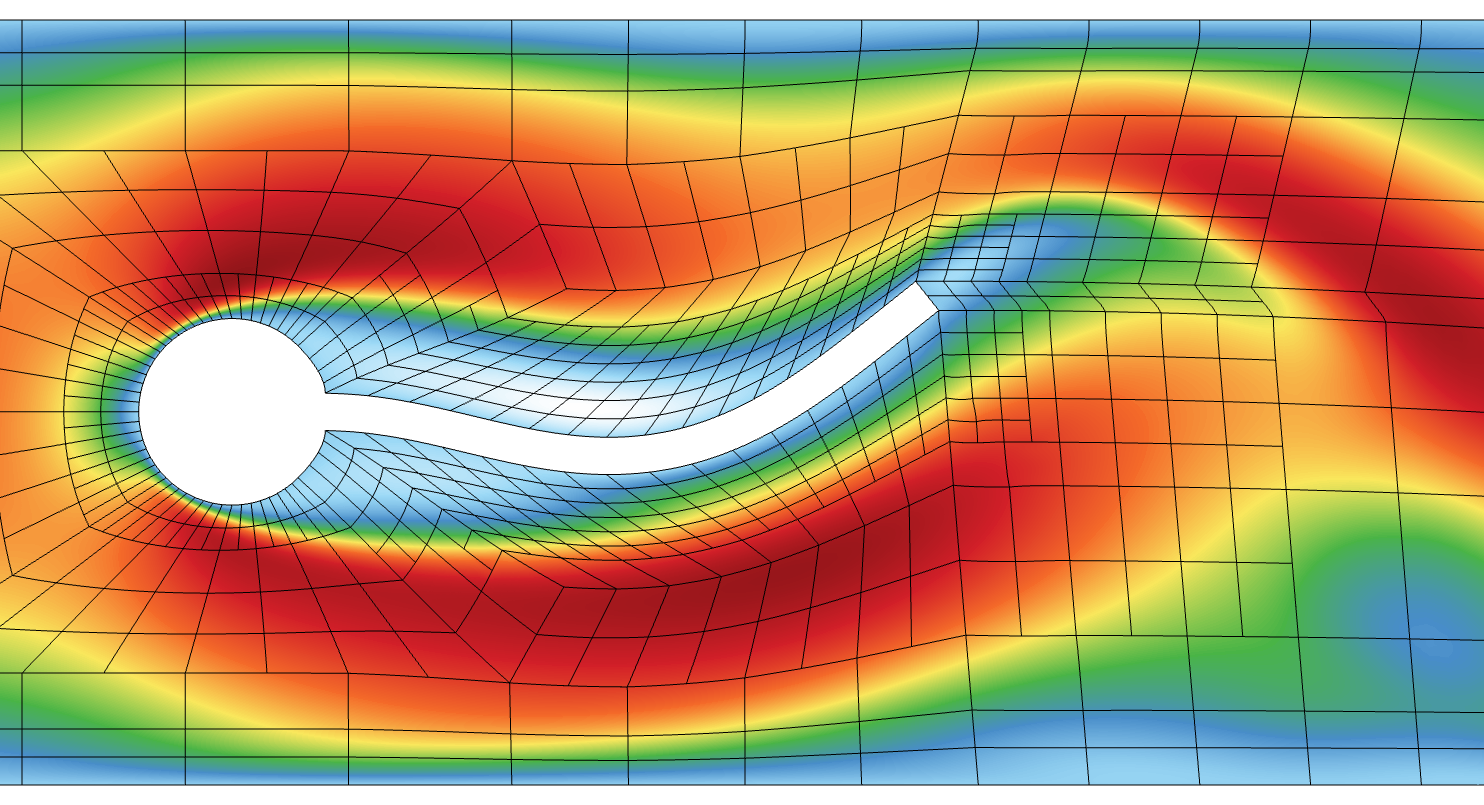}\\
    \includegraphics[width=0.9\textwidth]{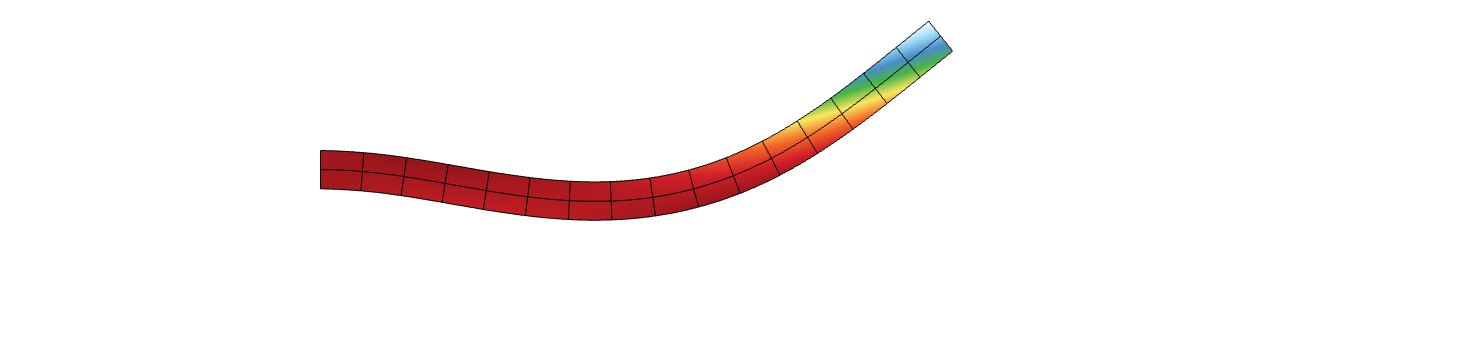}
    \caption{$t = t_0 + T/3$}
    \label{fig:img508}
  \end{subfigure}\hspace{-4mm}
  \begin{subfigure}[b]{0.52\linewidth}
    \centering
    \includegraphics[width=0.9\textwidth]{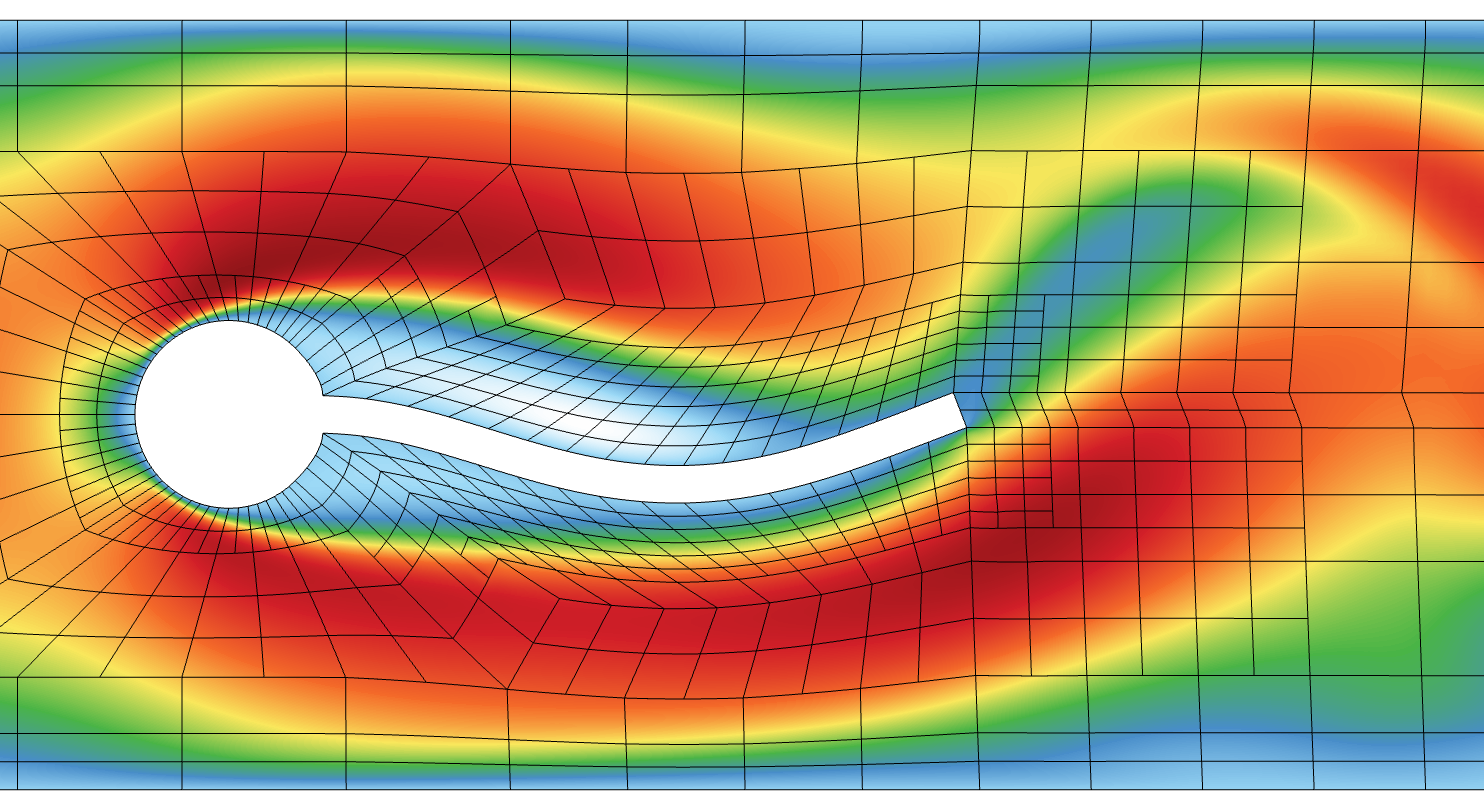}\hspace{-1mm}
    \includegraphics[width=0.09\textwidth]{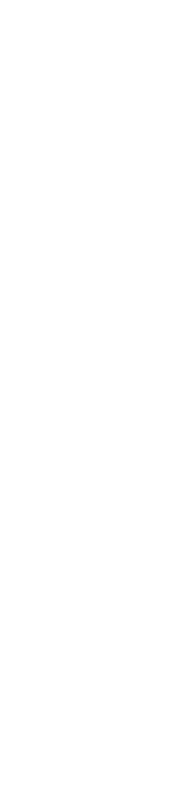}\\
    \includegraphics[width=0.9\textwidth]{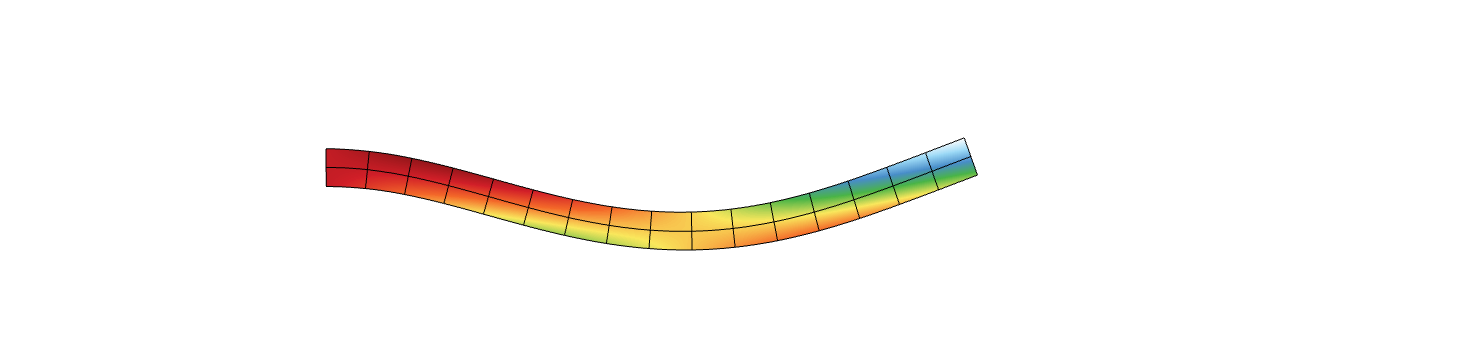}\hspace{-1mm}
    \includegraphics[width=0.09\textwidth]{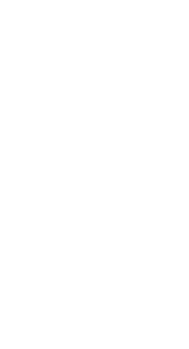}
    \caption{$t = t_0 + T/2$}
    \label{fig:img512}
  \end{subfigure}\\

  \begin{subfigure}[b]{0.52\linewidth}
    \centering
    \includegraphics[width=0.9\textwidth]{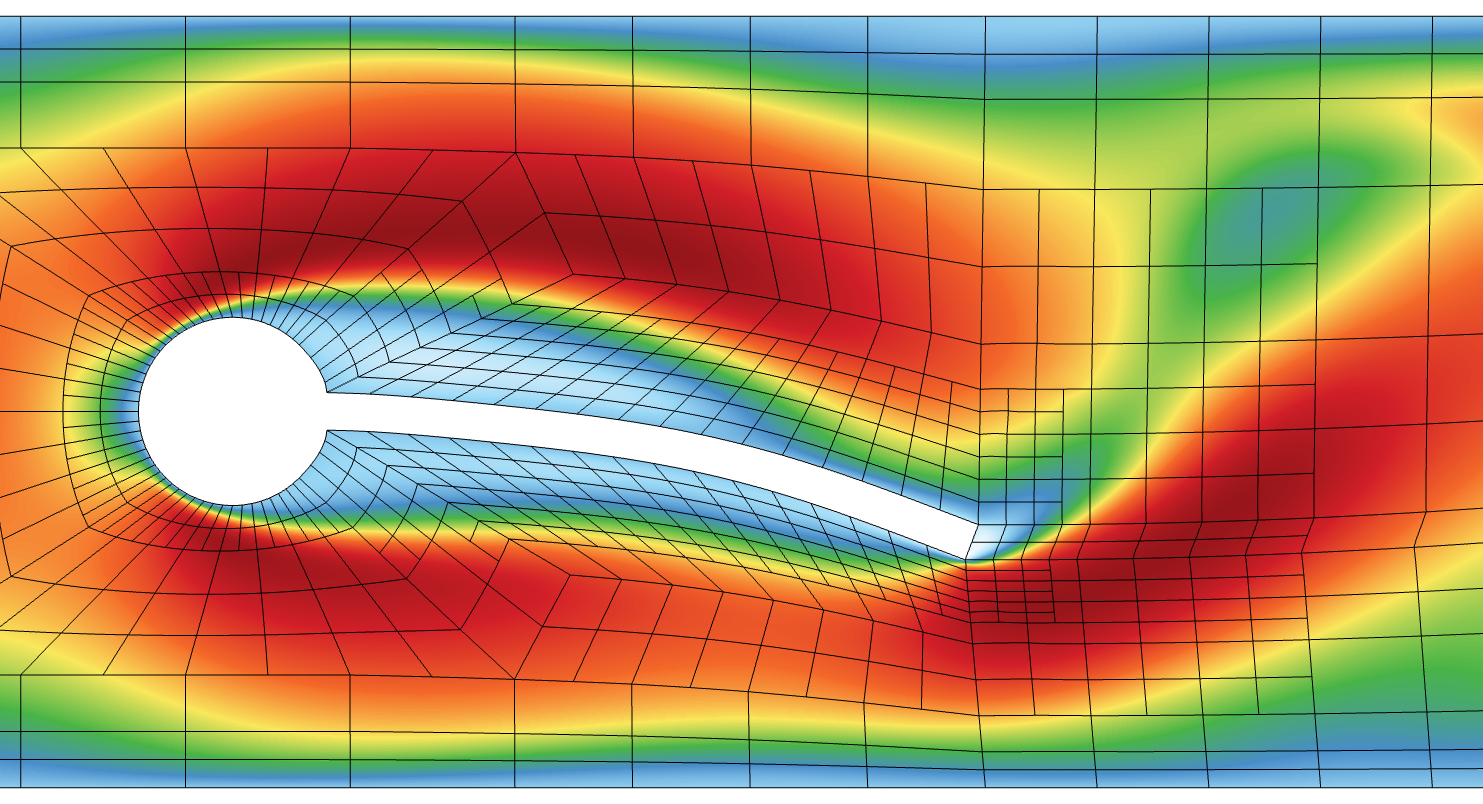}\\
    \includegraphics[width=0.9\textwidth]{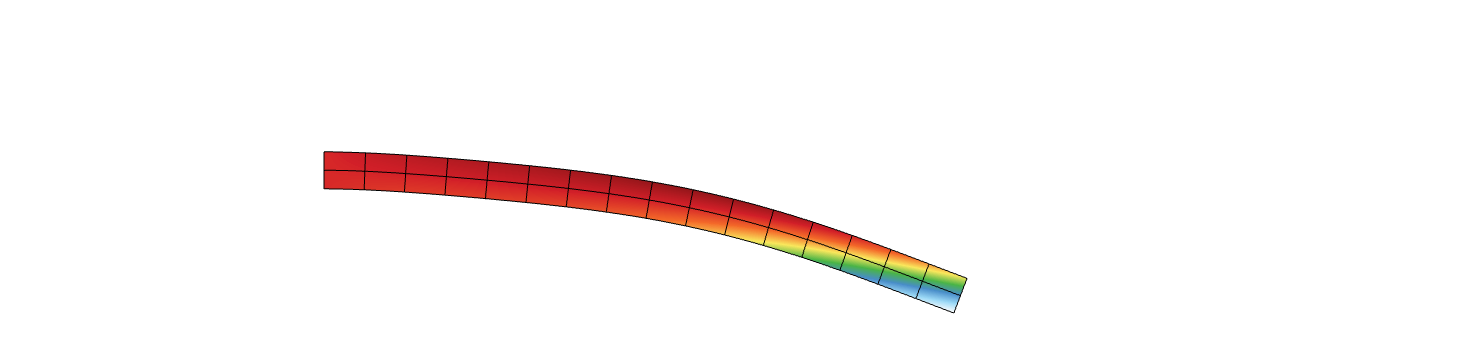}
    \caption{$t = t_0 + 2T/3$}
    \label{fig:img516}
  \end{subfigure}\hspace{-4mm}
  \begin{subfigure}[b]{0.52\linewidth}
    \centering
    \includegraphics[width=0.9\textwidth]{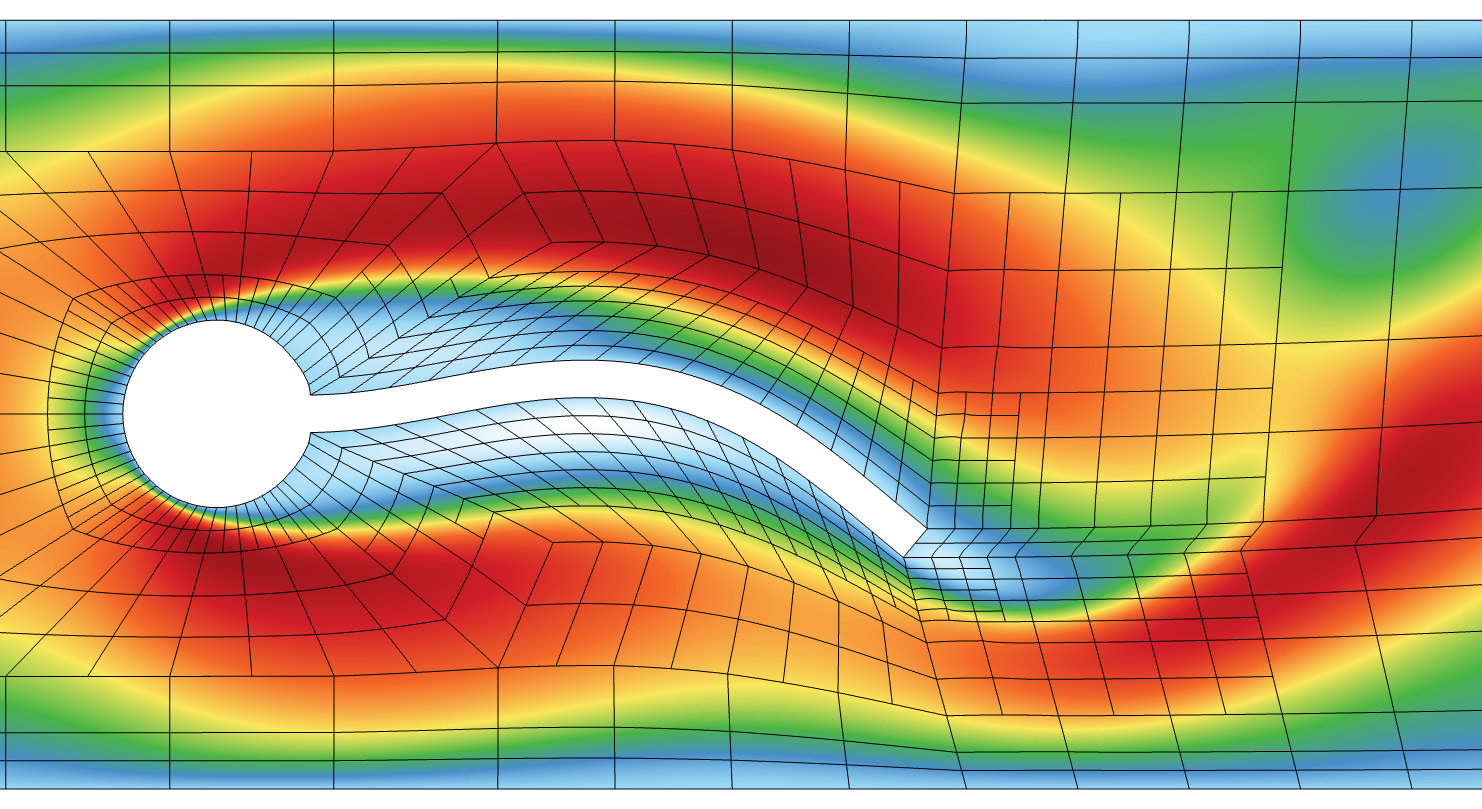}\hspace{-1mm}
    \includegraphics[width=0.09\textwidth]{blanc_turek1}\\
    \includegraphics[width=0.9\textwidth]{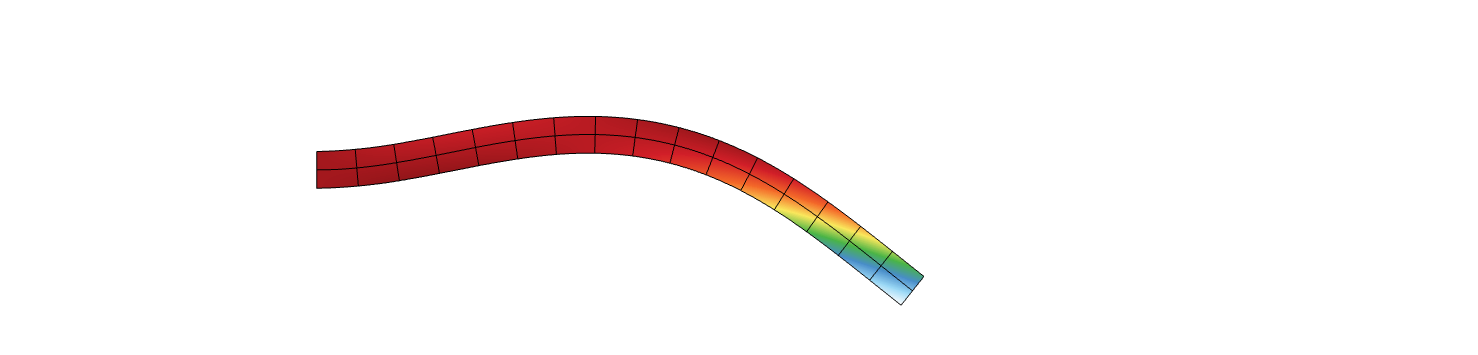}\hspace{-1mm}
    \includegraphics[width=0.09\textwidth]{blanc_turek2}
    \caption{$t = t_0 + 5T/6$}
    \label{fig:img520}
  \end{subfigure}
  \caption{Momentum-x (top) and displacement-x (bottom) fields with grids at some times}
  \label{fig:turek_field}
\end{figure}

Moreover, some snapshots of the velocity and displacement fields (first components) are shown in Fig.~(\ref{fig:turek_field}), including the deformed grids. This case corresponds to the medium grid in~Tab.(\ref{tab:grid_nb}), for six different time steps in a single solution period $T$. One can notice the smoothness of the quadratic solution field, although a discontinuous Galerkin method is employed for the flow solver. The smoothness of the deformation of the fluid grid can also be underlined, despite the large displacement imposed by the structure. This results from the regularity of the underlying NURBS surfaces that control the deformation of the fluid grid, shown in Fig.(\ref{fig:turek_mesh_ini}), although a simple explicit algorithm is actually used for the deformation.

\begin{figure}[!ht]
    \centering
    \includegraphics[width=0.7\textwidth]{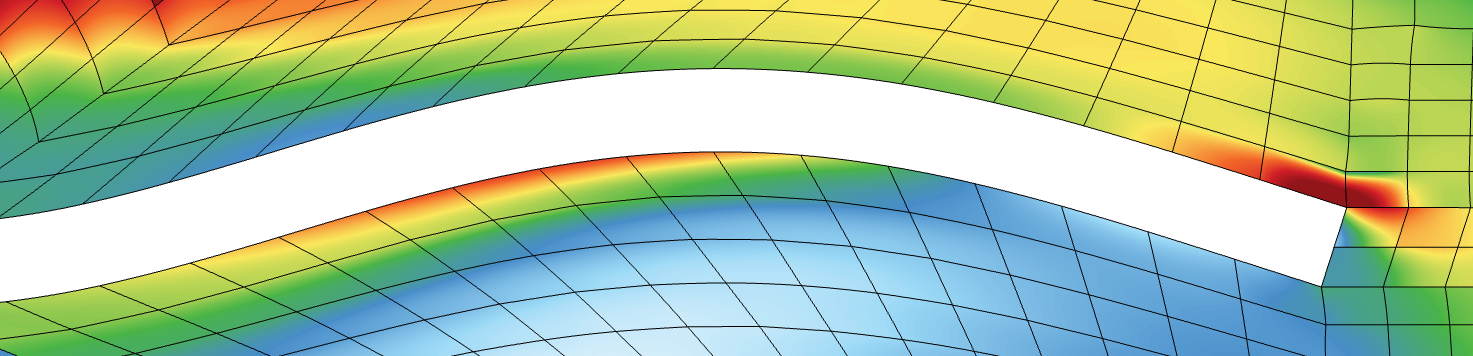}
    \includegraphics[width=0.06\textwidth]{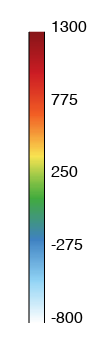}\\
    \includegraphics[width=0.7\textwidth]{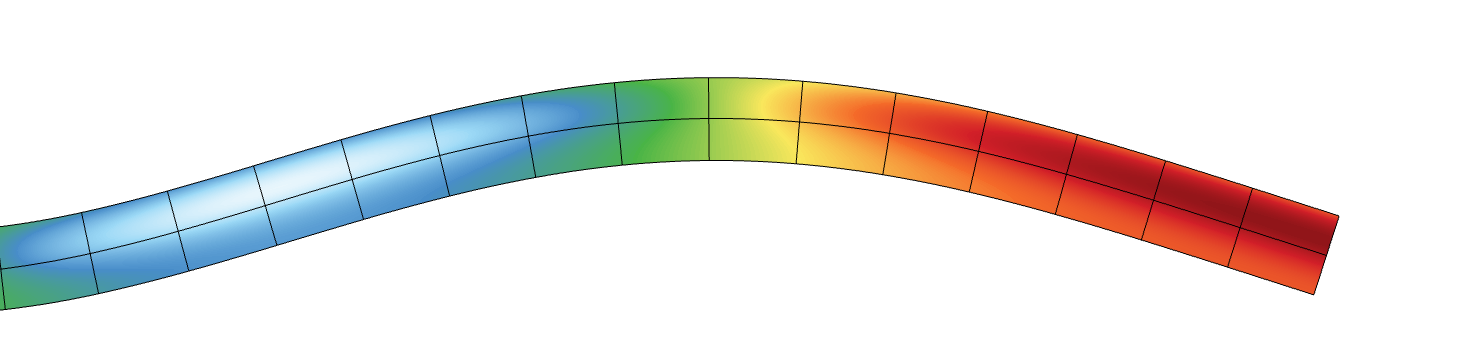}
    \includegraphics[width=0.06\textwidth]{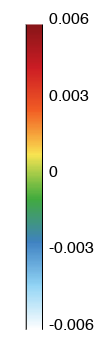}
  \caption{Momentum-x (top) and displacement-x (bottom) derivatives at $t_0$}
  \label{fig:turek_derivatives}
\end{figure}

Finally, it is interesting to examine the derivatives of the solution fields, since these quantities can be critical for several FSI simulations. Fig.~(\ref{fig:turek_derivatives}) shows the derivative of the first component of the momentum and displacement with respect to $\eta$ direction, in the vicinity of the coupling interface, at time $t_0$. One can underline the smoothness of the derivative of the displacement in the bar, which benefits from the $C^1$ continuity of the solution. The derivative of the flow momentum is less regular due to the discontinuities of the solution field at element faces.

\section{Conclusion}

A new approach for the simulation of fluid-structure interactions is presented, relying on NURBS bases. It combines continuous and discontinuous Galerkin formulations for the stucture and the fluid respectively. The use of NURBS bases enables an accurate description of the geometry, consistent with the CAD methodologies. An important property of this representation, namely Bézier extraction, is leveraged to switch between continuous and discontinuous representations at the interface between the fluid and the structure domains with an exact matching, whatever the degree of the basis. Consequently, high-order displacement fields of the structure can be transmitted to the fluid without approximation. The smoothness of this representation can also be exploited to define an explicit but robust grid deformation approach.

This methodology has been demonstrated for two problems, involving compressible viscous flows and non-linear elastic structures. The results provided are in line with reference computations from the literature. Especially, they seem to indicate that the high regularity of the representation is beneficial in terms of convergence, although it is difficult to provide a formal proof.

The extension of this work to more complex three-dimensional geometries raises various issues. In particular, the generation of the computational domain by tessellation of NURBS volumes is obviously a challenging task, which requires strong interactions with CAD community. From algorithmic viewpoint, the proposed coupling scheme can be extended to three-dimensional problems in a rather straightforward way, thanks to the tensorization of the NURBS representation. Nevertheless, implementation of local refinement will necessarily be more complex.

Future works will be focused on fluid-structure interaction problems involving transonic or supersonic flows and light structures, more demanding in terms of conservativity and for which the proposed mixed continuous / discontinuous Galerkin approach seems to be well suited.

\section*{Code Repository}
The developed methodology is implemented in the \textbf{Igloo} software suite, which has been employed to perform all the presented computations. The source code and data are available, under the GNU General Public Licence v3, at the following repository: \url{https://gitlab.inria.fr/igloo/igloo/-/wikis/home}.

\section*{Acknowledgements}
The authors are grateful to the OPAL infrastructure from Université Côte d'Azur for providing resources and support.

\bibliographystyle{plain}
\bibliography{references}

\begin{thebibliography}{10}

\bibitem{Apostolatos_19}
A.~Apostolatos, G.~{De Nayer}, K.-U. Bletzinger, M.~Breuer, and R.~Wüchner.
\newblock Systematic evaluation of the interface description for
  fluid–structure interaction simulations using the isogeometric mortar-based
  mapping.
\newblock {\em Journal of Fluids and Structures}, 86:368--399, 2019.

\bibitem{Bazilevs2006}
Y.~Bazilevs, L.~{Beir{\~{a}}o}~da Veiga, J.~A. Cottrell, T.~J.~R. Hughes, and
  G.~Sangalli.
\newblock Isogeometric analysis: Approximation, stability and error estimates
  for h-refined meshes.
\newblock {\em Mathematical Methods and Models in Applied Sciences},
  16:1031--1090, 2006.

\bibitem{Bazilevs_12}
Y.~Bazilevs, M.-C. Hsu, and M.A. Scott.
\newblock Isogeometric fluid–structure interaction analysis with emphasis on
  non-matching discretizations, and with application to wind turbines.
\newblock {\em Computer Methods in Applied Mechanics and Engineering},
  249-252:28--41, 2012.

\bibitem{Blom_98}
F.J. Blom.
\newblock A monolithical fluid-structure interaction algorithm applied to the
  piston problem.
\newblock {\em Computer Methods in Applied Mechanics and Engineering},
  167(3):369--391, 1998.

\bibitem{Chabannes_13}
V.~Chabannes.
\newblock {\em Vers la simulation num{\'e}rique des {\'e}coulements sanguins}.
\newblock PhD thesis, Universit\'e de Grenoble, 2013.

\bibitem{Cockburn_Shu_98}
B.~Cockburn and C.-W. Shu.
\newblock The local discontinuous {G}alerkin method for time-dependent
  convection-diffusion systems.
\newblock {\em SIAM Journal of Num. An.}, 35(6):2440--2463, 1998.

\bibitem{Cottrell_etal_09}
J.A. Cottrell, T.J.R. Hughes, and Y.~Bazilevs.
\newblock {\em Isogeometric analysis : towards integration of {CAD} and {FEA}}.
\newblock John Wiley \& sons, 2009.

\bibitem{Cottrell_etal_07}
J.A. Cottrell, T.J.R. Hughes, and A.~Reali.
\newblock Studies of refinement and continuity in isogeometric analysis.
\newblock {\em Computer Methods in Applied Mechanics and Engineering},
  (196):4160--4183, 2007.

\bibitem{DeBoer_07}
A.~{de Boer}, A.H. {van Zuijlen}, and H.~Bijl.
\newblock Review of coupling methods for non-matching meshes.
\newblock {\em Computer Methods in Applied Mechanics and Engineering},
  196(8):1515--1525, 2007.

\bibitem{DeBoor_78}
C.~De~Boor.
\newblock {\em A Practical Guide to Splines}.
\newblock Springer Verlag, 1978.

\bibitem{Degroote_13}
J.~Degroote.
\newblock Partitioned simulation of fluid-structure interaction.
\newblock {\em Archives of Computational Methods in Engineering},
  (20):185--238, 2013.

\bibitem{Duvigneau_18}
R.~Duvigneau.
\newblock Isogeometric analysis for compressible flows using a {D}iscontinuous
  {G}alerkin method.
\newblock {\em Computer Methods in Applied Mechanics and Engineering},
  333(443--461), 2018.

\bibitem{Duvigneau_20}
R.~Duvigneau.
\newblock {CAD}‐consistent adaptive refinement using a {NURBS}‐based
  discontinuous {G}alerkin method.
\newblock {\em Int. J. for Numerical Methods in Fluids}, February 2020.

\bibitem{Farhat_etal_98}
C.~Farhat, M.~Lesoinne, and P.~Le Tallec.
\newblock Load and motion transfer algorithms for fluid/structure interaction
  problems with non-matching discrete interfaces: momentum and energy
  conservation, optimal discretization and application to aeroelasticity.
\newblock {\em Computer Methods in Applied Mechanics and Engineering},
  157(1-2):95--114, 1998.

\bibitem{Farin_89}
G.~Farin.
\newblock {\em Curves and Surfaces for Computer-Aided Geometric Design}.
\newblock Academic Press, 1989.

\bibitem{Felippa_eta_01}
C.~Felippa, K.~Park, and C.~Farhat.
\newblock Partitioned analysis of coupled mechanical systems.
\newblock {\em Computer Methods in Applied Mechanics and Engineering},
  190(24-25):3247--3270, 2001.

\bibitem{Frohle_persson_14}
Bradley Froehle and Per-Olof Persson.
\newblock A high-order discontinuous galerkin method for fluid--structure
  interaction with efficient implicit--explicit time stepping.
\newblock {\em Journal of Computational Physics}, 272:455--470, 2014.

\bibitem{Gordnier_09}
R.E. Gordnier.
\newblock High fidelity computational simulation of a membrane wing airfoil.
\newblock {\em Journal of Fluids and Structures}, 25(5):897--917, 2009.

\bibitem{Hughes_etal_05}
T.J.R. Hughes, J.A. Cottrell, and Y.~Bazilevs.
\newblock Isogeometric analysis: {CAD}, finite elements, {NURBS}, exact
  geometry, and mesh refinement.
\newblock {\em Computer Methods in Applied Mechanics and Engineering},
  (194):4135--4195, 2005.

\bibitem{Kamensky_15}
D.~Kamensky, M.C. Hsu, D.~Schillinger, J.A. Evans, A.~Aggarwal, Y.~Bazilevs,
  M.S. Sacks, and T.J.R. Hughes.
\newblock An immersogeometric variational framework for fluid–structure
  interaction: Application to bioprosthetic heart valves.
\newblock {\em Computer Methods in Applied Mechanics and Engineering},
  284:1005--1053, 2015.

\bibitem{Li_etal_20}
G.~Li, B.~C. Khoo, and R.~K. Jaiman.
\newblock Computational aeroelasticity of flexible membrane wings at moderate
  reynolds numbers.
\newblock In {\em AIAA Scitech 2020 Forum, Orlando, FL}, 2020.

\bibitem{Michler_etal04}
C.~Michler, S.J. Hulshoff, E.H. {van Brummelen}, and R.~{de Borst}.
\newblock A monolithic approach to fluid–structure interaction.
\newblock {\em Computers \& Fluids}, 33(5):839--848, 2004.
\newblock Applied Mathematics for Industrial Flow Problems.

\bibitem{Molki_10}
M.~Molki and K.~Breuer.
\newblock Oscillatory motions of a prestrained compliant membrane caused by
  fluid--membrane interaction.
\newblock {\em Journal of Fluids and Structures}, 26(3):339--358, 2010.

\bibitem{Nordanger_16}
K.~Nordanger, A.~Rasheed, K.~M. Okstad, A.~M. Kvarving, R.~Holdahl, and
  T.~Kvamsdal.
\newblock Numerical benchmarking of fluid--structure interaction: An
  isogeometric finite element approach.
\newblock {\em Ocean Engineering}, 124:324--339, 2016.

\bibitem{Pezzano_Duvigneau_21}
S.~Pezzano and R.~Duvigneau.
\newblock {A NURBS-based Discontinuous Galerkin method for conservation laws
  with high-order moving meshes}.
\newblock {\em {Journal of Computational Physics}}, 434(1), 2021.

\bibitem{Piegl_95}
L.~Piegl and W.~Tiller.
\newblock {\em The {NURBS} book}.
\newblock Springer-Verlag, 1995.

\bibitem{Piperno_Farhat_01}
S.~Piperno and C.~Farhat.
\newblock Partitioned procedures for the transient solution of coupled
  aeroelastic problems. part ii: Energy transfer analysis and three-dimensional
  applications.
\newblock {\em Computer Methods in Applied Mechanics and Engineering},
  190(24-25):3147--3170, 2001.

\bibitem{RAO2011241}
S.S. Rao.
\newblock Chapter 7 - numerical solution of finite element equations.
\newblock In {\em The Finite Element Method in Engineering (Fifth Edition)},
  pages 241--274. Butterworth-Heinemann, 2011.

\bibitem{Rojratsirikul2010}
P.~Rojratsirikul, Z.~Wang, and I.~Gursul.
\newblock {\em Unsteady fluid-structure interactions of membrane airfoils at
  low Reynolds numbers}, pages 297--310.
\newblock Springer Berlin Heidelberg, 2010.

\bibitem{Serrano_16}
S.~Serrano-Galiano.
\newblock {\em Fluid-structure interaction of membrane aerofoils at low
  reynolds numbers}.
\newblock PhD thesis, University of Southampton, 2016.

\bibitem{Smith_Shyy_95}
R.~Smith and W.~Shyy.
\newblock {Computation of unsteady laminar flow over a flexible
  two‐dimensional membrane wing}.
\newblock {\em Physics of Fluids}, 7(9):2175--2184, 1995.

\bibitem{Sun_etal_17}
X.~Sun, XL. Ren, and JZ. Zhang.
\newblock Nonlinear dynamic responses of a perimeter-reinforced membrane wing
  in laminar flows.
\newblock {\em Nonlinear Dynamics}, 88:749--776, 2017.

\bibitem{Tiomkin_21}
S.~Tiomkin and D.~E. Raveh.
\newblock A review of membrane-wing aeroelasticity.
\newblock {\em Progress in Aerospace Sciences}, 126:100738, 2021.

\bibitem{Turek_Hron_2006}
S.~Turek and J.~Hron.
\newblock Proposal for numerical benchmarking of fluid-structure interaction
  between an elastic object and laminar incompressible flow.
\newblock In {\em Fluid-Structure Interaction}, pages 371--385. Springer Berlin
  Heidelberg, 2006.

\end{thebibliography}

\end{document}